\documentclass[11pt,leqno]{article}

\usepackage{amsmath,amsthm}
\usepackage {latexsym}
\usepackage{amssymb}

\newcommand{\les}{\lesssim}
\newcommand{\bea}{\begin{eqnarray}}\newcommand{\eea}{\end{eqnarray}}
\newcommand{\beq}{\begin{equation}}\newcommand{\ee}{\end{equation}}

\newcommand{\eps}{{\varepsilon}}\newcommand{\R}{{\mathbb R}}

\newcommand{\E}{{\cal E}}

\newcommand{\si}{\sigma}\renewcommand{\b}{\beta}

\def\nn{\nonumber}

\newtheorem{theorem}{Theorem}[section]\newtheorem{lemma}[theorem]{Lemma}
\newtheorem{cor}[theorem]{Corollary}\newtheorem{prop}[theorem]{Proposition}
\newtheorem{proposition}[theorem]{Proposition}

\theoremstyle{remark}
\newtheorem{remark}[theorem]{Remark}
\def\ve{\varepsilon}\def\a{\alpha}\def\ga{\gamma}\def\de{\delta}\def\Si{\Sigma}
\def\il{\int\limits}\def\bm{\left( \begin{array}{cc}}
\def\endm{\end{array}\right)}\newcommand{\eq}{\end{equation}}
\def\nab{\nabla}\def\tr{\text{tr}}\def\a{\alpha}\def\b{\beta}
\def\ga{\gamma}\def\de{\delta}\def\Box{\square}\def\pa{\partial}\def\pab{\bar\pa}

\def \rectangle#1#2{\hbox{\vrule\vbox to #2 {\hrule\hbox to #1{\hfil}\vfil\hrule}\vrule}}
\def\sq{\,\,\rectangle{7pt}{7 pt}\,\,}\def\Lb{\underline{L}}\def\nab{\nabla}
\def\tr{\text{tr}}\def\a{\alpha}\def\b{\beta}\def\ga{\gamma}
\def\de{\delta}\def\Box{\square}\def\Boxr{\widetilde{\square}}\def\pa{\partial}
\def\pab{\bar\pa}\def\Lb{\underline{L}}

\def\pas{\text{$\pa\mkern -9.0mu$\slash}}\def\pam{(\pa_{t\!}-_{\!}\pa_r)}\def\pap{(\pa_{t\!}+_{\!}\pa_r)}

\numberwithin{equation}{section}
\begin{document}

\title {Global stability of Minkowski space-time in harmonic gauge}
\author {Hans Lindblad\thanks
{Part of this work was done while H.L. was a Member of the
Institute for  Advanced Study, Princeton, supported by the NSF grant
DMS-0111298 to the Institute. H.L. was also partially supported
by the NSF Grant DMS-0200226. } \,\,and Igor Rodnianski\thanks{Part of this work was done while I.R. was a
Clay Mathematics Institute Long-Term
Prize Fellow.  His work was also partially supported by the NSF grant DMS--01007791.}\\
University of California at San Diego and  Princeton University}
\date{}
\maketitle

\section{Introduction}
In this paper we address the question of  stability of Minkowski space-time for the
system of the Einstein-scalar field equations\footnote{We use Greek indices
$\alpha,\beta,\mu,\nu...=0,...,3$, the summation convention over repeated
indices and the notation $\pa_\alpha=\pa/\pa x^\alpha$. The symbol $D$ denotes a covariant Levy-Civita
derivatives with respect to the metric $g$ }
\beq\label{eq:Einst-T}
R_{\mu\nu}-\frac 12 g_{\mu\nu} R=T_{\mu\nu}.
 \end{equation}
The equations connect the
gravitational tensor $G_{\mu\nu}=R_{\mu\nu}-\frac 12 g_{\mu\nu}$
given in terms of the Ricci $R_{\mu\nu}$ and scalar $R=g^{\mu\nu} R_{\mu\nu}$ curvatures of an unknown Lorentzian
metric $g_{\mu\nu}$ and the energy-momentum tensor $T_{\mu\nu}$ of
a matter field $\psi$:
\begin{equation}\label{eq:energymomentumscalarfield}
 T_{\mu\nu}=\pa_\mu\psi\, \pa_\nu\psi
 -\frac 12{g_{\mu\nu}} \big( g^{\alpha\beta} \pa_\alpha \psi\,
 \pa_\beta\psi \big)
 \end{equation}
The Bianchi identities
 $$
 D^\mu
 G_{\mu\nu}=0
 $$
imply that the scalar field $\psi$ satisfies the the covariant wave
equation
$$
\Box_g \psi = \frac 1{\sqrt{|\det g|}} \pa_\mu \Big (g^{\mu\nu}
\sqrt{|\det g|} \pa_\nu \psi\Big )=0
$$
The set $(m,\R^{3+1},0)$:
\,standard Minkowski metric $g=m=-dt^2+\sum_{i=1}^3 (dx^i)^2$
on $\R^{3+1}$ and vanishing scalar field $\psi\equiv 0$ describes the
Minkowski space-time solution of the system \eqref{eq:Einst-T}.

The problem of stability of Minkowski space appears in the {\it Cauchy
formulation of the Einstein equations in which
given a 3-d manifold $\Si_0$ with a
Riemannian metric $g_0$, a symmetric 2-tensor $k_0$ and the initial
data $(\psi_0,\psi_1)$ for the scalar field, one  needs to find
a 4-d manifold ${\cal M}$, with a Lorantzian metric $g$ and a scalar field
$\psi$ satisfying the Einstein equations \eqref{eq:Einst-T}, and an imbedding $\Sigma_0\subset
M$ such that $g_0$ is the restriction of $g$ to $\Sigma$, $k_0$
is the second fundamental form of $\Sigma$ and the restriction of
$\psi$ to $\Si_0$ gives rise to the data $(\psi_0,\psi_1)$.}

The initial value problem is over determined and the data must satisfy
the constraint equations:
$$
R_0- {k_0}_{j}^i\, {k_0}^{j}_i + {k_0}_{i}^i \,{k_0}_j^j=|\nab\psi_0|^2 + |\psi_1|^2, \qquad
\nab^j {k_0}_{ij} - \nab_i \,{k_0}_j^j=\nab_i\psi_0 \,\psi_1.
$$
Here $R_0$ is the scalar curvature of $g_0$ and $\nab$ is
covariant differentiation with respect to $g_0$.

The seminal result of Choquet-Bruhat \cite{CB1} followed  by the
work \cite{CB-G} showed existence and uniqueness (up to a
diffeomorphism)
 of a maximal globally
hyperbolic\footnote{A space-time
is called globally hyperbolic if every inextendable causal curve
intersects the initial surface $\Si$ once and only
once. A causal curve $x(s)$ is a curve such that $
g_{\alpha\beta} \dot{x}^\alpha\dot{x}^\beta\leq 0 $, where $
\dot{x}=d x/ds$. It is future directed if $\dot{x}^0>0$.}
smooth space-time arising from any set of smooth initial data.
The work of Choquet-Bruhat used the diffeomorphism invariance of the
Einstein equations which allowed her to choose a special
{\it harmonic} (also referred to as a {\it wave coordinate} or de Donder)
gauge, in which the Einstein equations become
a system of quasilinear wave equations on the  components of the unknown metric
$g_{\mu\nu}$
\begin{equation}\label{eq:RE1}
\Boxr_g \, g_{\mu\nu} = F_{\mu\nu} (g)(\pa g, \pa
g)+2\pa_\mu\psi\, \pa_\nu\psi,\qquad \quad \Boxr_g \psi=0,
\qquad\text{where}\qquad
\Boxr_g=g^{\alpha\beta}\pa_\alpha\pa_\beta
\end{equation}
with $F(u)(v,v)$ depending quadratically on $v$.
Wave coordinates $\{x^\mu\}_{\mu=0,...,3}$ are required
to be solutions of the wave equations
$
\Box_g\,  x^\mu=0,
 $
 where the  geometric wave operator
 is
 $
 \Box_g=D_\alpha D^\alpha
 =g^{\alpha\beta}\pa_\alpha\pa_\beta+g^{\alpha\beta}\Gamma^{\,\,\,
 \nu}_{\alpha\beta}\pa_\nu
$.
The metric $g_{\mu\nu}$ relative to wave coordinates $\{x^\mu\}$ satisfies
the {\it wave coordinate condition}
\begin{equation}\label{eq:wave-coord}
g^{\alpha\beta}g_{\nu\mu}\Gamma^{\,\,\,
 \nu}_{\alpha\beta}=
g^{\alpha\beta}\pa_{\beta} g_{\alpha\mu}-\tfrac 12
g^{\alpha\beta}\pa_{\mu} g_{\alpha\beta}=0.
\end{equation}
Under this condition the geometric wave operator $\Box_g$ is equal
to the reduced wave operator $\Boxr_g$. The use of harmonic
 gauge goes back to the work of Einstein on post-Newtonian and
 post-Minkowskian expansions.

In the PDE terminology the result of Choquet-Bruhat corresponds to the {\it local
well-posedness} of the Cauchy problem for the Einstein-vacuum (scalar field)
equations with smooth initial data. The term local is appropriate in the sense
that the result does not
guarantee that the constructed space-time is {\it causally geodesically complete} \footnote{i.e. any causal geodesics
$x(s)$,  $g_{\alpha\beta} \dot{x}^\alpha\dot{x}^\beta=const\leq 0$
can be extended to infinite parameter value $0\leq s<\infty$.} and thus
could "terminate" in a singularity. Our experience suggests that global results
require existence of conserved or more generally monotonic positive quantities.
The only known such quantity in the asymptotically flat case is the ADM mass,
whose positivity was established  by Schoen-Yau \cite{S-Y} and Witten \cite{Wi},
is highly supercritical
relative to the equations and thus not sufficient to upgrade a local space-time to a
global solution.
This leaves the questions related to the structure of maximal globally
hyperbolic space-times even for generic\footnote{Note that all the known explicit solutions,
with exception of the Minkowski space-time, are in fact incomplete.}
data firmly in the realm of the outstanding Cosmic
Censorship Conjectures of Penrose. Given this state of affairs the problem of stability
of special solutions, most importantly the Minkowski space-time, becomes of crucial
importance.

{\bf Stability of Minkowski space-time for the Einstein-vacuum (scalar field) equations:}

\noindent
{\it Show existence of a causally geodesically space-time asymptotically "converging"
to the Minkowski space-time for an arbitrary set of smooth asymptotically flat initial data
$(\Si_0,{g_0}_{ij},{k_0}_{ij})$ with $\Si_0\approx \R^3$,
\begin{equation}
{g_0}_{ij}=(1+\frac Mr) \delta_{ij} + o(r^{-1-\alpha}),\qquad
{k_0}_{ij}=o(r^{-2-\alpha}),\quad r=|x|\to \infty, \quad \alpha>0\label{eq:asflat}
\end{equation}
where $(g_0-\de)$ and $k_0$ satisfy global smallness assumptions .
The stability of Minkowski space-time for the Einstein-scalar field equations:
in addition requires a global smallness assumption on the scalar field data
$(\psi_0,\psi_1)$, which obey the asymptotic expansion}
\begin{equation}
\psi_0=o(r^{-1-\alpha}),\qquad \psi_1=o(r^{-2-\alpha}).\label{eq:asflat-sc}
\end{equation}
A positive parameter $M$ in the asymptotic expansion for the metric $g_0$ is the
ADM mass.

The stability of Minkowski space for the Einstein-vacuum equations
was shown in a remarkable work of Christodouolou-Klainerman for
strongly asymptotic initial data (the parameter $\alpha\ge 1/2$ in the
asymptotic expansion \eqref{eq:asflat}.)   The approach taken in that work
viewed the Einstein-vacuum equations as a system of equations
$$
D^\a W_{\a\b\ga\de}=0,\qquad D^\a*W_{\a\b\ga\de}=0
$$
for the Weyl tensor $W_{\a\b\ga\de}$ of the metric $g_{\a\b}$ and
used generalized energy inequalities associated with the
Bel-Robinson energy-momentum tensor, constructed from components of
$W$, and special geometrically constructed vector fields, designed
to mimic the rotation and the conformal Morawetz vector fields of
the Minkowski space-time, i.e., "almost conformally Killing" vector
fields of the unknown metric $g$. The proof was manifestly
invariant, in particular it did not use the wave coordinate
gauge. This approach was later extended to the Einstein-Maxwell equations
by N. Zipser, \cite{Z}.

Nevertheless the PDE appeal (for the other motivations see the
discussion below) of the harmonic gauge for the proof of
stability of Minkowski space-time lies in the fact that  the latter
can be simply viewed\footnote{This statement requires additional
care since a priori there is no guarantee that obtained "global in
time" solution $g_{\mu\nu}$ defines a causally geodesically complete
metric. However, the latter can be established provided one has good
control on the difference between $g_{\mu\nu}$ and the Minkowski
metric $m_{\mu\nu}$, see \cite{L-R2}.} as a {\it small data global
existence result} for the quasilinear system \eqref{eq:RE1}.
However, usefulness of the harmonic  gauge in this
context was questioned earlier and it was suspected that wave
coordinates are "unstable in the large", \cite{CB2}.  The conclusion is
suggested from the analysis of the iteration scheme for the system
\eqref{eq:RE1}:
$$
g_{\mu\nu}=m_{\mu\nu} + \ve g_{\mu\nu}^{(1)} + \ve^2  g_{\mu\nu}^{(2)}+ ...,
$$
where $g_{\mu\nu}^{(1)}, g_{\mu\nu}^{(2)}$ satisfy respectively homogeneous
and inhomogeneous wave equations on Minkowski background
\begin{equation}
\Box\, g_{\mu\nu}^{(1)}=0, \qquad \Box\, g_{\mu\nu}^{(2)}=
F(m_{\mu\nu})\big (\pa g_{\mu\nu}^{(1)}, \pa g_{\mu\nu}^{(1)}\big
),\qquad\text{where}\qquad
\Box=m^{\alpha\beta}\pa_\alpha\pa_\beta\label{eq:iter}
\end{equation}
As a solution of the homogeneous $3+1$-d wave equation with smooth
decaying initial data the functions $g^{(1)}_{\mu\nu}\approx
C\varepsilon\, t^{-1}$ in the so called wave zone $t\approx r$ as
$t\to \infty$. Integrating the second equation implies that
$g^{(2)}_{\mu\nu}\approx C'\varepsilon\, t^{-1}\ln t$. This suggests
that already the asymptotic behavior of the second iterate deviates
from the one of the free waves in Minkowski space-time and plants
the seeds of doubt about global stability of such a scheme. We
should note that this formal iteration procedure is termed the
post-Minkowskian expansion and plays an important role in the study
of gravitational radiation from isolated sources, see e.g.
\cite{Bl}, \cite{B-D}.

To understand some of the difficulties in establishing a small data global
existence result for the system \eqref{eq:RE1} let us consider a generic
quasilinear system of the form
\begin{equation}\label{eq:Quad}
\Box\, \phi_i = \sum_{} b_{i}^{jk\alpha\beta\,} \pa_\alpha\phi_j\,
\pa_\beta\phi_k + \sum_{} c_{i}^{jk\alpha\beta\,} \phi_j\,
\pa_\alpha\pa_\beta\phi_k +{\text{cubic terms}}
\end{equation}
The influence of cubic terms is negligible while
the quadratic terms are of two types, the
{\it semilinear terms} and the {\it quasilinear terms}, each of which
present their own problems. The semilinear terms can cause blow-up in finite
time for smooth arbitrarily small initial data, as was shown by
John\cite{J1}, for the equation $\Box\, \phi = (\pa_{t}\phi)^{2}$.
D. Christodoulou \cite{C1} and S. Klainerman \cite{K2} showed global
existence for systems of the form \eqref{eq:Quad} if the semilinear
terms satisfy the {\it null condition} and the quasilinear terms are
absent. The null condition, first introduced by S. Klainerman in \cite{K1}, was designed to detect systems for which
solutions are asymptotically free and decay like solutions of a linear
equation. It requires special algebraic cancellations in the
coefficients $ b_i^{jk\alpha\beta}$, e.g. $\Box\, \phi=(\pa_t
\phi)^2-|\nabla_x \phi|^2$. However, the semilinear terms for
the Einstein equations do not satisfy the null condition, see
\cite{CB3}. The quasilinear terms is another source of trouble.
The only non-trivial example of a quasilinear equation of the type
\eqref{eq:Quad}, for which the small data global existence result
holds, is the model equation $\Box\,
\phi=\phi\,\Delta \phi$, as shown in \cite{L2} (radial case) and
\cite{A3} (general case).  However the solutions only decay like $\phi\sim
\varepsilon t^{-1+C\varepsilon}$. This in particular implies that
 the characteristic
surfaces of the associated Lorentzian metric
$g=-dt^2 + (1+\phi) \sum_{i=1}^3(dx^i)^2$
diverge distance $\!\sim t^{\,C\varepsilon}$ from the
Minkowski cones.


In our previous work \cite{L-R1} we identified a criteria under
which it is more likely that a quasilinear system of the form
\eqref{eq:Quas} has global solutions\footnote{At this point, it is unclear
whether this criteria is sufficient for establishing a "small data global
existence" result for a {\it general} system of quasilinear hyperbolic equations.}.
We said that a system of the
form \eqref{eq:Quas} satisfy the {\it weak null condition} if the
corresponding {\it asymptotic system} (c.f. \cite{H1, H2}) has
global solutions\footnote
{The condition trivially holds for the class of equations satisfying the standard null condition.  For more discussion of the weak null condition and a general intuition
behind the proof see Section 2.}
We showed that the Einstein equations in wave
coordinates satisfy the weak null condition. In addition there is
some additional cancelation for the Einstein equations in wave
coordinates that makes it better than a general system satisfying
the weak null condition. The system decouples to leading
order, when decomposed relative to the Minkowski {\it null frame}.
An approximate model that describes the semilinear terms has the form
$$
\Box\, \phi_2=(\pa_t \phi_1)^2,\qquad \Box \, \phi_1=0.
 $$
While every solution of this system is global in time,  the system
fails to satisfy the classical null
condition and solutions are not asymptotically free: $\phi_2\sim
\varepsilon t^{-1}\ln{|t|}$. The semilinear terms in Einstein's
equations can be shown to either satisfy the classical null
condition or decouple in the above fashion when expressed in a
null frame. The quasilinear terms also decouple but in a more
subtle way.  The influence of quasilinear terms can be detected via
asymptotic behavior of the charachteristic surfaces of metric $g$.
It turns out that the main features of the characteristic surfaces at
infinity are determined by a particular {\it null} component of the metric.
The asymptotic flatness of the initial
data and the wave coordinate condition \eqref{eq:wave-coord} give
good control of this particular component, i.e.,  $\sim  M/r$, which
in turn implies that the light cones associated with the metric
$g$ diverge only
logarithmically $\sim M\ln t$ from the Minkowski cones.

In \cite{L-R2} we were able to establish the
"small data global existence" result for the Einstein-vacuum equations in harmonic
gauge for special restricted type of initial data. In addition to the standard constraint compatibility
and smallness conditions the initial data was assumed to coincide
with the Schwarzschild data
$$
{g_0}_{ij}=(1+\frac Mr)\de_{ij},\qquad {k_0}_{ij}=0
$$
in the complement of the ball of radius one centered at the origin.
The existence of such data was recently demonstrated in \cite{Co},
\cite{C-D}. The stability of Minkowski space-time for the
Einstein-vacuum equations for such data of course also follows from the
work of Christodoulou-Klainerman and yet another approach of
Friedrich \cite{Fr}. This choice of data allowed us to completely ignore the
problem of the long range effect of the mass and the exterior existence.

In this paper we prove stability of Minkowski space-time for the
Einstein-vacuum and the Einstein-scalar field equations in harmonic  gauge
 for general asymptotically flat
initial data (any $\a>0$ in
\eqref{eq:asflat}-\eqref{eq:asflat-sc}) close to the data for the Minkowski solution.


 The asymptotic behavior of null components
 of the Riemann curvature tensor $R_{\a\b\ga\de}$ of metric $g$\-- the so called
 "peeling estimates"\--  was discussed in the works of Bondi, Sachs and Penrose
 and becomes important in the framework of asymptotically
 simple space-times (roughly speaking, space-times which can be conformally
 compactified), see also the paper of Christodoulou \cite{C2} for further discussion
 of such space-times. The work of \cite{C-K} provided
  very precise, although not
 entirely consistent with peeling estimates, analysis of  the asymptotic behavior
 of constructed global solutions. However, global solutions obtained
 by Klainerman-Nicolo \cite{K-N1} in the problem of
 exterior\footnote{Outside of the domain
 of dependence of a compact set} stability of Minkowski space
 were shown to possess peeling estimates for special initial data, \cite{K-N2}.

 Our work is less precise about the asymptotic behavior and  is focused more on
 developing a relatively technically simple approach allowing us to prove stability
 of Minkowski space in a physically interesting wave coordinate gauge, for
 general asymptotically flat data, and simultaneously treating the case of the
 Einstein equations coupled to a scalar field.

\begin{theorem} \label{globalexist}
Let $(\Si, g_0, k_0, \psi_0, \psi_1)$ be initial data for the
Einstein-scalar field equations. Assume that the initial time slice
$\Si$ is diffeomorphic to ${\Bbb R}^3$ and admits a global
coordinate chart relative to which the data is close to the initial
data for the Minkowski space-time. More precisely, we assume that
the data $(g_0, k_0,\psi_0,\psi_1)$ is smooth asymptotically flat in
the sense of \eqref{eq:asflat}-\eqref{eq:asflat-sc} with mass $M$
and impose the following smallness assumption: Let
$$
g_0= \delta+ h^0_0+h^1_0,\qquad \text{where}\quad
{h^0_{ij}}_0=\chi(r)\frac{M}{r}\,  \delta_{ij}
$$
where $\chi(s)\in C^\infty$ is $1$ when $s\geq 3/4$ and $0$ when
$s\leq 1/2$. Set
\begin{align}
E_N(0)=\sum_{0\le|I|\leq N} &\|(1+r)^{1/2+\ga+|I|}\nabla\,\nabla^I
h^1_0\|_{L^2}+\|(1+r)^{1/2+\ga+|I|}\nabla^I k_0\|_{L^2}\label{eq:enerN0} \\+&
\|(1+r)^{1/2+\ga+|I|}\nabla\, \nabla^I \psi_0\|_{L^2}+
\|(1+r)^{1/2+\ga+|I|} \nabla^I \psi_1\|_{L^2}\nn
\end{align}
There is a constant
$\varepsilon_0>0$ such that for all  $\varepsilon\leq \varepsilon_0$
and initial data verifying the condition
\begin{equation}
\label{eq:enerNini} E_N(0)+M\leq \varepsilon
\end{equation}
for some $\ga>\ga_0(\eps_0)$ with $\ga_0(\eps_0)\to 0$
as $\eps_0\to 0$ and $N\ge 6$, the Einstein-scalar
field equations possess a future causally geodesically complete solution
$(g,\psi)$  asymptotically
converging to Minkowski space-time. More precisely, there exists
a global system of coordinates $(t,x^1,x^2,x^3)$ in which the
energy
\begin{equation}
\label{eq:enerNdef} \E_N(t)=\sum_{|I|\leq N} \|w^{1/2}\pa Z^I
h^1(t,\cdot)\|_{L^2} + \|w^{1/2}\pa Z^I \psi(t,\cdot)\|_{L^2},\qquad
 w^{1/2}=\begin{cases} (1+|r-t|)^{1/2+\gamma},\,\,\, &r>t\\
1,\quad &r\leq t\end{cases}
\end{equation}
with Minkowski vector fields  $Z\in \{\pa_\a, x_\a
\pa_\b-x_\b\pa_\a, x^\a \pa_\a\}$ of the solution
$$
g(t)=m+h^0(t)+ h^1(t),\qquad
h^0_{\a\b}=\chi(r/t)\chi(r)\frac{M}{r}\de_{\a\b}.
$$
obeys the estimate
\begin{equation}
\label{eq:globalenergy} \E_N(t)\leq C_N
\varepsilon\,(1+t)^{C_N\varepsilon}.
\end{equation}
Moreover,
\begin{equation}\label{eq:globaldecay}
\,\,\,\,\,\,| \pa Z^I h^1| + |\pa Z^I \psi|\leq \begin{cases}
C_N^\prime\varepsilon
(1+t+r)^{-1+C_N\varepsilon}(1+|t-r|)^{-1-\gamma},\qquad r>t,\\
C_N^\prime \varepsilon
(1+t+r)^{-1+C_N\varepsilon}(1+|t-r|)^{-1/2},\qquad\,\, r\le t,
\end{cases}\qquad  |I|\leq N-2
\end{equation}
and
\begin{equation}\label{eq:globaldecay}
\,\,\,\,\,\,| Z^I h^1| + |Z^I \psi|\leq \begin{cases}
C_N^{\prime\prime} \varepsilon (1+t+r)^{-1+
C_N\varepsilon}(1+|t-r|)^{-\gamma},\qquad
r\ge t, \qquad\quad |I|\leq N-2\\
 C_N^{\prime\prime}
\varepsilon (1+t+r)^{-1+2
C_N\varepsilon},\qquad\qquad\qquad\qquad\!\! r\le t,\qquad\quad
|I|\leq N-3\end{cases}
\end{equation}
\end{theorem}
 For $|I|=0$ we prove the stronger bound $C\varepsilon
t^{-1}\ln{|t|}$ and in fact show that all but one component of
$h^1$, expressed relative to a null frame,
 can be bounded by
$C\varepsilon t^{-1}$. This more precise information is needed in
the proof.

As we will show in a future paper, the decay estimates above are
sufficient to prove the peeling properties of the curvature tensor
up to order $t^{-3}$ along any forward light cone, but peeling
properties of higher order require stronger decay of initial data.

{\bf Acknowledgments}:\, The authors would like to thank Sergiu
Klainerman for valuable suggestions and discussions.

\section{Strategy of the proof and outline of the paper}
\subsection{Notations and conventions}
Coordinates:
\begin{itemize}
\item $\{x^\alpha\}_{\alpha=0,..,3}=(t,x)$ with $t=x^0$ and $x=(x_1,x_2,x_3),\,
r=|x|$ are the
standard space-time coordinates
\item $s=r+t, q=r-t$ are the null coordinates
\end{itemize}
Derivatives:
\begin{itemize}
\item $\nabla= (\pa_1,\pa_2,\pa_3)$ denotes spatial derivatives
\item $\pa=(\pa_t,\nabla)$ denotes space-time derivatives
\item $\pas_i$ denotes spatial angular components of the
derivatives $\pa_i$.
\item $\pab= (\pa_t+\pa_r, \pas)$ denotes derivatives tangent to the light
cones $t-r=$constant
\item $\pa_s=\frac 12(\pa_r+\pa_t),\, \pa_q=\frac 12(\pa_r-\pa_t)$ denote the
null derivatives
\end{itemize}
Metrics:
\begin{itemize}
\item $m=-dt^2 +\sum_{i=1}^3 (dx^i)^2$ denotes the standard Minkowski metric
on $\R^{3+1}$
\item $g$ denotes a Lorentzian metric \- solution of the Einstein equations
\item Raising and lowering of indices in this paper is always done with respect
to the metric $m$, i.e., for an arbitrary n-tensor
$\Pi_{\alpha_1,\alpha_2..,\alpha_n}$ we define
$\Pi^{\alpha_1}_{\,\,\,\,\alpha_2,..,\alpha_n} = m^{\alpha_1\beta}
\Pi_{\b,\alpha_2,..,\alpha_n}$. The only exception is made for the
tensor $g^{\a\b}$ which stands for the inverse of the metric
$g_{\a\b}$.
\end{itemize}
Null frame:
\begin{itemize}
\item $L=\pa_t +\pa_r$ denotes the vector field generating the forward Minkowski light cones
 $t-r=$constant
 \item $\Lb=\pa_t-\pa_r$  denotes the vector field transversal to
 the light cones  $t-r=$constant
 \item $S_1,S_2$ denotes orthonormal vector fields spanning the tangent
 space of the spheres $t=$constant, $r=$constant
 \item The collection ${\cal T} = \{L, S_1, S_2 \}$ denotes the frame
 vector fields tangent to
 the light cones
 \item The collection  ${\cal U}=  \{L,{\Lb}, S_1, S_2 \}$ denotes the full null frame.
 \end{itemize}
 Null forms
 \begin{itemize}
 \item
 $Q_{\a\b}(\pa\phi,\pa\psi)=\pa_\a\phi\pa_\b\psi-\pa_\a\phi\pa_\a\psi,\,\,
 Q_0(\pa\phi,\pa\psi)=m^{\a\b} \pa_\phi\pa_\b\psi$ are the standard null
 forms
 \end{itemize}
Null frame decompositions:
\begin{itemize}
\item For an arbitrary vector field $X$ and frame vector  $U$,
$X_U=X_\alpha U^\alpha$, where $X_\alpha=m_{\alpha\beta} X^\beta$.
\item For an arbitrary vector field $X= X^\a \pa_a = X^L L + X^{\Lb} \Lb
+ X^{S_1} S_1+X^{S_2} S_2$, where $X^L=-X_{\Lb}/2,\,
X^{\Lb}=-X_L/2$, $X^{S_i}=X_{S_i}$.
\item For an arbitrary pair of vector fields $X, Y$:
$X^\a Y_\a = -X_L Y_{\Lb}/2 -X_{\Lb} Y_L/2 + X_{S_1}
Y_{S_1}+X_{S_2}Y_{S_2}$
\item Similar identities hold for higher order tensors
\end{itemize}
Minkowski vector fields $\{Z\}={\cal Z}$:
\begin{itemize}
\item ${\cal  Z}=\{\pa_\a, \Omega_{\a\b}=x_\a\pa_\b-x_\b \pa_\a, S= x^\a \pa_a\}$.
\item Commutation properties: $[\Box, \pa_\a]=[\Box,\Omega_{\a\b}]=0$,
$[\Box,S]=2\Box$,\,\, $c_Z \Box:=[Z,\Box]$.
\end{itemize}
\subsection{Strategy of the proof}
Our approach to constructing solutions of  the Einstein (Einstein-scalar field) equations is
to  solve the corresponding system of the reduced Einstein equations \eqref{eq:RE1}
\begin{align}
&\Boxr_g \, g_{\mu\nu} = \tilde F_{\mu\nu} (g)(\pa g, \pa
g)+2\pa_\mu\psi\,  \pa_\nu\psi,\label{eq:Fullred}\\ &\Boxr_g \psi=0\nn
\end{align}
It is well known that any solution $(g,\psi)$ of the full system of Einstein equations, written in wave
coordinates, will satisfy
\eqref{eq:Fullred}. The wave coordinate condition can be recast as a requirement
that relative to a coordinate system $\{x^\mu\}_{\mu=0,...,3}$
the tensor  $g_{\mu\nu}$ verifies
\begin{equation}\label{eq:WE1}
g^{\a\b} \pa_\b g_{\a\mu} = \frac 12 g^{\a\b} \pa_\mu g_{\a\b}.
\end{equation}
Conversely, any solution of the system \eqref{eq:Fullred} with
initial data $(g_{\mu\nu}|_{t=0},\pa_t g_{\mu\nu}|_{t=0})$  and
$(\psi|_{t=0},\pa_t\psi|_{t=0})$ compatible with \eqref{eq:WE1} and
constraint equations, gives rise to a solution of the full Einstein
system. In addition, thus constructed  tensor $g_{\mu\nu}$ will
satisfy the wave coordinate condition \eqref{eq:WE1} for all times.
Verification of the above statements is straightforward and can be
found in e.g. \cite{Wa}.

\noindent {\bf Construction of initial data}

The assumptions of the Theorem assert that the initial data $({g_0}_{ij},{k_0}_{ij},
\psi_0,\psi_1)$ verify the global smallness condition \eqref{eq:enerN0}.
We define the initial data $(g_{\mu\nu}|_{t=0},\pa_t g_{\mu\nu}|_{t=0}, \psi|_{t=0},\pa_t\psi|_{t=0})$ for the system \eqref{eq:Fullred} as follows:
\begin{align}
&g_{ij}|_{t=0}={g_0}_{ij},\quad g_{00}|_{t=0}=-a^2,\quad g_{0i}|_{t=0}=0,
\quad\psi|_{t=0}=\psi_0\label{eq:init-red1}\\
&\pa_t g_{ij}|_{t=0}=-2a {k_0}_{ij},\quad \pa_t g_{00}|_{t=0}=
2 a^3 g^{ij}_0 {k_0}_{ij},\quad
\pa_t\psi|_{t=0}=a\psi_1,\label{eq:init-red2}\\
&\pa_t g_{0\ell}=a^2 g^{ij}_0\pa_j\, {g_0}_{i\ell}- \frac 12 a^2
g_0^{ij} \pa_\ell \,{g_0}_{ij}-a \pa_\ell\, a. \label{eq:init-red3}
\end{align}
The lapse function $a^2=(1-M\chi(r) r^{-1})$.
The data constructed above is compatible with the constraint equations and
satisfies the
wave coordinate condition \eqref{eq:WE1}. In fact, it simply corresponds to a choice
of a local system of coordinates $\{x^\mu\}_{\mu=0,...,3}$ satisfying the wave
coordinate condition $\Box_g x^\mu=0$ at $x^0=t=0$. This procedure is standard
and its slightly simpler version can be found  e.g. in \cite{Wa}, see also \cite{L-R2}.

The smallness condition
$E_N(0) + M\le \epsilon$ for the initial data $({g_0}_{ij},{k_0}_{ij},
\psi_0,\psi_1)$ implies the smallness condition
$$
\E_N(0) + M\le \epsilon.
$$
for the energy $\E_N(0)$:
\begin{align*}
 \E_N(0)&=\sum_{0\le|I|\leq N} \left(\|(1+r)^{1/2+\ga+|I|}\nabla\,\nabla^I
h^1|_{t=0}\|_{L^2}+ \sum_{0\le|I|\leq N} \|(1+r)^{1/2+\ga+|I|}\,\nabla^I
\pa_t h^1|_{t=0}\|_{L^2}\right) \\ &+
\sum_{0\le|I|\leq N} \left(\|(1+r)^{1/2+\ga+|I|}\nab\,\nab^I \psi|_{t=0}\|_{L^2}+\|(1+r)^{1/2+\ga+|I|}\,\nab^I \pa_t\psi|_{t=0}\|_{L^2}\right).
\end{align*}
defined for a tensor ${h^1}_{\mu\nu}$ which appears in the decomposition
$$
g= m+ h^0+h^1,\qquad \text{where}\quad
{{h^0}}_{\mu\nu}=\chi(\frac rt)\chi(r)\frac{M}{r}\, \delta_{\mu\nu}.
$$
This allows us to reformulate the problem of global stability of Minkowski space
as a small data global existence problem for the system
\begin{align}
&\Boxr_g \, h^1_{\mu\nu} = F_{\mu\nu} (g)(\pa h, \pa
h)+2\pa_\mu\psi\,  \pa_\nu\psi - \Boxr_g h^0,\label{eq:RE3}\\
&\Boxr_g \psi=0\nn
\end{align}
with $h^0_{\mu\nu} = \chi (r/t) \chi (r) \frac Mr \delta_{\mu\nu}$.
Arguing as in section 4 of \cite{L-R2} we can show that the tensor
$g_{\mu\nu}(t)=m+h^0_{\mu\nu} (t) + h^1_{\mu\nu}(t)$ obtained
by solving \eqref{eq:RE3} with initial data given in \eqref{eq:init-red1}-
\eqref{eq:init-red3} defines a solution of the Einstein-scalar field
equations. Moreover, the wave coordinate condition \eqref{eq:WE1}
is propagated in time.

The above discussion leads to the following result:
\begin{theorem}
A global in time solution of \eqref{eq:RE3} obeying the energy bound
$\E_N(t)\le C\ve (1+t)^\delta$ for some sufficiently small $\delta>0$
gives rise to a future causally geodesically complete solution of the Einstein
equations \eqref{eq:Einst-T} converging to the Minkowski space-time.
\end{theorem}
The proof of the geodesic completeness can be established by
arguments identical to those in \cite{L-R2} and will not be repeated
here.

\noindent {\bf The mass problem}

A first glance at the system \eqref{eq:Fullred} suggests that a natural approach
to recast the problem of global stability of Minkowski space as a small data
global existence question for a system of quasilinear wave equations is to
rewrite \eqref{eq:Fullred} as an equation for the tensor $h=g-m$:
  \begin{align}
&\Boxr_g \, h_{\mu\nu} = F_{\mu\nu} (\pa h, \pa h)+2\pa_\mu\psi\,  \pa_\nu\psi,\\\
&\Boxr_g \psi=0\nn
\end{align}
The initial data for $h$ possess the asymptotic expansion as
$r\to\infty$: \beq\label{eq:expand} {h}_{\mu\nu}|_{t=0} = \frac Mr
\delta_{\mu\nu}+ O(r^{-1-\alpha}),\quad {\pa_t h}_{\mu\nu}|_{t=0} =
O(r^{-2-\alpha})
\end{equation}
for some positive $\a>0$. While the data appears to be "small" it
does not have sufficient decay rate in $r$ at infinity due to the
presence of the term with positive mass $M$. This could potentially
lead to a  "long range effect" problem, i.e., the decay of the
solution at the time-like infinity $t\to \infty$ is affected by the
slow fall-off at the space-like infinity $r\to \infty$. In the work
\cite{C-K} this problem was resolved (albeit in a different
language) by taking advantage of the fact that the long range term
is spherically symmetric and the ADM mass $M$ is conserved along the
Einstein flow. Thus differentiating the solution (in \cite{C-K} this
means the Weyl field corresponding to the conformally invariant part
of the Riemann curvature tensor of metric $g$) with respect to
properly defined (non-Minkowskian) angular momentum and time-like
vector fields one obtains a new field still approximately satisfying
the Einstein field equations but with considerably better decay
properties at space-like infinity. We however pursue a different
approach by taking an educated "guess" about the contribution of the
long range term $M\delta_{\mu\nu}/r$ to the solution. Thus we set
$$
h^0_{\mu\nu} = \chi(\frac r{t}) \chi(r) \frac Mr \delta_{\mu\nu},
$$
split the tensor
$h=h^1+h^0$ and write the equation for the new unknown $h^1$.
The important cancelation occurring in a new inhomogeneous term
$\Boxr_g h^0$ is the vanishing of $\Box \frac Mr=0$ away from $r=0$.
The cut-off function $\chi(r/t)$ ensures that the essential contribution
of $\Boxr_g h^0$ comes from the "good" interior region
$\frac 12 t\le r\le \frac 34 t$ \-- the support of the derivatives of $\chi(r/t)$.
Finally, observe that in the presence of the term containing the mass $M$
the tensor $h=h^0+h^1$ has infinite energy
$\E_N(0)$.

The system \eqref{eq:RE3} is a system of quasilinear wave equations.
The results of \cite{C1}, \cite{K2} show that a sufficient condition
for a system of quasilinear wave equations to have global solutions
for all smooth sufficiently small data is the null condition.
However quasilinear problems where the metric depends on the
solution, rather than its derivatives, as in the problems arising in
elasticity, do not satisfy the null condition. Moreover, as shown in
\cite{CB1,CB2} even the semilinear terms  $F_{\mu\nu} (\pa h, \pa
h)$ violate the standard null condition.

\noindent
{\bf Connection with the weak null condition}

Consider a general system of quasilinear wave equations:
\begin{equation}\label{eq:Quas}
\Box \phi_I = \!\!\!\!\!\sum_{|\alpha|\le|\beta|\le 2,\,\, |
\beta|\geq 1} \!\!\!\!\!\! A_{I,\alpha\beta}^{JK}\,  \pa^\alpha
\phi_J\, \pa^\beta\phi_K + {\text{cubic terms}}
\end{equation}
The weak null condition, introduced in \cite{L-R1}  requires that
the asymptotic system for $\Phi_I\!=\!r\phi_I$ corresponding to
\eqref{eq:Quas}:
\begin{equation}\label{eq:asymp}
\pap\pam\Phi_I \sim\\{r}^{-1} \!\!\!\!\sum_{n\le m\le 2,\,\, m\geq
1} \!\!\!\! A_{I\!,nm}^{JK} \pam^n \Phi_{\!J} \,\, \pam^m \Phi_K,
\end{equation}
has global solutions for all small data. Here, the tensor
$$
A_{I\!,nm}^{JK}(\omega):=
{(-2)^{-m-n}}\!\!\!\!\sum_{|\alpha|=n,\,|\beta|=m} \!\!\!\!\!\!
A_{I\!,\alpha\beta\,}^{JK} \hat\omega^\alpha \hat
\omega^\beta\!\!,\quad \hat
\omega\!=\!(\!-\!1,\omega\!),\,\omega\!\in\!\bold{S}^2
$$
The standard null condition, which guarantees a small data global
existence result for the system \eqref{eq:Quas}, is that
$A_{\!I\!,nm}^{JK}(\omega) \equiv 0$ and clearly is included in the
weak null condition, since in that case the corresponding asymptotic
system \eqref{eq:asymp} is represented by a {\it linear} equation.
Asymptotic systems were introduced by H\"ormander, see \cite{H1,H2},
as a tool to find the exact blow-up time of solutions of scalar
equations violating the standard null condition with quadratic terms
independent of $\phi$, i.e. $|\alpha|\geq 1$ in \eqref{eq:Quas}.
This program was completed by Alinhac, see e.g. \cite{A1}. It was
observed in \cite{L1} that the asymptotic systems for quasilinear
equations of the form $\Box\, \phi=\sum c^{\,\alpha\beta}
\phi\,\,\pa_\alpha\pa_\beta \phi$ in fact have global solutions. In
other words these equations satisfy the weak null condition. It was
therefore conjectured that these wave equations should have global
solutions for small data. This has been so far only proven for the
equation
\begin{equation}\label{eq:involved}
 \Box \phi=\phi\, \Delta \phi,
 \end{equation}
   \cite{L1} (radial case),
\cite{A3} (general case). Note that the case $|\alpha|=|\beta|=0$,
in \eqref{eq:Quas} is excluded. The asymptotic system only predicts
the behavior of the solution close to the light cones, and for the
case $\Box \phi=\phi^2$ the blow-up occurs in the interior and much
sooner, see \cite{J2,L3}

The asymptotic system \eqref{eq:asymp} is obtained from the system
\eqref{eq:Quas} by neglecting derivatives tangential to the outgoing
Minkowski light cones and cubic terms, that can be expected to decay
faster. In particular,
\begin{align*}
&\Box \phi\!=\!r^{-1}\pap\pam (_{\!}r\phi_{\!})+\text{angular
derivatives},\\
&\pa_\mu=-\tfrac{1}{2}\hat{\omega}_\mu\pam+\text{tangential
derivatives}
\end{align*}
Recall that for solutions of linear wave homogeneous wave equation
derivatives tangential to the forward light cone $t=r$ decay
faster: $\,\,\!\!\Box\phi\!=\!0\!$  implies that $|\pa\phi|\!\leq \!C/t$ while
$|\pab\! \phi|\leq\! C/t^2$.

A simple example of a system satisfying the weak null condition,
violating the standard null condition and yet possessing global
solutions is
\begin{equation}\label{eq:simple}
\begin{aligned}
\square \phi_1&=\phi_3\cdot\pa^2\phi_1+(\pa \phi_2)^2,\\
\Box\phi_2&=0,\qquad\Box \phi_3=0.
\end{aligned}
\end{equation}
 Another, far less trivial example
  is provided by the equation \eqref{eq:involved}.

 The asymptotic system for Einstein's equations can be modeled by
that of \eqref{eq:simple}. To see this we introduce a null-frame
$\{L, \Lb,S_1, S_2\}$ decomposition of Einstein's equation. With
$h\!=\!g\!-\!m$ we have
\begin{equation}
\Boxr_g h_{\mu\nu} =F_{\mu\nu} (h) (\pa h, \pa
h)+2\pa_\mu\psi\,\pa_\nu\pa\psi,\qquad\quad \Boxr_g \psi=0
\end{equation}
where
$$
F_{\mu\nu} (h) (\pa h, \pa h)=\frac 14 \pa_\mu h^{\a}_\a \pa_\nu h^\b_\b-
\frac 12 \pa_\mu h^{\a\b} \pa_\nu h_{\a\b}
+Q_{\mu\nu}(\pa h,\pa h) +G_{\mu\nu}(h)(\pa h,\pa h),
 $$
with $Q_{\mu\nu}$ - linear combinations of the standard
null-forms and $G_{\mu\nu}(h)(\pa h,\pa h)$  contains only cubic terms.

The asymptotic system for the reduced Einstein equations has the
following form:
 \begin{align}
& (\pa_t + \pa_r) \pa_q D_{\mu\nu}= H_{LL} \pa_q^2 D_{\mu\nu} - \tfrac{1}{4r}
 \hat{\omega}_\mu \hat{\omega}_\nu \big(P(\pa_q D,\pa_qD)+2\pa_q\Phi\,\pa_q\Phi\big),
 \quad\pa_q=\pa_t-\pa_r\label{eq:Asym-C1}\\
 & (\pa_t + \pa_r) \pa_q \Phi = H_{LL} \pa_q^2\Phi,\label{eq:Asym-C2}\\
&D_{\mu\nu}\sim r h_{\mu\nu},\quad
\Phi\sim r \psi ,\quad H^{\a\b} = g^{\a\b}-m^{\a\b} ,\quad H_{LL}=H_{\a\b}
L^\alpha L^\beta,\quad L^\alpha \pa_\a=\pa_t+\pa_r \nn
\end{align}
Here
 \beq
 P(D,E)=\frac{1}{4} D_{\alpha}^{\,\, \, \alpha}
E_{\beta}^{\,\, \, \beta}-\frac{1}{2} D^{\alpha\beta}
E_{\alpha\beta}.
 \eq
 On the other hand the asymptotic form of the wave coordinate
condition \eqref{eq:WE1}  is
\begin{equation}\label{eq:asymptoticwavec}
\pa_q D_{L T}\sim 0,\qquad T\in{ \cal T}=\{L,S_1,S_2\}
\end{equation}
Observe that \eqref{eq:asymptoticwavec} combined with the initial asymptotic expansion
\eqref{eq:expand} of the metric $g$ suggests that asymptotically,
as $t\to \infty$ and $|q|=|t-r|\le C$,
\begin{equation}\label{eq:AsymL}
h_{L\alpha} \sim H_{L\alpha} \sim \frac {M}t c_\alpha,
\qquad c_\alpha=\begin{cases} 0,\qquad\,\,\,\, \alpha=0,\\ {x_i}/{|x|},\quad\,\, \alpha=i
\end{cases}
\end{equation}
Decomposing the system with the help of the null frame $\{L, S_1,
S_2\}={\cal T}$ and $\,\Lb:\,\,\Lb^\a\pa_\a=\pa_r-\pa_t$ and using
that $L^\mu\hat\omega_\mu=A^\mu\hat\omega_\mu=0$, we obtain that
 \begin{align}
(\pa_t + \pa_r) \pa_q D_{\underline{L}\,\underline{L}}& =
H_{LL} \pa_q^2 D_{\Lb\Lb} -  2r^{-1}P(\pa_q D,\pa_q D)-r^{-1}\pa_q\Phi\,\pa_q\Phi
\label{eq:Asym1}\\
(\pa_t + \pa_r) \pa_q D_{TU}&= H_{LL} \pa_q^2 D_{TU},\qquad\qquad T\in{ \cal T}
=\{L,S_1,S_2\}, \,\,\, U\in{\cal U}=\{L,\Lb,S_1,S_2\},\label{eq:Asym2}\\
 (\pa_t + \pa_r) \pa_q \Phi& =
H_{LL} \pa_q^2 \Phi\label{eq:Asym-Phi}
 \end{align}
 In view of \eqref{eq:AsymL} the equations \eqref{eq:Asym2} allows us
 to find the asymptotic behavior of the components $\pa_q D_{TU}\sim const$
 and consequently  $h_{TU}\sim t^{-1}$. The asymptotic behavior $\Phi\sim t^{-1}$
 can be similarly determined from \eqref{eq:Asym-Phi}.
On the other hand,
$$
P(\pa_q D,\pa_q D) = \frac 14 \pa_q D_{LL}\pa_q D_{\Lb\Lb} +
\pa_q D_{TU} \cdot \pa_q D_{TU}
$$
It is the the absence of the quadratic term $(\pa_q D_{\Lb\Lb})^2$ and the boundedness
of the components $\pa_q D_{TU}$ that
allows us to solve the equation \eqref{eq:Asym1} for $D_{\Lb\Lb}$ although
the suggested asymptotic behavior is different
$h_{\Lb\Lb}\sim H_{\Lb\Lb}\sim t^{-1} \ln t $.

It turns out that the asymptotic system indeed correctly predicts the asymptotic
behavior of the tensor $g_{\mu\nu}=m_{\mu\nu} + h_{\mu\nu}$ \--
solution of the Einstein equations in wave coordinates

Before proceeding to explain the strategy of the proof of our result we review
the ingredients of a proof of a typical small data global existence for quasilinear
wave equations in dimensions $n\ge 4$ or equations satisfying the standard
null condition. For simplicity we consider a semilinear equation
\beq\label{eq:semil}
\Box\phi = N(\pa\phi,\pa\phi)
\eq
with a quadratic nonlinearity $N$. We first note that almost every small data
global existence result was established under the assumption of compactly
supported (or rapidly decaying) data, which by finite speed of propagation
ensures that a solution
is supported in the interior of a light cone $r=t+C$ for some sufficiently large
constant $C$. The proof is based on generalized energy estimates
\beq\label{eq:ener-in-sem}
E_N(t) \le E_N(0) +
\sum_{Z\in {\cal Z}, |I|\le N}\int_0^t \|Z^I F(\tau)\|_{L^2} E_N^{\frac 12}(\tau)\,d\tau
\eq
for solutions
of an inhomogeneous wave equation $\Box\phi=F$. The generalized energy
$$
E_N(t)=\sum_{Z\in {\cal Z}, |I|\le N} \|\pa Z^I \phi(t)\|^2_{L^2},
$$
with vector fields ${\cal Z}=\{\pa_\a,
\Omega_{\alpha\beta}=-x_\a\pa_\b+x_\b\pa_\a, S=x^\a\pa_\a\}$, of a
solution of \eqref{eq:semil} can be shown to satisfy the inequality
\beq\label{eq:ener-sem} E_N(t) \le \exp{\big [C \sup_{Z\in {\cal Z},
|I|\le N/2+2} \int_0^t\|\pa Z^I \phi(\tau)\|_{L^\infty_x}\big ]}
E_N(0) \eq The other crucial component is the Klainerman-Sobolev
inequality which asserts that for an arbitrary smooth function
$\psi$ \beq\label{eq:K-S-toy} |\pa\psi(t,x)|\le {C}{(1+t+r)^{-\frac
{n-1}2}(1+|q|)^{-\frac 12}}\,\,\, E_4(t),\qquad r=|x|,\,\,\, q=r-t.
\eq
 Combining \eqref{eq:K-S-toy} with the energy inequality leads to the
 proof of the small data global existence result for a generic semilinear
 wave equation in dimension $n\ge 4$. It also elucidates the difficulty
 of proving such a result in dimension $n=3$. In fact, as was shown
 in \cite{J1}, the result can be false in dimension $n=3$, e.g., the equation
 $\Box\phi=\phi_t^2$ admits small data solutions with finite time
 of existence. It is interesting to note that the corresponding asymptotic system
 for such equation leads to a Ricatti type ODE.

 In the case when the semilinear terms $N(\pa\phi,\pa\phi)$ obey the
 standard null condition, i.e., it is a linear combination of the quadratic
 null forms $Q_{\a\b}, Q_0$, it is possible to refine the energy inequality
 \eqref{eq:ener-sem} so that a combination of the energy and Klainerman-Sobolev
 inequalities still yields a small data global existence result in dimension
 $n=3$. This can be traced to the following pointwise estimate on a null
 form $Q$:
 $$
 |Q(\pa\phi,\pa\phi)|\le C
 (1+t+r)^{-1}\sum_{Z\in {\cal Z}} |\pa \phi|\, |Z\phi|
 $$
 In the absence of the standard null condition in dimension $n=3$
 the inequalities \eqref{eq:ener-sem}-\eqref{eq:K-S-toy} are just barely
 insufficient for the desired result.
 An illuminating example is provided by a semilinear version of the system \eqref{eq:simple}
 $$\square \phi_1=(\pa \phi_2)^2,\qquad
\Box\phi_2=0
$$
A combination of \eqref{eq:ener-sem}-\eqref{eq:K-S-toy}
applied to the vector $\phi=(\phi_1,\phi_2)$ would incorrectly suggest a
possible finite time blow-up. This analysis however reflects the following
phenomena:\,\,{\it  small data solutions of the system above have a polynomially
growing energy $E_N(T)\sim t^\de$ and the asymptotic behavior
$\pa\phi(t)\sim t^{-1}\ln t$.} Of course the small data global existence result
 in this case by applying \eqref{eq:ener-in-sem}-\eqref{eq:K-S-toy} {\it separately} to
each of the components of $\phi=(\phi_1,\phi_2)$. Recall however
that the system \eqref{eq:simple} models the Einstein equations only
after the latter is decomposed relative to its null frame components
$h_{LL}, h_{\Lb.\Lb}....$. Such decompositions do not commute with
the wave equation and thus prevent one from deriving separate energy
estimates for each of the null components of $h$. In this paper we
are able to solve this problem by adding another ingredient: {\it an
"independent" decay estimate}. The discussion below will be focused
on the tensor $h^1=g-m-h^0$ obtained from the original metric $g$ by
subtracting its "Schwarzschild part". We reluctantly allow for the
fact that the energy $E_N(t)$ of the tensor $h^1$ might be growing
with the rate of $\ve^2 t^\delta$ as $t\to\infty$ for some small
constant $\de>0$ dependent on the smallness of the initial data. The
Klainerman-Sobolev inequality \eqref{eq:K-S-toy} would then imply
that \beq\label{eq:decay-toy} \sup_{Z\in {\cal Z}, |I|\le N-3}|\pa
Z^I h^1(t,x)|\le C\ve (1+t+r)^{-1+\de} (1+|q|)^{-\frac 12} \eq In
order to close the argument, i.e., verify that the energy $E_N(t)$
indeed grows at the rate of at most $\ve^2 t^\de$ we must, according
to \eqref{eq:ener-sem}, upgrade the decay estimates
\eqref{eq:decay-toy} to the decay rate of $\ve t^{-1}$. To do that
we invoke the asymptotic system
\eqref{eq:Asym-C1}-\eqref{eq:Asym-C2}. It turns out that merely
assuming that the energy $E_N(t)\le \ve (1+t)^\de$ and consequently
\eqref{eq:decay-toy} allows us to show that the asymptotic system
provides an effective approximation of the full nonlinear system,
i.e., the discarded terms containing tangential derivatives do not
influence the asymptotic behavior of the field $h^1$. The crucial
property of the asymptotic system
\eqref{eq:Asym-C1}-\eqref{eq:Asym-C2} is that, as opposed to the
full nonlinear equation, it does commute with the null frame
decomposition of $h$ thus leading to the system
\eqref{eq:Asym1}-\eqref{eq:Asym-Phi}. Once the validity of the
asymptotic system is established the asymptotic behavior of its
solutions provides a new decay estimate and, as we described above,
shows that for ${\cal T}=\{L,S_1,S_2\},\, {\cal
U}=\{L,\Lb,S_1,S_2\}$
$$
|\pa h|_{\cal TU}\le C \ve t^{-1},\qquad
|\pa h|_{\cal \Lb\Lb}\le C\ve t^{-1}\ln t,\qquad |\pa \psi|\le C\ve t^{-1}
$$
The decay rate of $\ve t^{-1}\ln t$ can be potentially disastrous since
\eqref{eq:ener-sem} would imply that the energy $E_N(t)$ could grow
at the super-polynomial rate of $\ve^2\exp[C\ve \ln^2 t]$.
However the remarkable structure of the Einstein equations comes to the
rescue. The analysis of the semilinear terms $F_{\mu\nu}(\pa h,\pa h)$
in the equation
$$
\Boxr_g h^1_{\mu\nu}= F_{\mu\nu}(\pa h,\pa h) - \Boxr_g h^0_{\mu\nu}+2\pa_\mu\psi\,\pa_\nu\psi
$$
shows that the sharp $t^{-1}$ decay is only required for the $\pa h_{\cal TU}$
components
in order for the energy estimate $E_N(t)\le C\ve^2 t^\de$ to hold.
A glimpse of this structure has already appeared in our discussion of the weak
null condition. Recall that
$$
F_{\mu\nu}(\pa h,\pa h)=P(\pa_\mu h,\pa_\nu h) + Q_{\mu\nu}(\pa h,\pa h)
+ G_{\mu\nu}(h)(\pa h,\pa h).
$$
The cubic term $G$ and the quadratic form $Q$ satisfying the standard null condition
 are consistent with the uniform boundedness of the energy $E_N(t)$ and thus
 irrelevant for this discussion. The quadratic form $P(\pa_\mu h,\pa_\nu h) $
 can be decomposed relative to the null frame and obeys the following
 estimate
 $$
 |P(\pa_\mu h,\pa_\nu h) |\le |\pa h|_{\cal TU} |\pa h|_{\cal UU},
 $$
i.e., the most dangerous term $|\pa h|^2_{\cal \Lb\Lb}$, which would lead
to the damaging estimate
$$
E_0(t) \le \exp [C\int_0^t |\pa h(\tau)|_{\cal \Lb\Lb}] E_0(0),
$$
is absent!

The picture painted above is clearly somewhat simplified: \,\,we have only
indicated the argument implying that the lowest order energy
$E_0(t)\le C\ve^2 t^\de$ and we have not even begun to address the
effect of the quasilinear terms. In what follows we describe the building blocks
of our result and explain challenges of the quasilinear structure.

\noindent
{\bf  Energy inequality with weights}
\begin{equation}\label{eq:basicenergy}
 \int_{\Si_{T}} |\pa\phi|^{2} w + \int_{0}^{T} \int_{\Si_{t}}
|\pab\phi|^{2} w'\leq
 8\!\int_{\Si_{0}}\!\! |\pa \phi|^{2}w +  C\varepsilon\!
\int_0^T\!\!\int_{\Si_{t}}\! \frac {|\pa\phi|^{2}w}{1\!+\!t}
 +16\!\int_0^T \!\!\int_{\Sigma_t}\!\! |\Boxr_g\phi||\pa_t\phi|w,
 \end{equation}
 where $\pab$ denotes derivatives tangential to the outgoing
 Minkowski light cones $q=r-t$.
 The weight function
 $$
w(q)=\begin{cases} 1+(1+|q|)^{1+2\gamma} ,\qquad q\ge 0\\
                  1+(1+|q|)^{-2\mu} ,\qquad q<0 \end{cases}
$$
with $\mu\ge 0, \ga\ge -1/2$
serves a double purpose: \,\,it generates an additional {\it positive} space-time
integral giving an a priori control of the tangential derivatives $\pab\phi$
(we mention in passing that the use of such space-time norms leads to a very
simple proof of the small data global existence result for semilinear equations
$\Box\phi=Q(\pa\phi,\pa\phi)$ satisfying the standard null condition. To
our knowledge this argument has not appeared in the literature
before), it
provides the means to establish additional decay in $q$ via the Klainerman-Sobolev
inequality.
The energy estimate \eqref {eq:basicenergy} is established under very weak
general assumptions on the background metric $g^{\a\b}=m^{\a\b}+H^{\a\b}$,
e.g.,
$$
|\pa H|\le C\ve (1+t)^{-\frac 12} (1+|q|)^{-\frac 12}
(1+q_-)^{-\mu}, \qquad q_-=|\min (q,0)|,
$$
which are consistent with our expectations that if the energy
$E_N(t)\le \ve^2(1+t)^\de$ then the tensor $h$ and consequently
$H=-h +O(h^2)$ decay with a rate of at least $t^{-1+\de}$.
However, as is the case with virtually every estimate in this paper,
special stronger conditions are required for the $H_{LL}$ component
of the metric:
$$
|\pa H|_{LL}\le C\ve (1+t)^{-1} (1+|q|)^{-1}.
$$
Once again we note that such decay is consistent with the behavior predicted
by the asymptotic system for the Einstein equations, yet even better control
for this component is provided by the wave coordinate condition, which we will
turn to momentarily.

The interior estimate \eqref{eq:basicenergy}, i.e., with $w(q)\equiv 0$ for $q\ge 0$,
 in the constant coefficient case basically follows
 by averaging the energy estimates on light cones used e.g. in \cite{S}.
 We also note that the interior energy estimates with space-time quantities involving
 special derivatives of a solution were also considered and used in the work of Alinhac, see
 e.g. \cite{A2}, \cite{A3}). In \cite{L-R2} we proved the
 interior estimate \eqref{eq:basicenergy} under natural assumptions on the
 metric $g$ following in particular from the wave coordinate condition.
 The use of the weights $w(q)$ in the exterior $q\ge 0$ in energy estimates
 for the space part
 $\int_{\Si_T} |\pa\phi|^2 w$ originates in \cite{K-N2}.

Motivated by the energy inequality \eqref{eq:basicenergy} we define the
energy $E_N(t)$ associated with the tensor
$$
h^1_{\mu\nu}= g_{\mu\nu}-m_{\mu\nu}-h^0_{\mu\nu}, \qquad
h^0_{\mu\nu}=\chi(\frac rt)\chi (r) \frac Mr \de_{\mu\nu},
$$
solution of the reduced Einstein equations $\Boxr_g
h^1_{\mu\nu}=F_{\mu\nu} -\Boxr_g h^0+\pa_\mu\psi\,\pa_\nu\psi$ and a collection ${\cal
Z}=\{\pa_\a, \Omega_{\a\b}=-x_\a\pa_\b+x_\b\pa_\a, S=x^\a\pa_\a\}$
of commuting Minkowski vector fields, as follows:
$$
E_N(t) = \sum_{Z\in {\cal Z}, |I|\le N} \left (\|w^{\frac 12} \,\pa Z^I h^1(t,\cdot)\|_{L^2}^2+
\|w^{\frac 12} \,\pa Z^I \psi(t,\cdot)\|_{L^2}^2\right)
$$
\noindent
 {\bf The  Klainerman-Sobolev inequalities with weights}
\begin{equation*}
w^{1/2}\left (|\pa Z^I h^1(t,x)|+|\pa Z^I \psi(t,x)|\right)\leq \frac{C E_N(t)
}{(1\!+\!t\!+\!r)(1\!+\!|q|)^{1/2}},\qquad |I|\le N-2
\end{equation*}
is the fundamental tool allowing us to derive first preliminary, sometimes
referred to as weak, decay estimates.

\noindent
{\bf Decay estimate}

The following decay estimate
\begin{equation}\label{eq:decay-simple}
\varpi(q)(1\!+\!t)|\pa \phi(t,x)| \!\les\!
\!\!\int_0^t\! (1\!+\!\tau)\| \varpi\Boxr_g
\phi(\tau,\cdot) \|_{L^\infty} d\tau  +
\!\!\!\!\sum_{\,\,|I|\leq 2\!\!} \,\,\,\,\int_0^t
 \!\!\!\| \varpi Z^I
\phi(\tau,\cdot)\|_{L^\infty} \frac{d\tau}{1\!+\!\tau}
\end{equation}
is the additional "independent" estimate designed to boost the weak
decay estimates derived via the Klainerman-Sobolev inequality.
The weight function
$$
\varpi(q)=\begin{cases} 1+(1+|q|)^{1+\gamma'} ,\qquad q\ge 0\\
                  1+(1+|q|)^{\frac 12-\mu'} ,\qquad q<0 \end{cases}
$$
with $\gamma'\ge -1, \mu'\le 1/2$ is chosen in harmony with the
weight function $w(q)$ of the energy estimates. The unweighted
version of this estimate was used in \cite{L1} in the constant
coefficient case. We generalize this estimate to the variable
coefficient operator $\Boxr_g=g^{\a\b}\pa_\a\pa_\b$ under very weak
general assumptions on the metric $g^{\a\b}=m^{\a\b} + H^{\a\b}$,
e.g.,
$$
\int_0^t \|H(\tau,\cdot)\|_{L^\infty} \frac {d\tau}{1+\tau}\le \ve.
$$
However, once again the $H_{\cal LT}$ components are required to obey
the stronger condition
$$
|\pa H|_{\cal LT} \le C\ve (1+t)^{-1} (1+|q|)^{-1}
$$
The estimate above is closely connected with
the asymptotic equation in the sense that it
is obtained by treating the angular derivatives as lower order
and integrating the equation
$$
\pap\pa_q (r\phi)=r\Box
\phi+\tfrac{1}{r}\triangle_\omega \phi
$$
An important virtue of \eqref{eq:decay-simple} is that when
applied to systems $\Boxr_g \phi_{\mu\nu}= F_{\mu\nu}$
it can be derived for each {\it null} component of $\phi$.
This property is indispensable in view of the fact that the
weak null condition for the system of reduced Einstein equations
becomes transparent only after decomposing the tensor
$h_{\mu\nu}$ relative to a null frame
$\{L,\Lb,A,B\}$. In particular, \eqref{eq:decay-simple}
will lead to the estimates
$$
|\pa\psi|+|\pa h|_{\cal TU} \le C\ve (1+t+|q|)^{-1},
\qquad |\pa h|\le C\ve (1+t)^{-1}\ln (1+t).
$$
\noindent
{\bf Commutators}

As was mentioned above the combination: "energy estimate -
Klainerman-Sobolev inequality-decay estimate",  results in sharp
decay estimates for the first derivatives of the tensor $h^1$  and
merely polynomial growth of the energy $E_0(t)$. However, both the
Klainerman-Sobolev and the additional decay estimate require control
of the higher energy norms $E_N(t)$ involving the collection ${\cal
Z}$ of Minkowski vector fields. To obtain that control one needs to
commute vector fields $Z\in {\cal Z}$ through the equation $\Boxr_g
h^1 = F-\Boxr_g h^0$ and reapply the energy estimate
\eqref{eq:basicenergy}. This is by far the most difficult task in
dealing with quasilinear equations. The collection ${\cal Z}$ enjoys
good commutation properties with the wave operator $\Box$ of
Minkowski space: $[Z,\Box]=-c_Z\Box$, where $c_Z=2$ for $Z=S$ and is
zero otherwise, however its commutator with the wave operator
$\Boxr_g$ can create undesirable terms. It is exactly for this
reason that both in the work \cite{C-K} on stability of Minkowski
space and \cite{A3} on small data global existence for the equation
$\Box\phi=\phi\Delta\phi$ a different collection of {\it
geometrically modified} vector fields ${\cal Z}$ was used. The new
vector fields are adapted to the true characteristic surfaces of the
metric $g$ (this originates and is especially manifest in the
beautiful construction of \cite{C-K}) and better suited for
commuting through $\Boxr_g$, this construction however adds another
complicated layer to the proof.

In this work we employ the collection ${\cal Z}$ of original
Minkowski vector fields and argue that the bad commutation
properties of ${\cal Z}$ with $\Boxr_g$ can be improved if one takes
into account the wave coordinate condition \beq\label{eq:wave-intro}
\pa_\mu\Big (g^{\mu\nu} \sqrt{|\det g|}\Big )=0 \eq satisfied by the
metric $g$ in a coordinate system $\{x^\mu\}_{\mu=0,...,3}$. To
explain this phenomenon consider the commutator $[Z,\Boxr_g]$ with
one of the Minkowski vector fields $Z\in {\cal Z}$. For simplicity
assume that $Z\ne S$ so that $[Z,\Box]=0$. Then
$$
[Z,\Boxr_g]=[Z,H^{\a\b}\pa_\a\pa_\b]=\big (ZH^{\a\b} + H^{\a\b}_Z
\Big )\pa_\a\pa_\b,\qquad
H^{\a\b}_Z:=H^{\alpha\gamma} c_{\gamma}^{\,\,\, \beta}
+H^{\gamma\beta}c_{\gamma}^{\,\,\,\alpha}
$$
and the coefficients $c_{\a}^{\,\,\,\b}=\pa_\a Z^\beta$. It turns
out that since $Z$ is either Killing or conformally Killing vector
field ($Z=S$) vector field of Minkowski space the coefficients
$c_{\a}^{\,\,\,\b}$ have the property that $c_{L}^{\,\,\,\Lb}=0$.,
which implies that $|H^{\Lb\Lb}_Z|\le  |H|_{\cal LT}$. In accordance
with the usual arguments the worst term generated by the commutator
$[Z,\Boxr_g]$ contains two derivatives $\pa_q^2$ transversal to the
light cones $q=r-t$ and any modification to the vector fields $Z$
targets to eliminate such a term. In our case this term comes with a
coefficient $(ZH^{\Lb\Lb} +H^{\Lb\Lb}_Z)$ and thus can be estimated
by
$$
\frac{|ZH|_{\cal LL} + |H|_{\cal LT}}{1+|q|} |\pa Zh^1|.
$$
 Examining the equation
 $$
 \Boxr_g Zh^1= [\Boxr_g, Z] h^1 + ZF-Z\Boxr_g h^0+Z(\pa\psi\,\pa\psi)
 $$
 it is not too difficult to see that in order to control (i.e. establish a polynomial bound)
 the energy $E_1(t)$ one needs
 to at least show that
 $$
 |ZH|_{\cal LL} + |H|_{\cal LT}\le C\ve (1+t)^{-1} (1+|q|)
 $$
To explain the difficulty of such estimate we note that in the
region $|q|\le C$ this estimate is saturated by the Schwarzschild
part $\chi(r/t) \chi(r) M \delta_{\mu\nu}/r$ of $H$. While the
desired estimate for $H_{\cal LT}$ can be obtained from the sharp
decay estimates on the first derivatives of $h^1$, the estimate for
$ZH_{LL}$ is a much more subtle issue. An attempt to return to the
asymptotic system for the Einstein equations or alternatively the
decay estimate \eqref{eq:decay-simple} and show that $ZH_{LL}$ still
satisfies the sharp decay estimate requires once again a commutator
argument, this time in the context of the decay estimates, and
fails. In fact, we can only show that the higher $Z$ derivatives of
$h^1$ decay at the rate of $C\ve (1+t)^{-1+C\ve}$. It is at this
point that we recall that the tensor $g^{\a\b}=m^{\a\b} +H^{\a\b}$
verifies the wave coordinate condition \eqref{eq:wave-intro} which
after simple manipulations imply that \beq\label{eq:intro-H-wave}
\pa_q H_{\cal LT} = \pab H + O(H\cdot \pa H) \eq Here $\pab$ stands
for a tangential derivative, which in particular means that $|\pab
H|\le (1+t)^{-1} |ZH|$. Thus if we expect $H$ and its higher $Z$
derivatives to decay with the rate $t^{-1+\de}$ consistent with the
$t^\de$ energy growth and the Klainerman-Sobolev inequality we can
conclude that $\pa_q H_{\cal LT}$ will decay with the rate of at
least $t^{-2+\de}$. Integrating with respect to the $q=r-t$ variable
from the initial data at $t=0$ we can get the desired rate of decay
for $H_{\cal LT}$. We should note that in order to exploit the
additional decay in $t$ we need to keep track and use the decay in
$|q|$. This, in particular directly applies to the integration
described above. The necessity of utilizing the decay in $q$
explains our desire to work with weighted energy and decay
estimates. Amazingly \eqref{eq:intro-H-wave} is preserved after
commuting with a $Z$ vector field, but only for the $H_{LL}$
component, i.e.,
$$
\pa_q ZH_{\cal LL} = \pab ZH + Z O(H\cdot \pa H).
$$
Moreover, further commutations would destroy this precise structure.
Luckily, the control of higher energies $E_k(t)$ still only requires the estimate
$$
 |ZH|_{\cal LL} + |H|_{\cal LT}\le C\ve (1+t)^{-1} (1+|q|),
 $$
since the principal term in the commutator
$[Z^I,\Boxr_g]$ is still $(ZH^{\Lb\Lb} +H^{\Lb\Lb}_Z)\pa_q^2 Z^I$.
We should also mention that from the point of view of the energy estimates
the commutator $[Z^I,\Boxr_g] h^1$  contains
another dangerous term:
$$
Z^IH_{\cal LL} \pa^2_q h^1.
$$
Since most of the $Z$ derivatives fall on the fist factor one is forced to take it
in the energy norm and use the decay estimates for $\pa^2_q h^1$.
Note the two essential problems:\,\, the term $Z^IH_{\cal LL}$ does not
contain a derivative $\pa$, which is required in order to identify with the
energy norm, the decay of $\pa^2_q h^1$ falls short of the needed decay rate of $t^{-1}$. We resolve these problems as follows. First, a Hardy type inequality allows to convert a weighted $L^2$ norm of $Z^IH_{\cal LL}$
into a weighted $L^2$ norm of its derivatives. We then  use the
wave coordinate condition, which even after commutation of $Z^I$ preserves
some of its strength:
$$
|\pa Z^I|_{\cal LT} \le \sum_{|J|\le |I|} |\pab Z^I H| +
\sum_{|J|<|I|} |\pa Z^J H| +
Z^{I} O(H\cdot \pa H)
$$
The crucial fact here is that the principal term $\pab Z^I H$ contains
tangential derivatives $\pab$ and thus can be compared to the positive
space-time integral on the left hand-side of the energy inequality \eqref{eq:basicenergy},
which in this case will be
$$
\int_0^t \int_{\Si_\tau} |\pab Z^I h^1|^2 w',
$$
which means that we no longer have to be concerned about the decay of
 $\pa^2_q h^1$ in $t$ but rather trade it for the decay in $q$, also needed in the
Hardy type inequality, which we always have in abundance.
\vskip 1pc

\noindent
{\bf Scalar field}

The scalar field $\psi$ satisfies the wave equation $\tilde\Box_g\psi=0$ and 
contributes the quadratic term of the form $\pa_\mu\psi\,\pa_\nu\psi$ to the right
hand side of the equation for the tensor $h_{\mu\nu}$. Its asymptotic behavior 
is determined by the decay rate $|\pa\psi|\le C\epsilon (1+t+|q|)^{-1}$. All of these 
properties indicate that $\psi$ is very similar to the ``good" components 
$h_{\mathcal {TU}}$ of the tensor $h$. The only substantial difference, which in fact 
works to the advantage of $\psi$, is that unlike $h_{\mathcal {TU}}$ function $\psi$ 
satisfies its own wave equation and thus admits its own independent energy estimates.
This discussion indicates that in the problem of stability of Minkowski space in harmonic 
gauge, coupling 
the scalar field to the Einstein equations does not lead to fundamental changes in the 
structure of the equations and requires only superficial modification of our analysis in the
vacuum problem. As a consequence, starting from Section 9, the proof will be given only for
the vacuum case with $\psi\equiv 0$.

\vskip 2pc

We conclude this section by giving the plan of the
paper. For the sake of simplicity and clarity of the
exposition we will only deal with the case of a vanishing scalar field. The
generalization to the case of non-vanishing scalar field is at
each point immediate. Theorem \ref{globalexist} for the case of
vanishing scalar field is stated in section \ref{section:exist},
where the actual proof of the theorem starts.

{In \bf Section 3} we write the Einstein equations as a system of
quasilinear wave equations for the tensor $h_{\mu\nu}=g_{\mu\nu}-
m_{\mu\nu}$ relative to a system of wave coordinates
$\{x^\mu\}_{\mu=0,...,3}$.

In {\bf Section 4} we define the Minkowski null frame
$\{L,\Lb,S_1,S_2\}$ and describe null frame decompositions. We
estimate relevant tensorial quantities, including the $\Boxr_g$ and
the quadratic form $P$ appearing on the right hand side of the
reduced Einstein equations, in terms of their null decompositions.

{\bf Section 5} introduces the collection of Minkowski vector fields
${\cal Z}$ and records relations between the standard derivatives
$\pa$ and vector fields $Z$. It contains an important proposition
giving the estimate for the commutator between $\Boxr_g$ and the
powers $Z^I$. More details are contained  in the Appendix A.

In {\bf Section 6} we derive our basic weighted energy estimate.

{\bf Section 7} deals with the decay estimates for solutions of an inhomogeneous
wave equation $\Boxr_g \phi =F$. The key result is contained in Corollary \ref{decaywaveeq3}.

In {\bf Section 8} we discuss the wave coordinate condition. In particular,
we state the estimates for the $H_{\cal LL}$ and $H_{\cal LT}$ components
of the tensor $H^{\a\b}=g^{\a\b}-m^{\a\b}$ and its $Z$ derivatives.
Some of the details are provided in Appendix D.

{\bf Sections 9-11} contain the statement and the proof of the small data global
existence result for the system of reduced Einstein equations.

In {\bf section 9} we state our main result and set up our inductive argument.
We assume that the energy $\E_N(t)$ of the tensor $h^1_{\mu\nu}=h_{\mu\nu}-h_{\mu\nu}^0$  obeys the estimate  $E_N(t)\le C\ve^2 (1+t)^\de$ and
on its basis derive the weak decay estimates for $h$ and $h^1$ and estimate
the inhomogeneous terms $F_{\mu\nu}$, $F^0_{\mu\nu}=\Boxr_g h^0_{\mu\nu}$
appearing on the right hand side of the equation for $h^1$.
The weak decay estimates are obtained by means of the Klainerman-Sobolev
inequality proved in Appendix C.

In {\bf section 10} we use the decay estimates derived in section 7 to upgrade
the weak decay estimate of section 9.

{\bf Section 11} uses the energy estimate of section 6 to verify the inductive
assumption on the energy $\E_N(t)$. The results of section 11.3 heavily rely on
Hardy type inequalities established in Appendix B.

\section {The Einstein equations in wave
coordinates}\label{section:einstwavec}

For a Lorentzian metric  $g$ and a system of coordinates
$\{x^\mu\}_{\mu=0,..,3}$ we
denote
\beq \label{eq:christof} \Gamma_{\mu\,\,\, \nu}^{\,\,\,
\lambda} =\frac{1}{2} g^{\lambda\delta} \big( \pa_\mu
g_{\delta\nu} + \pa_\nu g_{\delta\mu}-\pa_\delta g_{\mu\nu}\big),
\eq
the Christoffel symbols of $g$ with respect to $\{x^\mu\}$.
We recall that
\beq \label{eq:ctensor}
R_{\mu\,\,\, \nu\delta }^{\,\,\, \lambda}= \pa_\delta
\Gamma_{\mu\,\,\, \nu}^{\,\,\, \lambda} -\pa_\nu\Gamma_{\mu\,\,\,
\delta}^{\,\,\, \lambda} +\Gamma_{\rho\,\,\, \delta}^{\,\,\,
\lambda} \Gamma_{\mu\,\,\, \nu}^{\,\,\, \rho} -\Gamma_{\rho\,\,\,
\nu}^{\,\,\, \lambda} \Gamma_{\mu\,\,\, \delta}^{\,\,\, \rho} \eq
is the Riemann curvature tensor of $g$ and  $R_{\mu\nu}=R_{\mu\,\,\,
\nu\alpha}^{\,\,\, \alpha}$ is the Ricci tensor.

We assume that the metric $g$ together with a field $\psi$ is a solution of the Einstein-scalar field
equations
\beq \label{eq:Einst} R_{\mu\nu}=\pa_\mu\psi\,\pa_\nu\psi.
\eq
A system of coordinates $\{x^\mu\}$ is called the wave coordinates
iff the Cristoffel symbols $\Gamma$ verify the condition:
\beq \label{eq:wavecord}
 g^{\alpha\beta}\, \Gamma_{\alpha\,\,\, \beta}^{\,\,\, \lambda}=0,\qquad
 \forall \lambda=0,..,3
\eq
Each of the following three equations is equivalent to \eqref{eq:wavecord}, which we will refer to as the {\it wave coordinate condition},
\beq \label{eq:wavec2}
\pa_\a (g^{\a\b} \sqrt{|g|})=0,\qquad
g^{\alpha\beta}\pa_\alpha
g_{\beta\mu} =\frac{1}{2}g^{\alpha\beta}\pa_\mu
g_{\alpha\beta}, \qquad \pa_\alpha
g^{\alpha\nu}=\frac{1}{2} g_{\alpha\beta} g^{\nu\mu}
 \pa_\mu g^{\alpha\beta}
\eq
The  covariant wave operator $\Box_g$
 in wave coordinates, in view of
\eqref{eq:wavec2}, coincides with the reduced wave operator $\Boxr_g=g^{\a\b}\pa^2_{\a\b}$,
\beq \label{eq:boxg}
\Boxr_g={\square}_g=\frac{1}{\sqrt{|g|}}\pa_\alpha
 g^{\alpha\beta}\sqrt{|g|}\pa_\beta
\eq
Define a 2-tensor $h$ according to the decomposition
$$
g_{\mu\nu}=m_{\mu\nu}+h_{\mu\nu},
$$
where $m$ is the standard Minkowski metric on $\R^{3+1}$.
Let $m^{\mu\nu}$ be the inverse of $m_{\mu\nu}$. Then for small
$h$
$$
H^{\mu\nu}:=g^{\mu\nu}-m^{\mu\nu}=-h^{\mu\nu}+O^{\mu\nu}(h^2),
\qquad\text{where}\qquad
h^{\mu\nu}=m^{\mu\mu^\prime}m^{\nu\nu^\prime}h_{\mu^\prime\nu^\prime}
$$
and $O^{\mu\nu}(h^2)$ vanishes to second order at $h=0$. Recall that
according to our conventions the indices of tensors $h_{\mu\nu},
H^{\mu\nu}$ are raised/lowered with respect to the metric $m$. The
following proposition was established in the vacuum case in \cite{L-R1} and amounts to
a rather tedious calculation.
\begin{proposition} \label{Einstwavecquad}
Let $(g,\psi)$ be a solution of the Einstein-scalar field  equations then relative to
the wave coordinates the tensor $h$ and $\psi$ solve a system of wave equations
\beq \label{Eired}
{\Boxr}_g h_{\mu\nu}=F_{\mu\nu}(h)(\pa h,\pa h) + 2\pa_\mu\psi\,\pa_\nu\psi,\qquad {\Boxr}_g\psi=0\eq
The inhomogeneous term $F_{\mu\nu}$ has the following structure:
\begin{align}
&F_{\mu\nu}(h)(\pa h,\pa h)=P(\pa_\mu h,\pa_\nu
h)+Q_{\mu\nu}(\pa h,\pa h) +G_{\mu\nu}(h)(\pa h,\pa h) \label{eq:defF} \\
&P(\pa_\mu h,\pa_\nu h)=\frac{1}{4}
m^{\alpha\alpha^\prime}\pa_\mu h_{\alpha\alpha^\prime} \,
m^{\beta\beta^\prime}\pa_\nu h_{\beta\beta^\prime}  -\frac{1}{2}
m^{\alpha\alpha^\prime}m^{\beta\beta^\prime} \pa_\mu
h_{\alpha\beta}\, \pa_\nu h_{\alpha^\prime\beta^\prime} \label{eq:defP},
\end{align}
quadratic form $Q_{\mu\nu}(U,V)$ for each $(\mu,\nu)$ is a linear
combination of the standard null forms $ m^{\a\b} \pa_\a U \pa_\b V$
and $\pa_\a U\pa_\b V -\pa_\b U \pa_\a V$ and $G_{\mu\nu}(h)(\pa
h,\pa h)$ is a quadratic form in $\pa h$ with coefficients smoothly
dependent on $h$ and vanishing for $h=0$: $G_{\mu\nu}(0)(\pa h,\pa
h)=0$, i,e., $G$ is a cubic term.
\end{proposition}
Observe that the quadratic terms in \eqref{eq:defP} do not satisfy the
classical null condition. However  the trace $\tr h= m^{\mu\nu}
h_{\mu\nu}$ satisfies a nonlinear wave equation with semilinear
terms obeying the  the null condition:
$$
g^{\alpha\beta}\pa_\alpha\pa_\beta m^{\mu\nu} h_{\mu\nu} =Q(\pa
h,\pa h)+G(h)(\pa h,\pa h)+Q(\pa\psi,\pa\psi).
$$

\section {The null-frame}\label{section:null}
 At each point $(t,x)$ we introduce a
pair of null vectors $(L,\Lb)$
$$
L^0=1, \quad L^i=x^i/|x|, \,\,\, \, i=1,2,3,\quad\text{and}\quad
\underline {L}^0=1,\quad \underline {L}^i=-x^i/|x|, \,\,\, \,
i=1,2,3.
$$
Let  $S_1$ and $S_2$ be two othonormal smooth tangent vectors fields
to the sphere $\bold{S}^2$.  $S_1,S_2$ are then orthogonal to the
normal $\omega=x/|x|$ to $\bold{S}^2$ and $(L,\Lb, S_1,S_2)$ form a
nullframe.

\begin{remark} The null frame described above is defined only locally. Replacing
orthonormal vector fields $S_1, S_2$ with the projections
\begin{equation}\label{eq:deftan}
\pas_i\!=\pa_i-\omega_i\omega^j\pa_j,\qquad \omega=\frac x{|x|}
\end{equation}
of the standard coordinate vector fields $\pa_i$ would define a
global frame. Set $\pab_0\!=\! L^\alpha\pa_\alpha$ and
$\pab_i=\pas_i$, for $i=1,2,3$. Then $\{\pab_0,...,\pab_3\}$ span
the tangent space of the forward light cone.
\end{remark}

We raise and lower indices with respect to the Minkowski metric;
$X_\alpha=m_{\alpha\beta} X^\beta$.  For a vector field/one form $X$
we define its components relative to a null frame according to
$$
X^\alpha=X^L L^\alpha+ X^{\underline{L}}\underline{L}^\alpha+X^A
A^\alpha,\qquad X_\alpha=X^L L_\alpha+
X^{\underline{L}}\underline{L}_\alpha+X^A A_\alpha.
$$
Here and in what follows $A,B,C..$ denotes any of the vectors $S_1,
S_2$, and we used the summation convention;
 $$
X^A A^\alpha= X^{S_1} S_1^\alpha+X^{S_2} S_2^\alpha
 $$
 The components can be calculated
using the following formulas
$$
X^{L}= -\frac 12 X_{\Lb},\quad X^{\Lb}=-\frac 12 X_{L},\quad
X^A=X_A,
$$
where
 $$
 X_Y=X_\alpha Y^\alpha.
 $$
 Similarly for a two form $\pi$ and two vector fields $X$ and $Y$ we define
 $$
 \pi_{XY}=\pi_{\alpha\beta} X^\alpha Y^\beta.
 $$
 The Minkowski metric $m$ has the following form relative to a
null frame
$$
m_{LL}=m_{\Lb\Lb}=m_{LA}= m_{\Lb A}=0, \qquad m_{L\Lb}=m_{\Lb
L}=-2, \qquad m_{AB}=\delta_{AB},
$$
i.e. $m_{\alpha\beta} X^\alpha Y^\beta=-2(X^L
Y^{\underline{L}}+X^{\underline{L}}Y^L) +\delta_{AB} X^A Y^B$, where
$\delta_{AB} X^A Y^B=X^{S_1} Y^{S_1}+X^{S_2}Y^{S_2}$. The inverse of
the metric has the form
$$
m^{LL}=m^{\Lb\Lb}=m^{LA}= m^{\Lb A}=0, \qquad m^{L\Lb}=m^{\Lb
L}=-1/2,
 \qquad m^{AB}=\delta^{AB}
$$

We also define the tangential trace \beq
\overline{\operatorname{tr}}\,k=\delta^{AB}k_{AB},
 \eq and record
the identity
$$
\pi^{\alpha\beta}\pa_\alpha =\pi^{L\beta} L+\pi^{\underline{L}\beta}
{\underline{L}} +\pi^{A\beta} A
=-\frac{1}{2}\pi_{\underline{L}}^{\,\,\,\beta} L -\frac{1}{2}
\pi_L^{\,\,\,\beta} {\underline{L}} +\pi_A^{\,\,\,\beta} A.
$$

We introduce the following notation. Let ${\cal T}=\{L,S_1,S_2\}$,
${\cal U}=\{\Lb,L,S_1,S_2\}$, ${\cal L}=\{L\}$ and ${\cal
S}=\{S_1,S_2\}$. For any two of these families ${\cal V}$ and ${\cal
W}$ and an arbitrary two-tensor $p$ we define
\begin{align}
|p|_{{\cal V W}} &=\sum_{V\in{\cal V},W\in{\cal W} , }
|p_{\beta\gamma} V^\beta W^\gamma|,\label{eq:norm1}\\
|\pa p|_{{\cal V W}} &=\sum_{U\in {\cal U},V\in{\cal V},W\in{\cal W}
, } |(\pa_\alpha p_{\beta\gamma})U^\alpha V^\beta W^\gamma|,
\label{eq:norm2}\\
|\overline{\pa} p|_{{\cal V W}} &=\sum_{T\in {\cal T},V\in{\cal
V},W\in{\cal W} , } |(\pa_\alpha p_{\beta\gamma}) T^\alpha V^\beta
W^\gamma| \label{eq:norm3}
\end{align}
Note that the contractions with the frame is outside the
differentiation so they are not differentiated.

 Let $Q$ be one of the
null quadratic form, i.e. $Q_{\alpha\beta}(\pa \phi,\pa \psi) =
\pa_\alpha  \phi\, \pa_\beta \psi-\pa_\beta \phi\, \pa_\alpha \psi$
if $\alpha\neq \beta$ and $Q_{0}(\pa \phi,\pa \psi)=m^{\alpha\beta}
\pa_\alpha  \phi\, \pa_\beta \psi$. Motivated by \eqref{eq:defP} we
define the quadratic form $P$ \beq\label{eq:P-quad}
P(\pi,\theta)=\frac 12 m^{\a\a'} m^{\b\b'}
\pi_{\a\b}\theta_{\a'\b'}- \frac 14 m^{\a\b} m^{\a'\b'}
\pi_{\a\b}\theta_{\a'\b'} \eq The proof of the following result is a
simple exercise and we leave it to the reader.
\begin{lemma}\label{tander}
For an arbitrary 2-tensor $\pi$ and a scalar function $\phi$
\begin{align}
&|\pi^{\alpha\beta}
\pa_\alpha\phi \, \pa_\beta \phi|\les
\big(|\pi|_{{\cal L}{\cal L}}|\pa\phi|^2+|\pi|\,
|\pab \phi||\pa\phi|\big),\label{eq:dernullframeZ0} \\
&|L_\alpha \pi^{\alpha\beta} \pa_\beta\phi
|\les \big(|\pi|_{{\cal L}{\cal L}}|\pa\phi|+|\pi|\,|\pab \phi|\big),\label{eq:dernullframeZ1} \\
&|(\pa_\alpha \pi^{\alpha\beta})
\pa_\beta\phi |\les \big(|\pa \pi|_{{\cal L}{\cal L}}+|\overline{\pa}
\pi|\big)\,|\pa\phi| +|\pa \pi|\, |\overline{\pa} \phi| \label{eq:dernullframeZ2}\\
&|\pi^{\alpha\beta} \pa_\alpha
\pa_\beta \phi|\les \big(|\pi|_{\cal LL}|\pa^2\phi|+|\pi|\,|\pab\pa
\phi|\big) \label{eq:dernullframeZ}
\end{align}
For the quadratic form  $P$ defined in \eqref{eq:P-quad} and a null form $Q$
\begin{align}
&|P(\pi,\theta)|\les |\pi\,|_{\cal
TU}|\theta|_{\cal TU} +|\pi\,|_{\cal LL}|\theta|+|\pi\,|\,|\theta|_{\cal LL},
\label{eq:tanP} \\
&|Q(\pa \phi,\pa \psi)|\les |\overline{\pa} \phi| |\pa \phi| + |\pa
\phi||\overline{\pa} \psi| \label{eq:nullformtan}
\end{align}
\end{lemma}

For a function $\phi$, $|\pab \phi|=\sum_{\alpha=0}^3
|\pab_\alpha\phi|$, where $\pab_\alpha$ the tangential derivatives
defined after \eqref{eq:deftan}, is equivalent to $\sum_{T\in\,{\cal
T}}| T^\alpha\pa_\alpha \phi|$.
 However
 $|\pab^2 \phi|=\sum_{\alpha,\beta=0}^3
 |\pab_\alpha\pab_\beta\phi|$
 is not equivalent to $\sum_{S,\,T\in \,{\cal T}} |T^\alpha
S^\beta\pa_\alpha\pa_\beta \phi|$, the difference is of order
$|\pa\phi|/r$ which is a main term. Furthermore $|\pab^2\phi|$ need
not even be equivalent to $\sum_{S,T\in{\cal T}}
|T^\alpha\pa_\alpha(S^\beta\pa_\beta \phi)|$, the difference is
however bounded by $|\pab \phi|/r$ which is a lower order term.

\begin{lemma}\label{tander2}
For an arbitrary symmetric 2-tensor $\pi^{\alpha\beta}$ and a
function $\phi$ \beq \label{eq:nuZ} \big|\pi^{\alpha\beta}
\pa_\alpha\pa_\beta \phi- \pi_{ LL}\pa_{q}^2\phi - 2 \pi_{ L\Lb}
\pa_s \pa_{q} \phi -r^{-1}\,\overline{\operatorname{tr} } \,\pi\,
\pa_q \phi\big|\les |\pi|_{\cal LT}|\pab\pa\phi|+|\pi|\,\big(|\pab^2
\phi| + r^{-1} |\bar\pa\phi|\big). \eq
\end{lemma}
\begin{proof}
Using that $\pa_q, \pa_s$ derivatives commute with the frame
$\{L,\Lb, S_1, S_2\}$ we obtain that
$$
\pi^{\alpha\beta} \pa_\alpha\pa_\beta \phi- \pi_{LL}\pa_{q}^2\phi -
2 \pi_{L\Lb} \pa_s \pa_{q} \phi= \pi_{\Lb\Lb} \pa_s^2 \phi +
2\pi_{\Lb A} \pa_s (A^i\pa_i \phi)
 -\pi_{L A}  {\Lb}^\alpha A^i\pa_i \pa_\alpha\phi +
\pi_{AB} A^i B^j \pa_i \pa_j \phi
$$
Furthemore,
$$
\pi_{AB} A^i B^j \pa_i \pa_j \, \phi=\pi_{AB} A^i \pa_i (B^j \pa_j
\,\phi ) - \pi_{AB} (A^i \pa_i B^j) \pa_j \, \phi
$$
Decomposing with respect to the null frame we obtain
$$
 (A^i \pa_i B^j) \pa_j = (A^i \pa_i B)^L L +  (A^i \pa_i B)^{\Lb}\, {\Lb} +
  (A^i \pa_i B)^C C
$$
Note that $|A^i \pa_i B^j|\le C r^{-1}$, since $B^j$ are smooth
functions of $\omega=x/|x|\in\bold{S}^2$. Now, since
$B^j\omega_j=0$;
$$
2(A^i \pa_i B)^{\Lb}=-(A^i\pa_i B^j)\,\omega_j
 = A^i B^j\pa_i \,\omega_j=A^i B^j \frac{1}{r}(\delta_{ij}-\omega_i\,\omega_j)
 =\frac{1}{r}\delta_{ij} \, A^i B^j
=\frac 1r \de_{AB}
$$
is precisely the null second fundamental form of the outgoing null
cone $t-r=$const. The lemma now follows since as pointed out before
the lemma
$$
\sum_{S,\, T\in\, {\cal T}} |T^\alpha
\pa_\alpha(S^\beta\pa_\beta\phi)|\leq
C\big(|\pab^2\phi|+r^{-1}|\pab\phi|\big).
$$
\end{proof}
\begin{cor}\label{cor:BoxH}
Let $\phi$ be a solution of the reduced wave equation $\Boxr_g\phi
=F$ with a metric $g$ such that $H^{\a\b}=g^{\a\b}-m^{\a\b}$
satisfies the condition that $|H^{L\Lb}|<\frac 14$. Then \beq
\label{eq:phiwave}
 \Big|\Big(4\pa_s -
\frac{H_{LL}}{2g^{L\Lb}}\pa_q-\frac{\overline{\operatorname{tr}}\, H
+H_{L\Lb}}{2g^{L\underline{L}}\,\, r} \Big) \pa_q
(r\phi)+\frac{rF}{2g^{L\underline{L}}}\Big|\les r|\triangle_\omega
\phi|+ |H|_{L{\cal T}}\left (r\, |\pab\pa\phi|+|\pa\phi|\right)+ |H|\,\big(r\,|\pab^2
\phi| + |\bar\pa\phi|+r^{-1}|\phi|\big) \eq where
$\triangle_\omega=\bar\triangle=\delta^{ij}\bar\pa_i \bar\pa_j$.
\end{cor}
\begin{proof}
Define the new metric
$$
\tilde g^{\a\b} = \frac {g^{\a\b}}{-2 g^{L\Lb}}.
$$
The equation $g^{\a\b}\pa_{\a}\pa_{\b}\phi = F$ then takes the
form
$$
\tilde g^{\a\b} \pa_{\a}\pa_{\b}\phi = \frac {F}{-2g^{L\Lb}},
$$
which also can be written as
$$
\Box \phi + (\tilde g^{\a\b}-m^{\a\b}) \pa_{\a}\pa_{\b}\phi =
\frac {F}{-2g^{L\Lb}}
$$
Let $\pi^{\a\b}$ be the tensor $\pi^{\a\b} = (\tilde
g^{\a\b}-m^{\a\b})$ Observe that
\begin{align*}
\pi^{\a\b} = (-2g^{L\Lb})^{-1} \big(g^{\a\b} + 2 m^{\a\b}
g^{L\Lb}\big ) & = (-2g^{L\Lb})^{-1} \big ( H^{\a\b} + m^{\a\b}
(2g^{L\Lb} + 1)\big )\\ &= (-2g^{L\Lb})^{-1} \big ( H^{\a\b} + 2
m^{\a\b} H^{L\Lb} \big )
\end{align*}
Thus,
\begin{equation}
\label{eq:Lbdisappear} \pi_{L\Lb}=0,\qquad \quad \pi_{L{\cal T}}=
(-2g^{L\Lb})^{-1} H_{L{\cal T}}, \qquad \quad
\overline{\operatorname{tr}}\, \pi \, =(-2g^{L\Lb})^{-1}
\big(\overline{\operatorname{tr}}\,H \,+H_{L\Lb}\big)
\end{equation}
Moreover, $|\pi|\les |H|$, since $g^{L\Lb}= H^{L\Lb} - \frac 12 $
and by the assumptions of the Corollary $|H^{L\Lb}|<\frac 14$.

Now using \eqref{eq:nuZ} of Lemma \ref{tander2}, with the
condition that $\pi_{L\Lb}=0$, together with the decomposition
$$
\Box \phi = -\pa_t^2 \phi + \triangle \phi = \frac 1r
(\pa_t+\pa_r)(\pa_r-\pa_t) r \phi + \triangle_\omega\phi= \frac 4r
\pa_s \pa_q r \phi + \triangle_\omega\phi.
$$
we find that the identity $\Box \phi + \pi^{\a\b}
\pa_{\a}\pa_{\b}\phi = ({-2g^{L\Lb}})^{-1} F$ leads to the
inequality
$$
\big|4 \pa_{s}\pa_{q} r \phi + r \pi_{LL}\pa_{q}^2\phi
+\overline{\operatorname{tr}}\, \pi
\,\pa_q\phi+(2g^{L\underline{L}})^{-1} r F\big|\les
r|\triangle_\omega \phi| + r |\pi|_{L{\cal T}}|\pab\pa\phi|+
|\pi|\,\big(r\, |\pab^2 \phi| +|\bar\pa\phi|\big)
$$
Finally, identity \eqref{eq:Lbdisappear} and a crude  estimate
$|\pi|\les |H|$ yield the desired result.
\end{proof}

\section {Vector fields}\label{section:commute}
The family of vector fields
$${\cal Z}=\{\pa_\a,
\Omega_{\a\b} = -x_{\a}\pa_{\b} + x_{\b}\pa_{\a},
 S = t\pa_{t} + r\pa_{r}\}
 $$
 plays a special role in the study of the wave
 equation in Minkowski space-time. We denote the above vector fields by  $Z^\iota$
with an 11-dimensional integer index $\iota=(0,...1,..,0)$. Let
$I=(\iota_1,...,\iota_k)$, where $|\iota_i|=1$, be a multi-index of
length $|I|=k$ and let $Z^I=Z^{\iota_1}\cdot\cdot\cdot Z^{\iota_k}$
denote a product of $k$ vector fields from the family ${\cal Z}$. By
a sum $I_1+I_2=I$ we mean a sum over all possible order preserving
partitions of the multi-index $I$ into two multi-indices $I_1$ and
$I_2$, i.e., if $I=(\iota_1,...,\iota_k)$, then
$I_1=(\iota_{i_1},...,\iota_{i_n})$ and
$I_2=(\iota_{i_{n+1}},...,\iota_{i_k})$, where $i_1,...,i_k$ is any
reordering of the integers $1,...,k$ such that $i_1<...<i_n$ and
$i_{n+1}<...<i_k$.
 With this convention
Leibnitz rule becomes $Z^I(fg)=\sum_{I_1+I_2=I} (Z^{I_1}
f)(Z^{I_2} g)$.

We recall that the family {\cal Z} possesses  special commutation
properties ${\cal Z}$: for any vector field $Z\in {\cal Z}$
$$
[Z,\Box]=-c_Z \Box,
$$
where the constant $c_Z$ is only different from zero in the case of
the scaling vector field $c_S=2$. We also record the following
expressions for the coordinate vector fields:
\begin{align}
&\pa_{t} = \frac {tS - x^{i}\Omega_{0i}}{t^{2}-r^{2}},
\label{eq:tZ}\\
&\pa_{r} = \omega^i\pa_i =\frac {t\omega^{i}\Omega_{0i} -
rS}{t^{2}-r^{2}},\label{eq:rZ}\\
&\pa_{i} = \frac {-x^{j}\Omega_{ij} + t \Omega_{0i} - x_{i}
S}{t^{2}-r^{2}}= -\frac {x_{i} S} {t^{2}-r^{2}} + \frac {x_{i}
x^{j}\Omega_{0j}}{t(t^{2}-r^{2})} + \frac {\Omega_{0i}}t
\label{eq:iZ}
\end{align}
In particular,
\begin{equation}
\label{eq:somega} \pab_0=\pa_{s} =\frac{1}{2}\big(\pa_t+\pa_r\big) =
\frac {S + \omega^{i}\Omega_{0i}}{2(t+r)},\qquad \pab_{i}=
\pa_{i}- \omega_{i}\pa_{r} =\frac{\omega^j\Omega_{ij}}{r} = \frac
{-\omega_{i}\omega^{j}\Omega_{0j} + \Omega_{0i}}t.
\end{equation}
Recall that $\bar\pa$ denotes the tangential derivatives, i.e.,
$\text{Span}\{\bar\pa_0,\bar\pa_1,\bar\pa_2,\bar\pa_3\} =
\text{Span}\{\pa_{s}, {S_1},{S_2}\}$.
\begin{lemma}\label{tanderZ}
For any function $\phi$ and a symmetric 2-tensor
$\pi$
\begin{align}
&(1+t+|q|)|\bar \pa  \phi|+(1+|q|)|\pa \phi|\les C
\sum_{|I|=1}|Z^I \phi|,\label{eq:tanZ}\\
&|\bar\pa^2 \phi|+r^{-1}|\pab \phi|\les \frac{C}{r}\sum_{|I|\leq 2}
\frac{|Z^I \phi|}{1+t+|q|}, \quad{\text {where}} \quad |\bar\pa^2
\phi|^2=\sum_{\alpha,\beta=0,1,2,3}
|\bar\pa_\alpha\bar\pa_\beta \phi|^2, \label{eq:2tanZ} \\
&|\pi^{\alpha\beta} \pa_\alpha\pa_\beta
\phi|\leq C\bigg (\frac{|\pi|}{1+t+|q|}+ \frac{|\pi|_{\cal LL}}{1+|q|}\bigg
) \sum_{|I|\le 1}|\pa Z^I \phi| \label{eq:derframeZ}
\end{align}
\end{lemma}
\begin{proof} First we note that if $r+t\leq 1$ then
\eqref{eq:tanZ} holds since the standard derivatives $\pa_\alpha$ are
included in the sum on the right. The inequality for $|\bar\pa \phi|$
in \eqref{eq:tanZ} follows directly from \eqref{eq:somega}.
The inequality for $|\pa f|$ in \eqref{eq:tanZ} follows from
\eqref{eq:tZ} and the first identity in \eqref{eq:iZ}.

The proof of \eqref{eq:2tanZ} follows immediately from  \eqref{eq:somega} and the inequality $|\pa_i\omega_j|\leq C
r^{-1}$.
The inequality \eqref{eq:derframeZ} follows from Lemma
\ref{tander2}, \eqref{eq:tanZ} and the commutator identity
$[Z,\pa_i]=c_i^\alpha\pa_\a$, which holds with constants
$c_i^\a$.
\end{proof}
Next we state the result following from the above lemma and Corollary \ref{cor:BoxH}.
\begin{lemma}\label{waveeqframe} Let $\phi$ be a solution of the
equation $\Boxr_g \phi=F$ with $g^{\a\b}=m^{\a\b} + H^{\a\b}$.
Then \beq\label{eq:waveoper}
 \Big|\Big(4\pa_s
-\frac{H_{LL}}{2g^{L\Lb}}\pa_q -\frac{\overline{\operatorname{tr}}\,
H \,+H_{L\Lb}}{2g^{L\underline{L}}\,\, r} \Big) \pa_q
(r\phi)+\frac{rF}{2g^{L\underline{L}}}\Big| \les \Big(1+
\frac{r\,|H|_{\cal LT}}{1+|q|}+ |H|\Big) r^{-1} \sum_{|I|\le 2}|Z^I
\phi|
 \eq
\end{lemma}
\begin{proof} By Corollary \ref{cor:BoxH}
\begin{multline*}
\Big|\Big(4\pa_s -\frac{H_{LL}}{2g^{L\Lb}}\pa_q
-\frac{\overline{\operatorname{tr}}\, H \,+H_{L\Lb}}
 {2g^{L\underline{L}}\,\, r}\Big) \pa_q
(r\phi)+\frac{r F}{2g^{L\underline{L}}}\Big|  \\
\les
r|\triangle_\omega \phi|+  |H|_{L{\cal T}}\left (r |\pab\pa\phi|+|\pa\phi|\right )+ |H|\,\big(
r\,|\pab^2 \phi| + |\bar\pa\phi|+r^{-1}|\phi|\big)
\end{multline*}
where $\triangle_\omega=\delta^{ij}\bar\pa_i \bar\pa_j$. Here all
the the derivatives can be reexpressed in terms of the vector
fields $Z$ and $\pa_q$ using \ref{tanderZ}, yielding the
expression \eqref{eq:waveoper}. Note that
$$
|\bar\pa \pa\phi| \les \frac{\sum_{|I|=1}|Z^I\pa \phi|}{1+t+|q|}
\les \frac{\sum_{|I|\leq 1} |\pa Z^I \phi|}{1+t+|q|} \les
\frac{\sum_{|I|\leq 2} |Z^I\phi|}{(1+|q|)(1+t+|q|)}.
$$
\end{proof}
The last result of this section records an important statement
concerning commutation between the reduced wave operator $\Boxr_g$
and the family of Minkowski vector fields ${\cal Z}$. As we already
explained in Section 2 our small data data global existence result
for the system of reduced Einstein equations $\Boxr_g h^1= F-\Boxr_g
h^0$ is based on controlling the energy and pointwise norms of the
quantities $Z^I h^1$ with vector fields $Z\in {\cal Z}$ \-- the
family of {\it Minkowski vector  fields}. The above control is
achieved via  the energy and decay estimates for solutions of the
inhomogeneous wave equation
$$
\Boxr_g Z^I h^1=\hat Z^I \Boxr_g h^1 - \Boxr_g Z^I h^1+
\hat Z^I  F- \hat Z^I\Boxr_g h^0,\qquad \hat Z=Z+c_Z
$$
and therefore requires good estimates on the commutator
$\hat Z^I \Boxr_g  - \Boxr_g Z^I$. To underscore the importance
of this commutator we recall that the only two known examples of
quasilinear hyperbolic systems, with the metric dependent on the solution
rather than its derivatives, possessing small data global solutions
 required the use of {\it modified} vector fields primarily due to the lack
 of good estimates for the commutator, \cite{A3}, \cite{C-K}.
 The proof of the proposition below together with other commutator
 related statements is contained in Appendix A.
\begin{proposition}\label{prop:commut3}
Let $\Boxr_g=\Box + H^{\a\b} \pa_\a\pa_\b$. Then for any vector
field $Z\in{\cal Z}$ we have with $\hat Z=Z+c_Z$
\begin{align}
&|\Boxr_g Z^I \phi-\hat{Z}^I \Boxr_g \phi| \les \frac 1{1+t+|q|}
\,\,\,\sum_{|K|\leq |I|,}\,\, \sum_{|J|+(|K|-1)_+\le |I|} \,\,\,
|Z^{J} H|\,\, {|\pa Z^{K} \phi|}
\label{eq:curvwaveeqcommutest2-5} \\
+ \frac 1{1+|q|}&
 \sum_{|K|\leq |I|}\Big(\sum_{|J|+(|K|-1)_+\leq |I|} \!\!\!\!\!|Z^{J} H|_{LL}
+\!\!\!\!\!\sum_{|J^{\prime}|+(|K|-1)_+\leq
|I|-1}\!\!\!\!\!|Z^{J^{\prime}} H|_{L\cal T}
+\!\!\!\!\!\sum_{|J^{\prime\prime}|+(|K|-1)_+\leq |I|-2}\!\!\!\!\!
|Z^{J^{\prime\prime}} H|\Big) {|\pa Z^{K} \phi|}\nn
\end{align}
where $(|K|-1)_+=|K|-1$ if $|K|\geq 1$ and $(|K|-1)_+=0$ if
$|K|=0$.
\end{proposition}

\section{Energy estimates in curved space-time}\label{section:energywave}
In this section we establish basic weighted energy identities and estimates for solutions
of the inhomogeneous wave equation
\beq
\label{eq:quasinh}
\Boxr_{g} \phi =F
\eq
The weights under consideration will depend only on the distance $q=r-t$
to the light cone $t=r$ and will be defined by a function $w$
$$
w=w(q)\geq 0 ,\quad w^\prime(q)\geq 0
$$
The weights will serve a two-fold purpose: 1) to provide additional decay
of $\phi$ in the exterior region $q\ge 0$, which will follow from the energy
estimates via a Klainerman-Sobolev type inequality\,\, 2) to establish an
additional a priori bound on a {\bf space-time} integral involving
tangential derivatives $\pab\phi$ of the solution.
\begin{lemma}
\label{lemma:Energy}
Let $\phi$ be a solution of the equation \eqref{eq:quasinh}
decaying sufficiently fast as $|x|\to\infty$. Assume that the background metric
$g$ is such that the tensor $H^{\a\b}=g^{\a\b}-m^{\a\b}$ satisfies
$|H|\le \frac 12$. Then with $\omega=x/|x|$
\begin{align}
 \int_{\Si_{t_{2}}} \big (|\pa_{t}\phi|^{2}
+|\nab \phi|^{2}\big ) w(q) dx
&+ 2\int_{t_{1}}^{t_{2}} \int_{\Si_{\tau}}
|\pab\phi|^{2} w^\prime(q)dx dt \le
 4 \int_{\Si_{t_{1}}} \big (|\pa_{t}\phi|^{2}
+|\nab \phi|^{2}\big )w(q)\, dx  \label{eq:ER}\\ &+ 2
\int_{t_{1}}^{t_{2}} \int_{\Si_{\tau}} \big |(2\pa_{\a} H^{\a\b})\pa_{\b}\phi
\pa_{t}\phi -(\pa_{t} H^{\a\b})\pa_{\a}\phi\pa_{\b}\phi
+ 2F\pa_{t}\phi \big | w(q)  dx dt
\label{eq:energy}\\ &+ 2\int_{t_{1}}^{t_{2}} \int_{\Si_{\tau}}
\big |
H^{\a\b}\pa_{\a}\phi\pa_{\b}\phi +
2( \omega_i H^{i\b}-H^{0\b})\pa_{\b}\phi\pa_{t}\phi \big | w^\prime(q)\, dx dt \nn
\end{align}
\end{lemma}
\begin{proof} Let $\phi_i=\pa_i \phi$, $i=1,2,3$, and $\phi_t=\pa_t\phi$.
Differentiating under the
integral sign and integrating by parts we get
\begin{multline}
\frac{d}{dt} \int\big(-g^{00}\phi_t^2+g^{ij}\phi_i\phi_j\big) w(q)\,
dx-\int 2\pa_j\big(g^{0j}\phi_t^2 w(q)\big)\, dx
=2\int w(q) \big(-g^{00} \phi_t\phi_{tt}+g^{ij}\phi_i\phi_{tj}-2g^{0j}\phi_t\phi_{tj}\big)\, dx \\
+\,\int  w(q) \big(-(\pa_t g^{00}) \phi_t^2+(\pa_t
g^{ij})\phi_i\phi_{j}-2(\pa_j g^{0j})\phi_t^2\big)\, -w^\prime(q)
\big(-g^{00}\phi_t^2+g^{ij}\phi_i\phi_j
+2\omega_j g^{0j}\phi_t^2\big)\, dx\\
=2 \int w(q) \big(-g^{00} \phi_t\phi_{tt}-g^{ij}\phi_t\phi_{ij}-2g^{0j}\phi_t\phi_{tj}\big)\, dx \\
+\,\int  w(q) \big(-(\pa_t g^{00}) \phi_t^2+(\pa_t g^{ij})\phi_i\phi_{j}-2(\pa_j g^{0j})\phi_t^2-2(\pa_i g^{ij})\phi_t\phi_j\big)\, dx\\
-\int w^\prime(q) \big(-g^{00}\phi_t^2+g^{ij}\phi_i\phi_j +2\omega_j
g^{0j}\phi_t^2+2\omega_i g^{ij}\phi_t \phi_j\big)\, dx
\end{multline}
Hence,
\begin{multline}
\frac{d}{dt} \int\big(-g^{00}\phi_t^2+g^{ij}\phi_i\phi_j\big) w(q)\,
dx =-\int w(q)\big( 2\phi_t \,\square_g \phi-(\pa_t g^{\alpha\beta})
\phi_\alpha \phi_\beta
+2(\pa_\alpha g^{\alpha\beta})\phi_\beta \phi_t\big)\, dx \\
-\int w^\prime(q) \big(g^{\alpha\beta}\phi_\alpha\phi_\beta
+2\big(\omega_i g^{i\alpha} -g^{0\alpha}\big)\phi_t \phi_\alpha\, dx
\end{multline}
Furthermore, with $\phi_r=\omega^i\phi_i=\pa_r \phi$ and $\overline{\phi}_i=\phi_i-\omega_i\phi_r=\overline{\pa}_i\phi$
$$
m^{\alpha\beta}\phi_\alpha\phi_\beta+2\phi_t(\omega_i
m^{i\alpha}-m^{0\alpha})
\phi_\alpha=-\phi_t^2+\delta^{ij}\phi_i\phi_j
+2\phi_t(\omega^i\phi_i+\phi_t)
=(\phi_t+\phi_r)^2+\delta^{ij}\overline{\phi}_i\overline{\phi}_j=|\pab\phi|^2
$$
Since $|H|<1/2$ we also have that
$$
\frac{1}{2} (\phi_t^2+\delta^{ij}\phi_i\phi_j)
\leq -g^{00}\phi_t^2+g^{ij}\phi_i\phi_j\leq 2(\phi_t^2+\delta^{ij}\phi_i\phi_j)
$$
The lemma follows.
\end{proof}

We now consider the following weight function:
\begin{equation}\label{eq:energyweight}
w=w(q)=\begin{cases} 1+(1+|q|)^{1+2\gamma},\quad\text{when }\quad q>0\\
      1+(1+|q|)^{-2\mu}\,\quad\text{when }\quad q<0\end{cases}
\end{equation}
for some $\mu\geq 0$, $\gamma\geq -1$.
We clearly have
 \beq\label{eq:w-prime}
w^\prime\leq 4w(1+|q|)^{-1}\leq 16\gamma^{-1} w^\prime(1+q_-)^{2\mu},
\eq
where $q_-=|q|$ for $q\le 0$ and $q_-=0$ otherwise.
\begin{proposition}
\label{prop:Decayenergy}
Let $\phi$ be a solution of the wave equation \eqref{eq:quasinh}
with the metric $g$ such that for $H^{\alpha\beta}=g^{\alpha\beta}-m^{\alpha\beta}$;
\begin{align}
&(1+|q|)^{-1} |H|_{LL} +|\pa H|_{LL}+|\overline{\pa} H|\leq
C\varepsilon' (1+t)^{-1},\nn\\
& (1+|q|)^{-1}\,|H|+ |\pa H|\leq C\varepsilon' (1+t)^{-\frac 12}
(1+|q|)^{-\frac 12}(1+q_-)^{-\mu}\label{eq:metricdecay}
\end{align}
Then for any $0<\ga\le 1$, and $0<\varepsilon'\leq \gamma/C_1$, we have
\beq\label{eq:firstenergy}
 \int_{\Si_{t}} |\pa\phi|^{2}\,w + \int_{0}^{t} \int_{\Si_{\tau}}
|\pab\phi|^{2}\,w^{\,\prime} \leq
 8\int_{\Si_{0}} |\pa \phi|^{2}\,w+
16\int_0^t\int_{\Si_{t}} \Big(\frac{C\varepsilon
\,|\pa\phi|^{2}}{1+t}
 + |F|\,|\pa\phi|\Big)\, w
\eq
\end{proposition}
\begin{remark}
Observe that by the Gronwall inequality the energy estimate of the
above proposition implies $t^{\varepsilon'}$ growth of the energy.
\end{remark}
\begin{remark}
We recall again that the interior estimate \eqref{eq:firstenergy}, i.e., with $w(q)\equiv 0$ for $q\ge 0$,
 in the constant coefficient case basically follows
 by averaging the energy estimates on light cones used e.g. in \cite{S1}.
 The interior energy estimates with space-time quantities involving
 special derivatives of a solution were also considered and used in the work of Alinhac, see
 e.g. \cite{A2}, \cite{A3}). In \cite{L-R2} we proved the
 interior version of \eqref{eq:firstenergy}.
 The use of the weights $w(q)$ in energy estimates in the exterior
 $q\ge 0$ for the space part
 $\int_{\Si_T} |\pa\phi|^2 w$ originates in \cite{K-N2}.
\end{remark}
\begin{proof}
The proof of the proposition relies on the energy estimate obtained
in Lemma \ref{lemma:Energy}.
We decompose the terms on the right hand-side of \eqref{eq:ER}
with respect to the null frame, using Lemma \ref{tander}
$$
|(\pa_{\a} H^{\a\b})\pa_{\b}\phi \pa_{t} \phi|+
|(\pa_{t} H^{\a\b})\pa_{\a}\phi \pa_{\b}\phi |\le
\big(|(\pa H)_{LL}| + |\pab H|\big ) \, |\pa \phi|^2 +
|\pa H|\,|\pab \phi|\,|\pa\phi|
$$

Therefore, using the assumptions \eqref{eq:metricdecay} on the metric
$g$, we obtain that
\beq\label{eq:spaceq}
|2(\pa_{\a} H^{\a\b})\pa_{\b}\phi \pa_{t} \phi -
(\pa_{t} H^{\a\b})\pa_{\a} \phi \pa_{\b} \phi | \les
\frac {\varepsilon'}{1+t} |\pa\phi|^2 + \frac {\varepsilon'}{(1+|q|)(1+q_-)^{2\mu}}
|\pab\phi|^2
\eq
Decomposing the remaining terms we infer that
$$
|H^{\a\b}\pa_{\a}\phi \pa_{\b} \phi |
+|L_\alpha H^{\a\b}(\pa_t\phi)\pa_{\b} \phi |
\le
\big |H|_{LL} |\pa\phi|^2 +
|H| |\pab\phi|\,|\pa\phi|
$$
Once again, using the assumptions \eqref{eq:metricdecay}, we have
\beq\label{eq:timeq}
|2 H^{\a\b} L_\a \pa_{\b} \phi \pa_{t} \phi +
H^{\a\b}\pa_{\a}\phi \pa_{\b} \phi| \les
\varepsilon' \frac {1+|q|}{1+t} |\pa\phi|^2 +\varepsilon' \frac {|\pab\phi|^2}
{(1+q_-)^{2\mu}}
\eq
Thus, with the help of \eqref{eq:w-prime},
$$
\int_{\Si_{t}}|\pa\phi|^2 w  + \int_{0}^{t} \int_{\Si_{\tau}}
|\pab\phi|^{2} w^\prime\les \int_{\Si_0}|\pa\phi|^2 w +
\varepsilon
\int_0^t \int_{\Si_\tau}\Big ( \frac {|\pa\phi|^2} {1+t} w +
|\pab\phi|^{2}\frac{w^\prime}{\gamma}\Big ) + \int_0^t \int_{\Si_\tau}
|F|\, |\pa\phi|w
$$
and the desired estimate follows.
\end{proof}

\section
{Decay estimates in curved space-time  } \label{section:decaywaveeq}
In this section we derive $L^\infty$ estimates for the first
derivatives of solutions of the equation \beq\label{eq:waveeqscal}
\Boxr_g \phi=F \eq These estimates are complimentary to the global
Sobolev inequalities derived in Appendix C and will provide a way to
improve upon the decay estimates derived from the weighted energy
estimates via global Sobolev inequalities. Estimates of these type
were obtained in \cite{L1} in the case of the wave equation in
Minkowski space-time, i.e., constant coefficient metric $g$, and
\cite{L-R1} for variable coefficients. Here, however, we introduce a
weighted version of the $L^\infty$ estimates deigned to capture
additional decay in $|q|=|r-t|$.

For $\gamma'\ge -1$, $\mu'\le 1/2$ define the weight
\begin{equation}\label{eq:decayweight}
\varpi=\varpi(q)=\begin{cases}
(1+|q|)^{1+\gamma^\prime},\quad\text{when }\quad q>0\\
      (1+|q|)^{1/2-\mu^\prime}\,\quad\text{when }\quad
      q<0\end{cases}.
\end{equation}
\begin{lemma} \label{decaywaveeq2} Let  $\phi$ be a solution of the
reduced scalar wave equation \eqref{eq:waveeqscal} on a curved
background with metric $g$. Assume that the tensor
$H^{\alpha\beta}=g^{\alpha\beta}-m^{\alpha\beta}$ obeys the estimates
\beq\label{eq:decaymetric1} |H|\leq \varepsilon^\prime, \qquad
\int_0^\infty \|H(t,\cdot)\|_{L^\infty(D_t)}\frac {dt}{1+t}\le \frac 14,\qquad |H|_{L{\cal T}}\leq
\varepsilon^\prime\frac{|q|+1}{1+t+|x|}, \eq
in the region $D_t=\{x:\,t/2<|x|<2t\}$. Then for
$\alpha=\max (1+\gamma',1/2-\mu')$
\begin{multline}\label{eq:decaywaveeq1}
(1+t)\varpi(q)\,|\pa \phi(t,x)| \leq  C\sup_{0\leq
\tau\leq t} \sum_{|I|\leq 1}\!
\|\varpi(q) \,Z^I\! \phi(\tau,\cdot)\|_{L^\infty}\\
+ C\int_0^t  \Big(\varepsilon^\prime \alpha\|\varpi(q)\,\pa
\phi(\tau,\cdot)\|_{L^\infty} +(1+\tau)\| \varpi(q)\,
F(\tau,\cdot)\|_{L^\infty(D_\tau)} +\sum_{|I|\leq 2} (1+\tau)^{-1}
\| \varpi(q)\, Z^I \phi(\tau,\cdot)\|_{L^\infty(D_\tau)}\Big)\, d\tau
\end{multline}
\end{lemma}
\begin{proof} Since by Lemma \ref{tanderZ}
\beq\label{eq:awaycone} (1+|t-r|) |\pa\phi|+(1+t+r)|\bar\pa\phi|
\leq C\sum_{|I|=1}|Z^I\phi|,\qquad r=|x|, \eq the inequality
\eqref{eq:decaywaveeq1} holds when $r<t/2+1/2$ or $r>2t-1$.

 By Lemma
\ref{waveeqframe} with $\pa_s=1/2(\pa_t+\pa_r)$ and
$\pa_q=1/2(\pa_r-\pa_t)$, \beq\label{eq:waveoper5} \big|(4\pa_s
-\frac{H_{LL}}{2g^{L\Lb}}\pa_{q})\pa_{q}( r \phi)\big| \les \Big(1+
\frac{r\,|H|_{\cal LT}}{1+|q|}+ |H|\Big) r^{-1}\sum_{|I|\le 2} |Z^I
\phi|+ |H|\, r^{-1}\, |\pa_q (r\phi) | + r |F| \eq Multiplying by
the weight $\varpi(q)$ and using that $\varpi'(q)\le C
\varpi(q)/(1+|q|)$ along with the assumptions
\eqref{eq:decaymetric1} we obtain \beq\label{eq:waveoper6}
\big|(4\pa_s -\frac{H_{LL}}{2g^{L\Lb}}\pa_{q})\varpi(q)\pa_{q}( r
\phi)\big| \les \Big (\frac{|H|}{1+t} +\alpha\frac
{H_{LL}}{1+|q|}\Big ) \varpi(q)|\pa_q (r\phi) |+\sum_{|I|\le 2}
\frac{\varpi(q)|Z^I \phi|}{1+t} + C (t+1)\varpi(q)|F|\eq in the
region $t/2+1/2<r<2t-1$. Let $(\tau,x(\tau))$ be the integral curve
of the vector field $\pa_s+H^{\Lb \Lb}(2g^{L\Lb})^{-1}\pa_q$ passing
through a given point $(t,x)$ contained in the region $t/2+1/2\leq
r\leq 2t-1$. Observe that, by the smallness assumption on $H$, any
such curve has to intersect the boundary of the set $t/2+1/2\leq
r\leq 2t-1$  at  $(\tau,y)$ such that $|y|=\tau/2+1/2$ or $|y|=2\,\tau-1$.

Then along such a curve the function
$\psi:=\varpi(q)\pa_q(r\phi)$ satisfies the following equation:
\beq \Big|\frac{ d}{dt}\psi\Big|\leq \hat{h}
|\psi|+ f
\eq
where
$$
\hat{h}=C\frac {|H|}{1+t},\qquad
f= \frac {\varepsilon' \alpha}{1+t} \varpi(q) |\pa_q r\phi|+
 C (1+t) \varpi(q)|F|+C\sum_{|I|\le 2} \frac{\varpi(q)|Z^I \phi|}{1+t}.
 $$
Thus using
the integrating factor $e^{-\hat{H}}$ with $\hat{H}=\int
\hat{h}(s)\, ds $ and  integrating
along the integral curve $(\tau,x(\tau))$ from any point $(t,x)$ in the
set $t/2+1/2\leq r\leq 2t-1$ to the first point of intersection $(t_0,x_0)$
with the boundary of the set $t/2+1/2\leq r\leq 2t-1$
we obtain
$$
|\psi(t,x)| \leq \exp\Big({\int_{\tau}^t
\|\hat{h}(\sigma,\cdot)\|_{L^\infty}\, d\sigma}\Big)
|\psi(t_0,x_0)| +\int_{\tau}^t \exp\Big({\int_{\tau^\prime}^t
\|\hat{h}(\sigma,\cdot)\|_{L^\infty}\,
 d\sigma}\Big)
\|f(\tau^\prime,\cdot)\|_{L^\infty}\, d\tau^\prime,
$$
with the $L^\infty$ norms are taken over the set $t/1+1/2\leq r\leq 2t-1$.

For the points $(t_0,x_0)$ such that
 $|x_0|=t_0/2+1/2$ or $|x_0|=2t_0-1$
we have by \eqref{eq:awaycone} that
$$
|\psi(t_0,x_0)|\leq Cr \varpi(q) |\pa_q\phi|+C \varpi(q)|\phi|\leq C\sum_{|I|\leq
1} \varpi(q) |Z^I\phi|.
$$
The desired inequality now follows from
\eqref{eq:decaymetric1}, which implies that
$\int_{0}^\infty
\|\hat{h}(\sigma,\cdot)\|_{L^\infty}\, d\sigma \leq \frac{1}{4}$
and
the inequality
$$
(1+t+r) |\pa\phi|\leq C\sum_{|I|\leq
1}|Z^I\phi|+C |\pa_q (r\phi)|.
$$
\end{proof}
We now state similar estimates for a system
\beq\label{eq:waveeqsyst}
\widetilde{\square}\phi_{\mu\nu}=F_{\mu\nu} \eq While it is trivial
to extend the estimates of Lemma \ref{decaywaveeq2} to each of the
components of $\phi_{\mu\nu}$ our interest lies in the estimates for
the {\bf null} components of $\phi$. Contracting
\eqref{eq:waveeqsyst} with the vector fields $\{L,\Lb,A,B\}$ is far
from straightforward. We instead exploit that the null derivatives
$\pa_s, \pa_q$ commute with any of the vector fields of the null
frame. We assume that the weight function $\varpi(q)$ is as in
\eqref {eq:decayweight} and $\alpha=\max (1+\gamma',1/2-\mu')$.
\begin{cor} \label{decaywaveeq3}
Let $\phi_{\mu\nu}$ be a solution of the
reduced wave equation system \eqref{eq:waveeqsyst} on a curved
background with a metric $g$. Assume that
$H^{\alpha\beta}=g^{\alpha\beta}-m^{\alpha\beta}$ satisfies
\beq\label{eq:decaymetric4} |H|\leq
\frac{\varepsilon^\prime}{4},\qquad\int_0^{\infty}\!\!
\|\,H(t,\cdot)\|_{L^\infty(D_t)}\frac{dt}{1+t} \leq
\frac{\varepsilon^\prime}{4},\qquad |H|_{L{\cal
T}} \leq\frac{\varepsilon^\prime}{4}\, \frac{|q|+1}{1+t+|x|}
\eq
in the region  $D_t=\{x\in\bold{R}^3;\,
t/2\leq |x|\leq 2t\}$.
 Then for any $U,V\in\{L,\Lb,A,B\}$ and an arbitrary point $x\in D_t$:
\begin{multline}\label{eq:decaywaveeq4}
(1+t+|x|)|\varpi(q)\pa \phi(t,x)|_{UV} \les\!\sup_{0\leq
\tau\leq t}
\sum_{|I|\leq 1}\|\,\varpi(q)Z^I\! \phi(\tau,\cdot)\|_{L^\infty}\\
+ \int_0^t\Big( \varepsilon^\prime \alpha\|\varpi(q) |\pa
\phi(t,\cdot)|_{UV}\|_{L^\infty} +(1+\tau)\| \varpi(q)
|F(\tau,\cdot)|_{UV}\|_{L^\infty(D_\tau)} +\sum_{|I|\leq 2} (1+\tau)^{-1}
\| \varpi(q)Z^I \phi(\tau,\cdot)\|_{L^\infty(D_\tau)}\Big)\, d\tau.
\end{multline}
\end{cor}
\begin{proof} By Lemma \ref{waveeqframe} for each component
we have the estimate \beq \Big|\Big(4\pa_s -
\frac{H_{LL}}{2g^{L\Lb}}\pa_q -\frac{\overline{\operatorname{tr}}\,
H \,+H_{L\underline{L}}}{2g^{L\underline{L}}\,\, r} \Big) \pa_q
(r\phi_{\mu\nu})+\frac{rF_{\mu\nu}}{2g^{L\underline{L}}}\Big| \les
\Big(1+ \frac{r\,|H|_{\cal LT}}{1+|q|}+ |H|\Big) r^{-1} \sum_{|I|\le
2}|Z^I \phi_{\mu\nu}| \eq and since $\pa_s$ and $\pa_q$ commute with
contraction with the frame vectors we get \beq \Big|\big(4\pa_s -
\frac{H_{LL}}{2g^{L\Lb}}\pa_q -\frac{\overline{\operatorname{tr}}\,
H \,+H_{L\underline{L}}}{2g^{L\underline{L}}\,\, r} \Big) \pa_q
(r\phi_{UV})+\frac{rF_{UV}}{2g^{L\underline{L}}}\Big| \les \Big(1+
\frac{r\,|H|_{\cal LT}}{1+|q|}+ |H|\Big) r^{-1} \sum_{|I|\le 2}|Z^I
\phi| \eq
 The proof  now proceeds as in Lemma \ref{decaywaveeq2}.
\end{proof}

\section{The wave coordinate condition} \label{section:decaywavec}
The results of previous sections underscore the special role played by
the $H_{LL}$ components of the tensor $H^{\a\b} =g^{\a\b}-m^{\a\b}$
in the energy and decay estimates for solutions of the wave equation
$\Boxr_g\phi=F$.
In this section we explain how the wave coordinate condition on the
tensor $g_{\mu\nu}$  provides additional information about $H_{LL}$.

Recall that the wave coordinate condition for metric $g$ in a coordinate
system $\{x^\mu\}_{\mu=0,...,3}$ takes the form
\beq\label{eq:wavecdef}
\pa_\mu\Big( g^{\mu\nu}
\sqrt{|\det{g}|}\Big)=0.
\eq
Expressing $g^{\mu\nu}$ in terms of the tensor $H^{\mu\nu}$
we obtain
$$
g^{\mu\nu} \sqrt{|\det{g}|}
=\big(m^{\mu\nu}+H^{\mu\nu}\big)\big(1-\frac{1}{2}\tr
H+O(H^2)\big)
$$
Therefore,
\beq\label{eq:wavecappr}
 \pa_\mu\Big( H^{\mu\nu}-\frac{1}{2}
m^{\mu\nu}\operatorname{tr}\, H+O^{\mu\nu}(H^2)\Big)=0 , \qquad
\text{where}\quad O^{\mu\nu}(H^2)=O(|H|^2).\eq Recall also the
family of, tangent to the outgoing Minkowski light cones, vector
fields ${\cal T}=\{L,A,B\}$.

The divergence of a vector field can be expressed relative to the
null frame as follows: \beq\label{eq:divnull} \pa_\mu
F^\mu=L_\mu\pa_q F^\mu -\underline{L}_\mu \pa_s F^\mu +A_\mu \pa_A
F^\mu \eq We can now easily prove
\begin{lemma} \label{decaywavec0} Assume that
$|H|\leq 1/4$. Then
\beq \label{eq:decaywavec1} |\pa H|_{\cal LT}
\les |\overline{\pa } H| + |H|\, |\pa H|
\eq
\end{lemma}
\begin{proof} It follows from \eqref{eq:wavecappr} and
\eqref{eq:divnull} that \beq
\big|L_\mu\pa\big(H^{\mu\nu}-\frac{1}{2}m^{\mu\nu}
\operatorname{tr}\, H\big)\big|\leq |\bar\pa H|+|H||\pa H| \eq
Contracting with $T\in {\cal T}$ and using that $m_{TL}=0$ gives the
desired result.
\end{proof}
We now state a generalization of the above result containing
estimates for the quantities $Z^I H_{LT}$ with vector fields $Z\in
{\cal Z}$\- our family of Minkowski Killing and conformally Killing
vector fields. The result is a rather tedious consequence of
commuting vector fields $Z$ through the wave coordinate condition
\eqref{eq:wavecdef} and we postpone the details of the proof until
Appendix D.
\begin{prop}\label{decaywavecZ}
Let $g$ be a Lorentzian metric satisfying the wave coordinate condition
 \eqref{eq:wavecdef} relative to a coordinate system $\{x^\mu\}_{\mu=0,...,3}$.
 Let $I$ be a multi-index and assume that the tensor $H^{\mu\nu}= g^{\mu\nu} -m^{\mu\nu}$ verifies the condition
$$
|Z^J H|\le C,\qquad \forall |J|\le |I|/2,\quad \forall Z\in {\cal Z}.
$$
Then for some constant $C'$
 \begin{align}
&|\pa Z^I H|_{ L\cal T}\le C'\Big (\sum_{|J|\leq |I|}|\overline{\pa} Z^J
H|+\!\!\!\!\sum_{|J|\leq |I|-1}\!\!\!|\pa Z^J H|\,\, +
 \sum_{\,\,\,\,|I_1|+||I_2| \leq |I|}|Z^{I_{2}} H| |\pa Z^{I_1} H|\Big )\label{eq:decaywavec6}\\
& |\pa Z^I H|_{ LL}\les  C'\Big (\sum_{|J|\leq |I|}|\overline{\pa} Z^J H|+
\sum_{|J|\leq |I|-2} |\pa Z^{J} H|
+ \sum_{|I_1|+|I_2| \leq |I|,\, m\geq
2} |Z^{I_{2}} H| |\pa Z^{I_1} H|\Big )\label{eq:decaywavec5}.
\end{align}
Similar estimates hold for the tensor $h_{\mu\nu}=g_{\mu\nu}-m_{\mu\nu}$.
\end{prop}

\section {Statement of the Main Theorem and beginning of the proof}\label{section:exist}
We consider the initial data $(g_{\mu\nu}|_{t=0}, \pa_t g_{\mu\nu}|_{t=0})$
for the Einstein-scalar field equations $R_{\mu\nu}=\pa_\mu\psi\,\pa_\nu\psi$, constructed in \eqref{eq:init-red1} -
\eqref{eq:init-red3}.
The spatial part ${g_{0}}_{ij}$ of $g_{\mu\nu}|_{t=0}$ together with the second fundamental form ${k_0}_{ij}=-1/2 \pa_t g_{ij}|_{t=0}$ satisfy the constraint equations
$$
R_0-|k_0|^2 + (\tr k_0)^2=\psi_1^2+|\nabla\psi_0|^2,\qquad
\nab^j {k_0}_{ij} -\nab_i \tr k_0=\psi_1\, \pa_i\psi_0.
$$
By construction the initial data  $(g_{\mu\nu}|_{t=0}, \pa_t g_{\mu\nu}|_{t=0})$
satisfies the wave coordinate condition
$$
\pa_\mu (g^{\mu\nu} \sqrt{|\det g|})=0.
$$
The reduced Einstein equations, written relative to a 2-tensor
$h_{\mu\nu}=g_{\mu\nu}-m_{\mu\nu}$ take the form
\begin{align}
&\Boxr_g h_{\mu\nu} = F_{\mu\nu}+2 \pa_\mu\psi\,\pa_\nu\psi,\qquad {\Boxr}_g\psi=0,\label{eq:Red-E}\\
& F_{\mu\nu}(h)(\pa h,\pa h)=P(\pa_\mu h,\pa_\nu h)+Q_{\mu\nu}(\pa
h,\pa h)
+G_{\mu\nu}(h)(\pa h,\pa h)\label{eq:F-prime},\\
& P(\pa_\mu h,\pa_\nu h):= \frac 14 \pa_{\mu} \tr h\, \pa_{\nu}\tr h-\frac 12 \pa_{\mu} h^{\a\b}\pa_{\nu}
h_{\a\b}, \label{eq:P-prime},
\end{align}
where $Q$ is a linear combination of standard quadratic null forms
and $G$ is a quadratic term in $\pa h$ vanishing in $h$ at $h=0$. We
assume that the initial data for $h(t)$ are such that at $t=0$,
$g_{\mu\nu} = m_{\mu\nu} + h_{\mu\nu}$ verifies the constraint
equations and the wave coordinate condition.

Define 2-tensor
$h^1_{\mu\nu}(t)$ \beq\label{eq:h1-h0}
h^1_{\mu\nu}(t):=h_{\mu\nu}(t)-h^0_{\mu\nu}(t),\qquad
\text{where}\quad h^0_{\mu\nu}(t)=\chi(r/t)\chi(r)
\frac{M}{r}\delta_{\mu\nu} \eq where $\chi(s)\in C^\infty$ is $1$
when $s\geq 3/4$ and $0$ when $s\leq 1/2$.  It is clear that we can
reinterpret \eqref{eq:Red-E} as the equation for a new unknown \--
2-tensor $h^1$: \beq\label{eq:Red-E1} \Boxr_g h^1_{\mu\nu} =
F_{\mu\nu} - \Boxr_g h^0_{\mu\nu} +2 \pa_\mu\psi\,\pa_\nu\psi\eq Set
\begin{equation}
\label{eq:enerNdef}
\E_N(t)=\sum_{|I|\leq N, Z\in {\cal Z}}
\left (\|w^{1/2}\pa Z^I h^1(t,\cdot)\|_{L^2}+\|w^{1/2}\pa Z^I \psi(t,\cdot)\|_{L^2}\right),
\end{equation}
{where}
$$
w=\begin{cases} 1+(1+|q|)^{1+2\gamma},\quad q>0\\
         1+(1+|q|)^{-2\mu},\quad q<0\end{cases}
$$
where $q=r-t$ and $\gamma,\mu>0$. We recall that the initial data for $\psi$ is given 
by $\psi|_{t=0}=\psi_0$ and $\pa_t\psi|_{t=0}= a \psi_1$ with $a^2=(1-\chi(r) M r^{-1})$.

We now state the main result.
\begin{theorem} \label{exist}
There exist a constant
$\varepsilon_0>0$ such that if $\varepsilon\leq \varepsilon_0$ and the
initial data $(h^1|_{t=0}, \pa_t h^1|_{t=0},\psi_0,\psi_1)$ are smooth, and obey
$\,\,\,\E_N(0)+M\leq \varepsilon\,\,\,$ together with the condition
\beq\label{eq:initialdecay}
\liminf\limits_{|x|\to \infty}\left(|h^1(0,x)|+|\psi_0(x)|\right)\to 0,
\eq
then the solution $h(t)$ of the reduced Einstein equations \eqref{eq:Red-E}
can be extended to a global smooth solution
 satisfying
\begin{equation}
\label{eq:enerNest} \E_N(t)\leq C_N \varepsilon\,
(1+t)^{c\varepsilon}
\end{equation}
where $C_N$ is a constant depending only on $N$ and $c$ is
independent $\ve$.
\end{theorem}
\begin{remark}\label{rem:Remark}
Recall that if the initial data for $(g_{\mu\nu},\psi)$ verifies the constraint
equations and the wave coordinate condition then there exists a local in time
classical solution $(h_{\mu\nu}(t),\psi(t))$ of \eqref{eq:Red-E} such that
$g_{\mu\nu}(t) = m_{\mu\nu}(t) + h_{\mu\nu}(t)$ obeys the wave coordinate
condition for any time $t$ in a maximum interval of existence.
\vskip 1pc
{\bf Therefore, in what follows we shall assume that a local in time solution
$g_{\mu\nu}(t)=m_{\mu\nu} + h_{\mu\nu}(t)$ obeys the wave coordinate
condition
\beq\label{eq:Wave-Let}
\pa_\mu \Big (g^{\mu\nu} (t) \sqrt{|\det g(t)|}\Big )=0
\eq
for any $0\le t\le T_0$.}
\vskip 1pc
We also note that
the maximum time of existence $T_0$ can be characterized by the blow-up of
the energy $\E_N(t)\to \infty$ as $t\to T_0^-$.
\end{remark}
\begin{remark}
As was explained in Section 2 the proof of the result for the Einstein-scalar field problem 
requires only superficial modifications as compared to the vacuum case.
\vskip 1pc
{\bf Therefore, in what follows all the arguments will be provided only for the vacuum problem with $\psi\equiv 0$.}
\end{remark}

For the proof we let $\delta$ be any fixed number $0<\delta<1/4$
and $\delta< \gamma$. We define the time $T<T_0$ to be the maximal time such
that the inequality
\beq \label{eq:enerNestdel}
\E_N(t)\leq 2 C_N\varepsilon (1+t)^\delta
\eq
holds for all $0\le t\le T$. Note that by the assumptions of the Theorem
$T>0$.
We will show
that if $\varepsilon>0$ is sufficiently small then this inequality
implies the same inequality with $2C_N$ replaced by $C_N$ for
$t\leq T$. Since the quantity is continuous this contradicts the
maximality of $T$ and it follows that the inequality holds for all
$T\le T_0$. Moreover, since the energy $\E_N(t)$ is now finite
at $t=T_0$ we can extend the solution beyond $T_0$ thus contradicting
maximality of $T_0$ and showing that $T_0=+\infty$.

The first step is to derive the preliminary decay estimates for $h^1$ under
the assumption \eqref{eq:enerNestdel}.
The estimate \eqref{eq:enerNestdel} can be explicitly written in the form
\begin{equation}\label{eq:apriorienergybound}
\sum_{Z\in {\cal Z}, |I|\leq N}\|w(q)^{1/2} \pa Z^I h^1(t,\cdot)\|_{L^2} \leq
C\varepsilon(1+t)^\delta ,\qquad 0<\delta< \gamma
\end{equation}
with some sufficiently large constant $C$.
The following result is a consequence of weighted global Sobolev inequalities proved in
Appendix C.
\begin{cor} \label{preest}
Let $h^1$ verify \eqref{eq:apriorienergybound} and $h^0$ be as in
\eqref{eq:h1-h0}.
 Then for $i=0,1$,
\begin{equation}\label{eq:weakdecayone}
|\pa Z^I h^i(t,x)|\leq \begin{cases}
C\varepsilon(1+t+|q|)^{-1+\delta} (1+|q|)^{-1-\delta^{\,\prime}},
\qquad q>0\\
C\varepsilon(1+t+|q|)^{-1+\delta}(1+|q|)^{-1/2},\qquad q<0
\end{cases},\qquad |I|\leq N-2.
\end{equation}
where $\delta^{\,\prime}=\delta$, if $i=0$ and
$\delta^\prime=\gamma>\delta$ if $i=1$. Furthermore,
\begin{equation}\label{eq:weakdecaytwo}
|Z^I h^i(t,x)|\leq \begin{cases}
C\varepsilon(1+t+|q|)^{-1+\delta}(1+|q|)^{-\delta^{\prime}},
\qquad q>0\\
C\varepsilon(1+t+|q|)^{-1+\delta}(1+|q|)^{1/2},\qquad q<0
\end{cases},\qquad |I|\leq N-2,
\end{equation}
and
\begin{equation}\label{eq:weakdecaythree}
|\overline{\pa} Z^I h^i(t,x)|\leq \begin{cases}
C\varepsilon(1+t+|q|)^{-2+\delta}(1+|q|)^{-\delta^{\prime}},
\qquad q>0\\
C\varepsilon(1+t+|q|)^{-2+\delta}(1+|q|)^{1/2},\qquad q<0
\end{cases},\qquad |I|\leq N-3
\end{equation}
\end{cor}
\begin{proof} We will only prove the estimates for $i=1$
since the estimates for $i=0$ follow by a direct calculation from the form of $h^0$.
Estimate \eqref{eq:weakdecayone} follows from the weighted Sobolev
inequality of Proposition \ref{prop:K-S}. We claim that
\beq\label{eq:initialdecayI}
\lim\limits_{|x|\to \infty}|Z^I h^1(0,x)|\to 0,\qquad |I|\le N-2
\eq
In fact, if $|I|=0$ this follows from \eqref{eq:initialdecay} and
\eqref{eq:weakdecayone}, and if $|I|\geq 1$ this follows from \eqref{eq:weakdecayone}, since $|Z\phi|\leq C(1+t+|x|)|\pa\phi|$. Estimate
\eqref{eq:weakdecaytwo} for $t=0$ follows by
integrating \eqref{eq:weakdecayone} from space-like infinity,  where
\eqref{eq:initialdecayI} hold.
\begin{remark}
The weighted Sobolev
inequality of Proposition \ref{prop:K-S} in fact implies the estimate
$$
|\pa Z^I h^1(t,x)|\leq \begin{cases}
C\varepsilon(1+t+|q|)^{-1} (1+|q|)^{-1-\gamma}(1+t)^\delta,
\qquad q>0\\
C\varepsilon(1+t+|q|)^{-1}(1+|q|)^{-1/2}(1+t)^\delta,\qquad q<0
\end{cases},\qquad |I|\leq N-2.
$$
In particular,
\begin{align}
&|\pa Z^I h^1(0,x)|\leq  C\varepsilon (1+|x|)^{-2-\gamma},\label{eq:Extra-1}\\
&|Z^I h^1(0,x)|\leq  C\varepsilon (1+|x|)^{-1-\gamma},\label{eq:Extra-2}
\end{align}
\end{remark}
The estimate  \eqref{eq:weakdecaytwo} for $r>t$ follows by integrating
\eqref{eq:weakdecayone} from the hyperplane $t=0$
along the lines with $t+r$ and
$\omega=x/|x|$ fixed:
\begin{multline}
|Z^I h^1(t,r\omega)|\leq  \int_r^{t+r}|\pa Z^I
h^1(t+r-\rho,\rho\omega)|\, d\rho+|Z^I
h^1\big(0,(t+r)\omega\big)|\\
\leq C\ve \int_r^\infty
\frac{d\rho}{(1+t+r)^{1-\delta}(1+|t-\rho|)^{1+\gamma}}+\frac{C\ve}{(1+t+r)^{1+\gamma-\delta}}
\leq\frac{C_\gamma\ve}{(1+t+r)^{1-\delta}(1+|t-r|)^\gamma}
\end{multline}
A  similar argument yields \eqref{eq:weakdecaytwo} for $r<t$.
 Inequality \eqref{eq:weakdecaythree}
follows from \eqref{eq:weakdecaytwo} using that $|\pab f|\les
\sum_{|I|=1}|Z^I f|/(1+t+|q|)$.
\end{proof}
The next subsection is devoted to the preparational estimates for
the inhomogeneous terms $F_{\mu\nu}$ and $F^0_{\mu\nu} =\Boxr_g
h^0_{\mu\nu}$ arising in the equation \eqref{eq:Red-E1} for the
tensor $h^1$.
\subsection{Estimates for the inhomogeneous terms}\label{section:decayinhom}
These estimates will play a key role in the derivation of the improved decay and energy estimates in the following two sections. As we have mentioned above the quadratic
terms in $F_{\mu\nu}$ in \eqref{eq:Red-E} do not satisfy the standard null condition. Nevertheless,
a special tensorial structure  and the wave coordinate condition will allow us to
obtain favorable estimates.

Recall that the inhomogeneous term $F_{\mu\nu}$ has the following structure:
\begin{align}
& F_{\mu\nu}(h)(\pa h,\pa h)=P(\pa_\mu h,\pa_\nu h)+Q_{\mu\nu}(\pa
h,\pa h)
+G_{\mu\nu}(h)(\pa h,\pa h)\label{eq:F},\\
& P(\pa_\mu \pi,\pa_\nu \theta):= \frac 12 \pa_{\mu} \pi^{\a\b}\pa_{\nu}
\theta_{\a\b}-\frac 14 \pa_{\mu} \tr \pi\, \pa_{\nu}\tr \theta,\label{eq:P}
\end{align}
The quadratic term $Q{\mu\nu}$ is a  linear combination of the null-forms and
$G_{\mu\nu}(h)(\pa h,\pa h)$ is a quadratic form in $\pa h$ with
the coefficients \-- a smooth function of $h$ vanishing at $h=0$.
The following result is an immediate consequence of \eqref{eq:P}, see also Lemma
\ref{tander}.
\begin{lemma}
\label{le:Pest} The quadratic form $P$ satisfies the following
pointwise estimate:
\begin{align}
&|P(\pa \pi,\pa \theta)|_{\cal TU}\les |\pab \pi\,|\, |\pa \theta| + |\pa \pi\,|\,
|\pab \theta|,
\label{eq:PTU}\\
&|P(\pa \pi,\pa \theta)|\les |\pa \pi\,|_{\cal TU} |\pa \theta|_{\cal TU} + |\pa
\pi\,|_{\cal LL} |\pa \theta| + |\pa \pi\,| \, |\pa \theta|_{\cal LL}\label{eq:Pall}
\end{align}
\end{lemma}
Using the additional estimates on the $h_{LL}$ component, derived
in Proposition \ref{decaywavecZ} under the assumption that the wave
coordinate condition holds, we obtain the following:
\begin{cor}
Let metric $g$ satisfy the wave coordinate condition \eqref{eq:wavecdef} relative to coordinates $\{x^\mu\}$, the quadratic form $P$ obeys the
following estimate on a 2-tensor $h_{\mu\nu}=g_{\mu\nu}-m_{\mu\nu}$:
\begin{align}
&|P(\pa h,\pa h)|_{\cal TU}\les |\pab h|\, |\pa h|,\label{eq:Ph}\\
&|P(\pa h,\pa h)|\les |\pa h|^2_{\cal TU} + |\pab h|\, |\pa h| +
|h|\,|\pa h|^2\label{eq:Pallcoord}
\end{align}
In addition, assuming that $|Z^J h|\le C$ for all multi-indices
$|J|\le |I|$ and vector fields $Z\in {\cal Z}$,
\begin{align*}
|Z^I P(\pa h, \pa h)|&\les \sum_{|J|+|K|\le |I|} \big (|\pa Z^J
h|_{\cal TU}\, |\pa Z^K h |_{\cal TU} + |\pab Z^J h|\, |\pa Z^K h|
\big )+ \sum_{|J|+|K|\le |I|-1} |\pa Z^J h|_{L \cal T} |\pa Z^K
h|\\ & + \sum_{|J|+|K|\le |I|-2} |\pa Z^J h|\, |\pa Z^K h| +
\sum_{|J_1|+|J_2|+|J_3| \leq |I|} |Z^{J_3} h|\,|\pa Z^{J_{2}} h| |\pa Z^{J_1} h|
\label{eq:PZI}
\end{align*}
\end{cor}
\begin{proof} The estimate \eqref{eq:Ph} follows directly
from \eqref{eq:PTU}. To prove \eqref{eq:Pallcoord} we use
\eqref{eq:Pall} and that by the wave coordinate condition $|\pa
h|_{LL}\les |\bar \pa h|+|h|\, |\pa h|$.

We now note that $Z^I P(\pa_\mu h,\pa_\nu h)$ is a sum of terms of
the form $P(\pa_\alpha Z^J h,\pa_\beta Z^K h)$ for some $\alpha$,
$\beta$ and $|J|+|K|\leq I$:
$$
|Z^I P(\pa h,\pa h)|\leq C\sum_{|J|+|K|\leq |I|} |P(\pa  Z^J h,\pa
Z^K h)|
$$
It follows from \eqref{eq:Pall} and Proposition \ref{decaywavecZ}
that
\begin{multline}
\!\!\!\!\!\sum_{|J|+|K|\leq |I|}\!\!\!\!\! |P(\pa Z^J h,\pa Z^K
h)| \les \!\!\!\!\!\sum_{|J|+|K|\leq |I|}\!\!\!\!\! |\pa Z^J
h|_{\cal TU}|\pa Z^K h|_{\cal TU}
+|\pa Z^J h|_{LL} \,|\pa Z^K h |\\
\les \sum_{|J|+|K|\leq |I|} |\bar\pa Z^J h|\, |Z^K h|
+|\pa Z^J h|_{\cal TU}\, |\pa Z^K h|_{\cal TU} \\
+\!\!\!\!\! \sum_{|J|+|K|\leq |I|}\!\!\!\!\Big(
\sum_{|J^\prime|\leq |J|-1} \!\!\!\!\!\! |\pa Z^{J^\prime}
h|_{\cal LT}\,\, +\!\!\!\!\!\!\sum_{|J^{\prime\prime}|\leq |J|-2}
\!\!\!\!\!|\pa Z^{J^{\prime\prime}} h|\,\,
+\!\!\!\!\!\sum_{|J_1|+|J_2|\leq |J|}
|Z^{J_2} h|\, |\pa Z^{J_1} h|\Big) |\pa Z^K h|
\end{multline}
which proves the result.
\end{proof}
We now state the complete estimates for the inhomogeneous term $F_{\mu\nu}$.
\begin{prop} \label{decayinhom}
Let $F_{\mu\nu}=F_{\mu\nu}(h)(\pa h,\pa h)$ be  as in \eqref{eq:F}
 and assume that the wave coordinate condition \eqref{eq:wavecdef}
holds for the metric $g_{\mu\nu}=m_{\mu\nu} + h_{\mu\nu}$ relative to
coordinates $\{x^\mu\}$. Then
\begin{align}
&| F|_{\cal T U}\les
|\overline{\pa} h|\, |\pa h|+|h|\,|\pa h|^2 ,\label{eq:decayinhom1} \\
&|F|\les |\pa h|^2_{\cal TU} + |\pab h|\,
|\pa h| + |h|\,|\pa h|^2 \label{eq:decayinhom2}
\end{align}
In addition, assuming that $|Z^J h|\le C$ for all multi-indices
$|J|\le |I|$ and vector fields $Z\in {\cal Z}$,
\begin{multline} \label{eq:decayinhom4}
|Z^I F| \les
 \sum_{|J|+|K|\le |I|}
\big (|\pa Z^J h|_{\cal TU}\, |\pa Z^K h |_{\cal TU} + |\pab Z^J
h|\, |\pa Z^K h| \big )+ \sum_{|J|+|K|\le |I|-1} |\pa Z^J h|_{L
\cal T} |\pa Z^K h|\\  + \sum_{|J|+|K|\le |I|-2} |\pa Z^J h|\,
|\pa Z^K h| + \sum_{|J_1|+|J_2|+|J_3| \leq |I|}
|Z^{J_3} h|\,|\pa Z^{J_{2}} h| |\pa Z^{J_1} h|
\end{multline}
\end{prop}
\begin{proof} First
$$
|Z^I G_{\mu\nu}(h)(\pa h,\pa h)| \leq C\sum_{ |I_1|+...+|I_k|\leq
|I|,\, k\geq 3}
 |Z^{I_k} h|\cdot\cdot\cdot |Z^{I_3} h|\,|\pa Z^{I_2} h|
\,|\pa Z^{I_1} h|.
$$
Since $Z Q(\pa u,\pa v)=Q(\pa u,\pa Zv)+Q(\pa Zu,\pa v)
+a^{ij}Q_{ij}(\pa u,\pa v)$, and $|Q_{\mu\nu}( \pa h,\pa k)|\leq
|\pa h|\, |\overline{\pa} k|+
 |\pa k|\, |\overline{\pa} h|$ it follows that
$$
|Z^I Q_{\mu\nu}(\pa h,\pa h)|\leq C\sum_{|J|+|k|\leq |I|}
|Q_{\mu\nu}(\pa  Z^J h,\pa Z^K h)|\leq  C\sum_{|J|+|k|\leq |I|}
|\pa Z^J h|\, |\overline{\pa} Z^K h|
$$
and the proposition follows.
\end{proof}
The next step is to estimate the extra inhomogeneous term $\Boxr_g h^0$
appearing in the wave equation $\Boxr_g h^1_{\mu\nu} = F_{\mu\nu}-
\Boxr_g h^0_{\mu\nu}$ via the decomposition $h=h^1+h^0$.
\begin{lemma} \label{subtractoffschwarz} Let
$$
F^0_{\mu\nu}=\Boxr_g h^0_{\mu\nu},\qquad h^0_{\mu\nu}=\chi(r)\chi(t/r)\frac{M}{r}\, \de_{\mu\nu}
$$
Then, if the weak decay estimates in Corollary \ref{preest} hold,
we have
$$
|Z^I F^{0}|\leq \begin{cases}
C\varepsilon^2(t+|q|+1)^{-4+\delta}(1+|q|)^{-\delta},\quad q>0,\\
C\varepsilon(t+|q|+1)^{-3} ,\quad q<0,
\end{cases},\qquad |I|\leq N-2
$$
More generally,
$$
|Z^I F^{0}|\leq \begin{cases} C_N\varepsilon^2(t+|q|+1)^{-4}
,\quad q>0,\\
C_N\varepsilon(t+|q|+1)^{-3} ,\quad q<0,
\end{cases}+\frac{C_N\varepsilon }{(t+|q|+1)^3}
\sum_{|J|\leq |I|}
  |Z^{J}h^1|,\qquad |I|\leq N.
$$
\end{lemma}
\begin{proof}
Set
$$
F^{0}=\Boxr_g h^0=F^{\,00}+F^{\,01},\qquad F^{\,00}=\Box \,
h^0,\quad F^{\,01}= H^{\,\alpha\beta}\pa_\alpha\pa_\beta h^0
$$
One can easily check that
\begin{align*}
&|Z^I F^{\,00}|\leq  \frac{C_N\varepsilon}{(t+|q|+1)^3},
\quad\text{and}\quad  F^{\,00}\equiv 0,\quad \text{for}\,\, r<t/2
\,\,{\text or}\,\,
r>3t/4\\
&|Z^I F^{\,01}|\les \frac{C_N\varepsilon
}{(t+|q|+1)^3}\sum_{|J|\leq |I|} |Z^J H|
\end{align*}
On the other hand, with the help of  \eqref{eq:weakdecaytwo} for $h^0$,
 $$
\sum_{|J|\leq k}|Z^J H|\les \frac{C_N\varepsilon }{t+|q|+1}
+\sum_{|J_1|+...+|J_n|\leq k} |Z^{J_1}h^1|\cdot\cdot\cdot|Z^{J_n}
h^1|,
 $$
 since $H=-h+O(h^2)$, $h=h^0+h^1$. Using \eqref{eq:weakdecaytwo}, this time
 for $h^1$, we obtain
  $$
|Z^I F^{\,01}|\leq \begin{cases}
C_N\varepsilon^2(t+|q|+1)^{-4+\delta} (1+|q|)^{-\delta},\quad q>0,\\
C_N\varepsilon^2(t+|q|+1)^{-4+\delta} (1+|q|)^{1/2} ,\quad q<0,
\end{cases},\qquad |I|\leq N-2
$$
\end{proof}

\section
{Decay estimates for the Einstein equations}\label{section:decayeinst}
In this section we put to use the wave coordinate condition
\eqref{eq:Wave-Let} for the metric $g(t)$ in coordinates $\{x^\mu\}_{\mu=0,...,3}$,
which, as explained in the Remark \ref{rem:Remark}, is satisfied  on the maximum
interval of existence $[0,T_0]$, and the decay estimates for solutions of the wave
equation $\Boxr_g\phi_{\mu\nu}=W_{\mu\nu}$, derived in Corollary \ref{decaywaveeq3}, to upgrade pointwise estimates of Corollary  \ref{preest}.
\begin{prop}[Estimates for $h$]\label{decayeinst}
Let $h=h^1+h^0$ be a solution of the reduced Einstein equations \eqref{eq:Red-E}.
Assume that $h^1$ verifies the energy estimate \eqref{eq:enerNestdel} on the time
interval $[0,T]$. Then for any $t\in [0,T]$,
 \begin{align}
|\pa h|_{L \cal T} +|\pa Z h|_{LL}\leq
&\begin{cases} C\varepsilon
(1+t+|q|)^{-2+\delta}(1+|q|)^{-\delta},
\quad q>0\\
C\varepsilon (1+t+|q|)^{-2+\delta}(1+|q|)^{1/2}, \quad q<0
\end{cases},\label{eq:sharpdecayone}\\
|h|_{L \cal T} +|Z h|_{LL}\leq
&\begin{cases}
C\varepsilon (1+t+|q|)^{-1},\qquad q>0.\\
C\varepsilon (1+t+|q|)^{-1}(1+|q|)^{1/2+\delta} ,\qquad q<0.
\end{cases} \label{eq:sharpdecaytwo}
\end{align}
Furthermore,
\begin{align}
&|\pa h|_{\cal T U}\leq  C\varepsilon (1+t+|q|)^{-1},\label{eq:sharpdecaythree}\\
&|\pa h|\leq  C\varepsilon t^{-1} \ln t .\label{eq:sharpdecayfour}
\end{align}
\end{prop}
\begin{prop}[Estimates for $h^1$]\label{prop:decayeinst}
Under the same assumptions
let $\gamma^\prime<\gamma-\delta$ and $\mu^\prime>\delta>0$ be
fixed. Then there exist constants $M_k$, $C_k$ and
$\varepsilon_k$, depending on $(\gamma^\prime,\mu^\prime,\delta)$,
such that
\begin{align}
| \pa Z^I h^1|\leq &\begin{cases} C_k\varepsilon
(1+t+|q|)^{-1+M_k\varepsilon}(1+|q|)^{-1-\gamma^\prime},\quad q>0\\
C_k \varepsilon
(1+t+|q|)^{-1+M_k\varepsilon}(1+|q|)^{-1/2+\mu^\prime},\quad q<0
\end{cases}
  ,\qquad  |I|=k\leq N/2+2 \label{eq:sharpdecayfive}\\
| Z^I h^1|\leq &\begin{cases} C_k\varepsilon
(1+t+|q|)^{-1+M_k\varepsilon}(1+|q|)^{-\gamma^\prime},\quad q>0\\
C_k\varepsilon
(1+t+|q|)^{-1+M_k\varepsilon}(1+|q|)^{1/2+\mu^\prime},\quad q<0
\end{cases}
  ,\qquad  |I|=k\leq N/2+2 \label{eq:sharpdecaysix}
\end{align}
 The same estimates hold for $h^0$ if we replace
$\gamma^\prime$ by $M_k\varepsilon$.
\end{prop}
\begin{remark}
Note the difference between the estimates \eqref{eq:sharpdecaythree},
\eqref{eq:sharpdecayfour} for $h$ and \eqref{eq:sharpdecayfive} with
$|I|=0$ for  $h^1$.
\end{remark}
The estimates \eqref{eq:sharpdecayone}-\eqref{eq:sharpdecaytwo} follow
from the wave coordinate condition combined with the weak decay estimates
of Corollary \ref{preest}. In the derivation of the sharp decay
estimates \eqref{eq:sharpdecaythree}-\eqref{eq:sharpdecaysix} we will use
decay estimates of  Corollary \ref{decaywaveeq3} with various weights.

\subsection{Proof of \eqref{eq:sharpdecayone}-\eqref{eq:sharpdecaytwo}.}
We rely on Proposition \ref{decaywavecZ} to establish an even more
general version of the desired estimates.
We remind the reader that under the assumptions of Proposition \ref{decayeinst}
both the wave coordinate condition \eqref{eq:Wave-Let} for tensor $g_{\mu\nu}$
and the weak decay estimates of Corollary \ref{preest} hold true.
\begin{lemma} \label{wavecdecay4} Under the assumptions of Proposition \ref{decayeinst}
\begin{align}
&\sum_{|I|\leq \, k}|\pa Z^I h|_{LL} + \sum_{|J|\leq \, k-1} |\pa
Z^{J} h|_{\cal LT}\les \sum_{|K|\leq \, k-2} |\pa Z^K h|
+\begin{cases} \varepsilon
(1+t+|q|)^{-2+2\delta}(1+|q|)^{-2\delta},\qquad q>0\\
\varepsilon (1+t+|q|)^{-2+2\delta}(1+|q|)^{1/2-\delta},\qquad
q<0\end{cases}  \label{eq:absentone}\\
& \sum_{|I|\leq \,
k}|Z^I h|_{LL} + \sum_{|J|\leq \, k-1} |Z^{J} h|_{\cal LT}\les
\sum_{|K|\leq \, k-2} \,\,\,\int\limits_{s,\omega=const}
|\pa Z^K h| +\begin{cases} \varepsilon
(1+t+|q|)^{-1},\qquad q>0\\
\varepsilon (1+t+|q|)^{-1}(1+|q|)^{1/2+\delta},\qquad
q<0\end{cases} \label{eq:absenttwo}
\end{align}
 Here the sums over $k-2$ are absent if $k\leq 1$,
 the sums over $k-1$ are absent if $k=0$ and $\int\limits_{s,\omega=const}$
 stands for an integral along the segment $\tau+|y|, y/|y|=const$ connecting
 a given point $(t,x)$ with a hyperplane $t=0$.
\end{lemma}
\begin{proof}
The proof follows immediately from Proposition  \ref{decaywavecZ}
and the estimates of Corollary \ref{preest} since
\beq
 \sum_{|I|\leq \,k}|\pa Z^I H|_{LL}
 +\sum_{|J|\leq \,k-1}|\pa Z^J H|_{L\cal T} \les
\sum_{|J|\leq |I|} |\pab Z^J H|+ \sum_{|K|\leq \,k-2}|\pa Z^K H|
+\sum_{I_{1}+I_2=I} |Z^{I_2}  H| \, |\pa Z^{I_1} H| .
\end{equation}
 Integrating \eqref{eq:absentone} along the lines on which the angle $\omega=y/|y|$
 and the null coordinate $s=\tau+|y|$ are constant, as in the
 proof of Corollary \ref{preest}, and using  \eqref{eq:weakdecaytwo}
 at $t=0$, yields \eqref{eq:absenttwo}.
\end{proof}

\subsection{Proof of \eqref{eq:sharpdecaythree}-\eqref{eq:sharpdecayfour}.}
We apply  the $L^\infty$
estimates, derived in Corollary \ref{decaywaveeq3}, for
the reduced wave equation
$$
{\Boxr}_g h_{\mu\nu}=F_{\mu\nu}
$$
with $F_{\mu\nu}$ given in \eqref{eq:F-prime} or
\eqref{eq:F}-\eqref{eq:P}.
Observe that the weak decay estimates of Corollary \ref{preest}
guarantee that the tensor $H^{\mu\nu}=g^{\mu\nu} - m^{\mu\nu}$
verifies the required assumptions of Corollary \ref{decaywaveeq3}.

First we derive the $L^\infty$ estimates for $F_{\mu\nu}$.
\begin{lemma} \label{inhomdecay}
Suppose that the assumptions of Proposition \ref{decayeinst} hold
and let $F_{\mu\nu}=F_{\mu\nu}(h)(\pa h,\pa h)$ be  as in
 \eqref{eq:F-prime}. Then
 \begin{align}
 &| F|_{\cal
T U}\leq
 C\varepsilon t^{-3/2+\delta}|\pa h|,\label{eq:decayinhom1} \\
&|F|\leq C\varepsilon
t^{-3/2+\delta}|\pa h| +C|\pa h|_{\cal TU}^2
\label{eq:decayinh}
\end{align}
\end{lemma}
\begin{proof} This follows from Proposition \ref{decayinhom} and
the weak decay estimates of Corollary \ref{preest}.
\end{proof}
The 2-tensor $h$ is a solution of the wave equation
$\Boxr_g h_{\mu\nu} = F_{\mu\nu}$.
We now recall that that according to Corollary \ref{decaywaveeq3}
with a trivial weight $\varpi(q)=1$, i.e., $\alpha=0$,
we have that for any $U,V\in\{L,\Lb,A,B\}$,
\begin{multline}
(1+t+|x|)| \,|\pa h(t,x)|_{UV} \les\!\sup_{0\leq
\tau\leq t}
\sum_{|I|\leq 1}\|\,Z^I\! h(\tau,\cdot)\|_{L^\infty}\\
+ \int_0^t\Big( (1+\tau)\|
\,|F(\tau,\cdot)|_{UV}\|_{L^\infty(D_\tau)} +\sum_{|I|\leq 2} (1+\tau)^{-1}
\| Z^I h(\tau,\cdot)\|_{L^\infty(D_\tau)}\Big)\, d\tau.
\end{multline}
Then with the help of Lemma \ref{inhomdecay} and the weak decay
estimates of Corollary \ref{preest}
we obtain
\begin{lemma} With a constant depending on $\gamma>0$ we have
\begin{align*}
&(1+t)\|\, |\pa h|_{\cal T U} (t,\cdot)\|_{L^\infty} \leq
C\varepsilon + C\varepsilon \int_0^t(1+\tau)^{\delta-1/2}\|\pa
h(\tau,\cdot)\|_{L^\infty}\, d\tau,\\
&(1+t)\|\pa h(t,\cdot)\|_{L^\infty} \leq  C\varepsilon + C\int_0^t
\Big(\varepsilon (1+\tau)^{\delta-1/2}\|\pa
h(\tau,\cdot)\|_{L^\infty}+ (1+\tau)\||\pa h|_{\cal T
U}(\tau,\cdot)\|_{L^\infty}^2\Big)\, d\tau.
\end{align*}
\end{lemma}
Here, by assumption $\delta<1/4$ so $\delta-1/2<-1/4$.
The estimates  \eqref{eq:sharpdecaythree} and \eqref{eq:sharpdecayfour}
 are now consequences of
the above lemma and the following technical result applied to
the functions $b(t):=(1+t)\|\, |\pa h|_{\cal T U} (t,\cdot)\|_{L^\infty}$
and $c(t):=(1+t)\|\pa h(t,\cdot)\|_{L^\infty}$:
\begin{lemma} \label{ode1}
Assume that the functions $b(t)\geq 0$ and $c(t)\geq 0$ satisfy
\begin{align}
b(t)&\leq C\varepsilon \Big(\int_0^t (1+s)^{-1-a} c(s)\, ds+1\Big) \label{eq:ode1}\\
c(t)&\leq  C\varepsilon \Big(\int_0^t (1+s)^{-1-a}
c(s)\, ds+1\Big) +C\int_0^t  (1+s)^{-1} b^2(s)\, ds\label{eq:ode2}
\end{align}
for some positive constants such that
$a\ge C^2\varepsilon$ and $a\ge 4C\varepsilon/(1-2C\varepsilon)$.
Then \beq\label{eq:odebound}
b(t)\leq 2C\varepsilon, \qquad\text{and}\qquad c(t)\leq
2C\varepsilon \big(1+a \ln{(1+t)}\big) \eq
\end{lemma}
\begin{proof} Let $\tau_0$ be the largest time such that
\label{eq:odebound} holds.
Substituting these bounds into \eqref{eq:ode1}-\eqref{eq:ode2}
and taking into account that
$$
\int_0^\infty (1+s)^{-1-a}   \big(1+a\ln{(1+s)}\big)\, ds\le
a^{-1}\int_0^\infty (1+\tau) \, e^{-\tau}\, d\tau =2 a^{-1}+1
$$
we obtain that for any $0\le t\le \tau_0$
\begin{align*}
&b(t)\le C\varepsilon (2C\varepsilon(1+ 2a^{-1}) +1)<2C\varepsilon,\\
&c(t)\le C\varepsilon \big (2C\varepsilon(1+ 2a^{-1})+ C^2\varepsilon\ln(1+t)\big )
<2C\varepsilon \big(1+a \ln{(1+t)}\big),
\end{align*}
which implies that $\tau_0=\infty$, as desired.
\end{proof}

\subsection{Proof of \eqref{eq:sharpdecayfive}-\eqref{eq:sharpdecaysix}}
The proof of the estimates for the tensor $h^1$ proceeds by induction.
We assume that \eqref{eq:sharpdecayfive}-\eqref{eq:sharpdecaysix} hold
for all values of multi-index $|I|\le k$ and prove the estimate for $|I|=k+1$.
The argument below will also apply unconditionally to the base of the
induction $k=0$.

Once again the first step is to establish $L^\infty$ estimates for
$Z^I F$.
\begin{lemma} \label{inhomdecay}
Suppose that the assumptions of Proposition \ref{decayeinst} hold
and let $F_{\mu\nu}=F_{\mu\nu}(h)(\pa h,\pa h)$ be  as in \eqref{eq:F-prime}.
Then \beq
 |Z^I F|\leq C\varepsilon  \sum_{|K|\leq |I|}\frac{|\pa Z^K h|}{1+t}
+\,\, C\sum_{|J|+|K|\leq |I|,\, |J|\leq |K|<|I|}|\pa Z^{J} h||\pa
Z^{K} h| \eq
\end{lemma}
\begin{proof} The result  follows from Proposition \ref{decayinhom} with the help of
\eqref{eq:weakdecayone}-\eqref{eq:weakdecaythree} and
\eqref{eq:sharpdecaythree}.
\end{proof}
Recall now that the 2-tensor $h^1=g-m-h^0$ is a solution of the reduced wave equation
\beq\label{eq:no-h}
\Boxr_g h^1_{\mu\nu}=F^1_{\mu\nu}=
F_{\mu\nu}-F^0_{\mu\nu},\qquad \text{where}\qquad F^0_{\mu\nu}=
\Boxr_g h^0_{\mu\nu}
\eq
and that the terms $F_{\mu\nu}$ and $F_{\mu\nu}^0$ have been treated
in Proposition \ref{decayinhom} and Lemma \ref{subtractoffschwarz} respectively.

To prove the estimates for $Z^I h^1$ with vector fields $Z\in {\cal
Z}$ we commute the equation \eqref{eq:no-h} with $Z^I$. By
Proposition \ref{prop:commut3}
\begin{align}
&|{\Boxr}_g Z^I h^1|\les |\hat{Z}^I F^1| +(1+t)^{-1}\sum_{|K|\leq
|I|,}\,\,\sum_{|J|+(|K|-1)_+\leq |I|} \!\!\!\!\!\!\!\!\!\!
|Z^{J} H| |\pa Z^{K} h^1|\label{eq:nas-T1}\\
+ \frac C{1+|q|}
 &\sum_{|K|\leq |I|}\Big(\sum_{|J|+(|K|-1)_+\leq |I|} \!\!\!\!\!|Z^{J} H|_{LL}
+\!\!\!\!\!\sum_{|J^{\prime}|+(|K|-1)_+\leq
|I|-1}\!\!\!\!\!|Z^{J^{\prime}} H|_{L\cal T}
+\!\!\!\!\!\sum_{|J^{\prime\prime}|+(|K|-1)_+\leq |I|-2}\!\!\!\!\!
|Z^{J^{\prime\prime}} H|\Big)|\pa Z^{K} h^1|\label{eq:nas-T2}
\end{align}
where $(|K|-1)_+=|K|-1$, if $|K|\geq 1$, and $0$, if $|K|=0$.
Using Lemma \ref{wavecdecay4} and \eqref{eq:sharpdecayfive}, which is
inductively assumed to be true for $|I|\leq k$, we get
\begin{multline} (1+|q|)^{-1}\!\!\!\!\sum_{|J|\leq
k,\,|J^\prime|\leq k-1,\,|J^{\prime\prime}|\leq k-2}
\Big ( |Z^J H|_{LL}+|Z^{J^\prime} H|_{L\cal T}
+|Z^{J^{\prime\prime}} H|\Big )\\
\les
\begin{cases} C_k\varepsilon
(1+t+|q|)^{-1+M_k\varepsilon}(1+|q|)^{-1-M_k
\varepsilon}\quad+\quad
\varepsilon (1+t+|q|)^{-1}(1+|q|)^{-1},\qquad q>0\\
C_k \varepsilon
(1+t+|q|)^{-1+M_k\varepsilon}(1+|q|)^{-1/2+\mu^\prime}\quad+\quad
\varepsilon(1+t+|q|)^{-1}(1+|q|)^{-1/2+\delta},\qquad
q<0\end{cases}
\end{multline}
while contribution of the terms in \eqref{eq:nas-T2} coupled to the highest
order term \-- $\pa Z^K h^1,\, |K|=|I|$, according to \eqref{eq:sharpdecaytwo}
amounts to
\begin{equation} (1+|q|)^{-1} \sum_{|J|\leq
1}
\Big ( |Z^J H|_{LL}+|H|_{L\cal T}\Big )
 \les \begin{cases}
\varepsilon (1+t+|q|)^{-1}(1+|q|)^{-1},\qquad q>0\\
\varepsilon(1+t+|q|)^{-1}(1+|q|)^{-1/2+\delta},\qquad
q<0\end{cases}
\end{equation}
Therefore splitting the sum in \eqref{eq:nas-T2} according to whether
$|K|<|I|$ or $|K|=|I|$ and using the inductive assumption
 \eqref{eq:sharpdecayfive} for $|K|<|I|=k+1$, we obtain
\beq |{\Boxr}_g Z^I h^1|\leq C\varepsilon \sum_{|K|\leq
|I|}\frac{|\pa Z^K h^1|}{1+t} +|\hat{Z}^I F^0|+
\begin{cases} C_k^2\varepsilon^2
(1+t+|q|)^{-2+2M_k\varepsilon}(1+|q|)^{-2-2M_k
\varepsilon},\qquad q>0\\
C_k^2 \varepsilon^2
(1+t+|q|)^{-2+2M_k\varepsilon}(1+|q|)^{-1+2\mu^\prime},\qquad
q<0\end{cases} \eq
Now by  Lemma \ref{subtractoffschwarz}
 $$
|Z^I F^{0}|\leq \begin{cases}
C\varepsilon^2(t+|q|+1)^{-4+\delta}(1+|q|)^{-\delta}
,\quad q>0,\\
C\varepsilon(t+|q|+1)^{-3} ,\quad q<0,
\end{cases},\qquad |I|\leq N-2
$$
Set
\beq
n_{k+1}(t)=(1+t)\sum_{|I|\leq k+1}\| \varpi(q)\pa Z^I
h^1(t,\cdot)\|_{L^\infty},\qquad
\varpi(q)=\begin{cases}
(1+|q|)^{1+\gamma^\prime},\quad q>0\\
      (1+|q|)^{1/2-\mu^\prime}\quad
      q<0\end{cases},
 \eq
 where $\mu^\prime>\delta$ and $\gamma^\prime<\gamma-\delta$.
 Then we have established that for
$|I|=k+1$:
\beq
\varpi(q)|\Boxr_g Z^I h^1|\les (1+t)^{-2}\big(
\varepsilon n_{k+1}(t)+ \varepsilon^2 C_k^2(1+t)^{2M_k\varepsilon}
+C\varepsilon (1+t)^{-1/2-\mu^\prime}\big)
\eq
The weak decay estimates of Corollary \ref{preest} imply that
 $$
 \varpi(q)|Z^I h^1(t,x)|\les \begin{cases}
\varepsilon(1+t+|q|)^{-1+\delta}(1+|q|)^{-\gamma}
(1+|q|)^{1+\gamma^\prime}
,\quad q>0\\
\varepsilon(1+t+|q|)^{-1+\delta}(1+|q|)^{1/2}
(1+|q|)^{1/2-\mu^\prime},\quad q<0
\end{cases}\les\,\, \varepsilon (1+t)^{-a}, \quad a>0
 $$
 provided that
 $a=\min{(\mu^\prime-\delta,\gamma-\delta-\gamma^\prime)}>0$.
The decay estimate proved in
Corollary \ref{decaywaveeq3} therefore shows that for some constant $C=C(k)$
we have  the inequality
\beq n_{k+1}(t)\leq C\varepsilon+ C\int_0^t (1+\tau)^{-1} \big(
\varepsilon n_{k+1}(\tau)+  \varepsilon^2
(1+\tau)^{C\varepsilon}+C\varepsilon (1+t)^{-1/2-\mu^\prime}\big)
\, d\tau \eq
The bound $n_{k+1}(t)\!\leq  2 C\varepsilon (1+t)^{2C\varepsilon}$
follows from the Gronwall inequality.
This proves \eqref{eq:sharpdecayfive}. The estimate \eqref{eq:sharpdecaysix}
then follows by integrating \eqref{eq:sharpdecayfive} along the line
$\omega=y/|y|, \tau+|y|=const$ from the hyperplane $t=0$ using
\eqref{eq:Extra-2}.

\section {Energy estimates for Einstein's equations}\label{section:energyeisnt}
Recall the definition of the weighted energy
\begin{align}
\E_{N}(t)=\sup_{0\le\tau\le t}\sum_{|I|\le N}\int_{\Si_{t}} |\pa
Z^{I} h^1|^{2}w(q),
\label{eq:energyN}
\end{align}
where \beq
w=\begin{cases} 1+(1+|q|)^{1+2\gamma},\quad q>0\\
         1+(1+|q|)^{-2\mu},\quad
         q<0\end{cases},\qquad
         w^\prime=\begin{cases} (1+2\gamma)(1+|q|)^{2\gamma},\quad q>0\\
         2\mu(1+|q|)^{1-2\mu},\quad q<0\end{cases}
         \eq
Recall also our decomposition
$$
g_{\mu\nu}(t)=m_{\mu\nu} + h_{\mu\nu}(t)=m_{\mu\nu} +
h^1_{\mu\nu}(t)+ h^0_{\mu\nu}(t),\qquad h^0_{\mu\nu}(t)=\chi(\frac
rt)\chi(r) \frac Mr \delta_{\mu\nu}
$$
of a local in time smooth solution $g_{\mu\nu}(t)$ of the reduced Einstein equations,
and the definition of the tensor $H^{\mu\nu}=g^{\mu\nu}-m^{\mu\nu}$.

 In this section we prove the following result.
\begin{theorem}\label{energyest}
Let $g_{\mu\nu}(t)=h_{\mu\nu}(t)+m_{\mu\nu}$
be a local in time solution of the reduced Einstein equations
\eqref{eq:Red-E}) satisfying the wave coordinate condition \eqref{eq:Wave-Let}
 on the interval  $[0,T)$. Suppose also that for some $0<\mu^\prime<1/2$ and
 $0<\gamma <1/2$ we have the following estimates for $0\leq t\leq T$,
all multi-indices
 $|I|\le N/2+2$ and the collections ${\cal T}=\{L,S_1,S_2\},\,\,
 {\cal U}=\{L,\Lb,S_1,S_2\}$:
\begin{align}
&|\pa H|_{\cal TU} + (1+|q|)^{-1} | H|_{\cal TL}+ (1+|q|)^{-1} |Z
H|_{\cal LL}\le C\varepsilon (1+t)^{-1},\label{eq:specialcomp}\\
&| \pa Z^I h|
+ \frac{|Z^{I} h|}{1+|q|}+\frac{1+t+|q|}{1+|q|} |\pab Z^I h|\les
\begin{cases} \label{eq:generalcomp}
C\varepsilon
(1+t+|q|)^{-1+C\varepsilon}(1+|q|)^{-1-C\varepsilon},\quad q>0\\
C\varepsilon
(1+t+|q|)^{-1+C\varepsilon}(1+|q|)^{-1/2+\mu^\prime},\quad q<0
\end{cases}\\
&E_N(0)+M^2\leq \ve^2.\label{eq:initialenergy}
 \end{align}
Then there is a positive constant $c$ independent of $T$ such
that if $\ve \leq c^{-2}$ we have the energy estimate
\beq\label{eq:energyestk} \E_N(t)\leq C_N\, \varepsilon^{2}
(1+t)^{c\varepsilon},
\end{equation}
for $0\leq t\leq T$. Here $C_N$ is a
constant that depends only on $N$.
\end{theorem}
Assuming the conclusions of Theorem \ref{energyest} for a moment
we finish the proof of the Main Theorem \ref{exist}.
\subsection{End of the proof of Theorem \ref{exist}}
Recall that $T$ was defined as the maximal time with the property
that the bound
$$
\E_N(t)\le 2C_N\ve (1+t)^\delta
$$
holds for all $0\le t\le T$. Assuming the energy bound above
we have establsihed in Propositions \ref{decayeinst}-\ref{prop:decayeinst}
the decay estimates for the
tensors $h=h^0+h^1$ and $h^1$ respectively. Direct check shows that
the estimates of Propositions \ref{decayeinst}-\ref{prop:decayeinst}
imply the assumptions \eqref{eq:specialcomp}-\eqref{eq:generalcomp}.
The conclusion of Theorem \ref{energyest} states that the energy
$$
\E_N(t)\le C_N \ve^2 (1+t)^{c\ve},\qquad \forall 0\le t\le T
$$
Thus choosing a sufficiently small $\ve>0$ we can show that
$\E_N(t)\le C_N\ve (1+t)^\delta$ thus contracting the maximality
of $T$ and consequently proving that $g_{\mu\nu}$ is a global solution.
It therefore remains to prove  Theorem \ref{energyest}.
\subsection{Proof of Theorem \ref{energyest}}
\begin{proof}
Recall that the components of the tensor $h_{\mu\nu} = g_{\mu\nu}
- m_{\mu\nu}$ satisfy the wave equations:
\begin{align}
&\Boxr_g h_{\mu\nu}=
F_{\mu\nu},\label{eq:equat},\\
&F_{\mu\nu}= P(\pa_{\mu} h,\pa_{\nu} h) + Q_{\mu\nu}(\pa h,\pa
h)+G_{\mu\nu}(h)(\pa h,\pa h),\nn\\
&P(\pa_{\mu}h,\pa_{\nu} h) =
\frac{1}{2}  m^{\alpha\alpha^\prime}m^{\beta\beta^\prime} \pa_\mu
h_{\alpha\beta}\, \pa_\nu h_{\alpha^\prime\beta^\prime}
-\frac{1}{4}  m^{\alpha\alpha^\prime}\pa_\mu
h_{\alpha\alpha^\prime} \,  m^{\beta\beta^\prime}\pa_\nu
h_{\beta\beta^\prime}. \label{eq:defiP}
\end{align}
Our goal is to compute the energy norms of $Z^I h^1$, where $Z\in {\cal Z}$
and $h^1$ is a solution of the problem
\begin{equation}\label{eq:h1waveeq}
\Boxr_{g} \, h^1_{\mu\nu}= {F}^1_{\mu\nu},\qquad\text{where}\qquad
F^1=F-F^0,\qquad F^0=\Boxr_g h^0
\end{equation}
Commuting with the vector fields $\hat Z^I,\, \, \hat Z=Z+c_Z$ we get
 \begin{equation}\label{eq:h1Iwaveeq}
\Boxr_{g} \, Z^I h^1_{\mu\nu}= F_{\mu\nu}^{1\,I},
 \end{equation}
where
 \begin{equation}\label{eq:F1I}
F^{1I}=\hat{Z}^I {F}-\hat{Z}^I F^0-D^{I}, \qquad
 D^I=\big(\hat{Z}^I \Boxr_g h^1-\Boxr_g Z^I h^1\big)
 \end{equation}

 We base our argument on the
energy estimate \eqref{eq:firstenergy} for a solution $\phi$ of the
wave equation $\Boxr_{g}\phi =F$ of Proposition
\ref{prop:Decayenergy}. Observe that the conditions of our
Theorem on the  tensor $H^{\mu\nu}=g^{\mu\nu}-m^{\mu\nu}$
imply the that the assumptions of
Proposition \ref{prop:Decayenergy} for the metric $g$ hold true.
In particular, we have
\beq\label{eq:energyuse}
 \int_{\Si_{t}} |\pa \phi|^{2}\,w + \int_{0}^{t} \int_{\Si_{\tau}}
|\pab\phi|^{2}\,w^{\,\prime} \leq
 8\int_{\Si_{0}} |\pa \phi|^{2}\,w+
16\int_0^t\int_{\Si_{t}} \Big(\frac{C\varepsilon
\,|\pa\phi|^{2}}{1+t}
 + |\Boxr_g \phi|\, |\pa\phi|\Big)\, w
\eq
 This applied to \eqref{eq:h1Iwaveeq} gives
\begin{multline}\label{eq:energyuseh1I}
 \int_{\Si_{t}} |\pa Z^I h^1|^{2}\,w + \int_{0}^{t} \int_{\Si_{\tau}}
|\pab Z^I h^1|^{2}\,w^{\,\prime} \leq
 8\int_{\Si_{0}} |\pa h^1|^{2}\,w+
16\int_0^t\int_{\Si_{t}} \Big(\frac{C\varepsilon \,|\pa Z^I
h^1|^{2}}{1+t}
 + |F_1^I|\, |\pa Z^I h^1|\Big)\, w\\
 \leq
  8\int_{\Si_{0}} |\pa h^1|^{2}\,w+
16\int_0^t\int_{\Si_{t}}\Big(\frac{C\varepsilon \,|\pa Z^I
h^1|^{2}}{1+t} \, w
 +\ve^{-1}
 \big(|\hat{Z}^I F|^2+\big|D^I\big|^2\big)\, (1+t)\,w\, + |Z^I F^{\,0}|\,|\pa Z^I h^1|\,
w \Big)
\end{multline}
We begin with the following estimate on the inhomogeneous term
$F$:
\begin{lemma}\label{lem:F}
Under the assumptions of Theorem \ref{energyest}
\begin{align}
 |Z^I F|&\les \sum_{|J|\le |I|}\Big(\frac{\ve |\pa Z^J h_1 |}{1+t} +
 \frac{\ve(1+|q|)^{\mu^\prime-1/2}}{(1+t+|q|)^{1-C\ve}}\, |\pab Z^J h_1
 |+\frac{\ve^2}{1+t+|q|}\, \frac{|Z^J h_1
 |}{1+|q|}\Big)\label{eq:Fest}
 \\ &+\sum_{|J|\le |I|-1}\frac{\ve\,|\pa Z^J h_1|}{(1+t)^{1-C\ve}}+\frac{\ve^2}{(1+t+|q|)^4}\nn
\end{align}
\end{lemma}
\begin{proof} According to  Proposition \ref{decayinhom}
\begin{align}
|Z^I F|& \les
 \sum_{|J|+|K|\le |I|}
\big (|\pa Z^J h|_{\cal TU}\, |\pa Z^K h |_{\cal TU} + |\pab Z^J
h|\, |\pa Z^K h| \big )+ \sum_{|J|+|K|\le |I|-1} |\pa Z^J h|\,
|\pa Z^K h|\nn\\ & + \sum_{|J_1|+|J_2|+|J_3| \leq |I|}
|Z^{J_3} h|\,|\pa Z^{J_{2}} h| |\pa
Z^{J_1} h|\label{eq:inhomogeneousestimate}
\end{align}
Since $h=h^1+h^0$ and  $h^0$ obeys the estimates $|\pa Z^J h^0|\les
\ve (1+t+|q|)^{-2}$ and $|Z^I h^0|\les \ve (1+t+|q|)^{-1}$ for all
$|I|\le N$, it follows that both $h^0$ and $h^1$ obey the
estimates \eqref{eq:specialcomp} and \eqref{eq:generalcomp} of the
theorem. Substituting $h=h^0+h^1$ into \eqref{eq:inhomogeneousestimate}
and expanding the products
we obtain terms in which either all factors contain $h^0$,
in which case we can estimate them by $\ve^2
(1+t+|q|)^{-4}$, or at least one factor contains $h^1$, in which case we
can simply estimate the other factors by  \eqref{eq:specialcomp}-\eqref{eq:generalcomp}.
It follows that
\eqref{eq:inhomogeneousestimate} leads to the estimate
\begin{align*}
 |Z^I F|&\les \sum_{|J|\le |I|}\Big(\frac{\ve |\pa Z^J h^1 |}{1+t} +
 \frac{\ve(1+|q|)^{\mu^\prime-1/2}}{(1+t+|q|)^{1-C\ve}}\, |\pab Z^J h^1
 |+\frac{\ve^2}{1+t+|q|}\, \frac{|Z^J h^1
 |}{1+|q|}\Big)\\
 &+\sum_{|J|\le |I|-1}\frac{\ve\,|\pa Z^J h^1|}{(1+t)^{1-C\ve}}+\frac{\ve^2}{(1+t+|q|)^4}
\end{align*}
\end{proof}
\begin{lemma}\label{lem:Fest}
Under the assumptions of Theorem \ref{energyest}
\begin{align*}
\ve^{-1} \int_0^T\int |Z^I F|^2\, (1+t)\,w\,  dx dt
 &\les \sum_{|J|\le |I|}\int_0^T\int \ve \Big(\frac{|\pa Z^J h^1 |^2}{1+t}\, w +
 |\pab Z^J h^1|^2 \, w^\prime\Big)dx dt\\
 &+\sum_{|J|\le |I|-1}\int_0^T\int
 \ve \frac{\,|\pa Z^J h^1|^2}{(1+t)^{1-2C\ve}}
 \, w\, dx dt +\ve^3.
\end{align*}
\end{lemma}
\begin{proof} The estimate  is a straightforward application of
Lemma \ref{lem:F}. We took into account that
$w\le w'(1+|q|) (1+q_-)^{2\mu}$ and the inequality $\mu<1-2\mu'$.
The estimate
$$
\int\frac{1}{1+t+|q|}\, \frac{|Z^J h^1|^2}{(1+|q|)^2}\,w\,dx \les
\int\frac{|\pa Z^J h^1|^2}{(1+t+|q|)}\,w\,dx
$$
is the Hardy type inequality established in Corollary
\ref{cor:Poinc} of Appendix B.
\end{proof}
Next we estimate $F^0=\Boxr_g h^0$:
\begin{lemma} \label{lem:Festiv}
The following inequality holds true:
\begin{equation*}
\int_0^T\int |Z^I F^{\,0}|\,|\pa Z^I h^1| w\, dx dt\leq C_N\,
\ve\!\! \sum_{|J|\leq |I|}\Big(\int_0^T\int \frac{|\pa Z^J
h^1|^2}{(1+t)^2} w \, dx dt+ \int_0^T\!\!\Big(\int|\pa Z^J h^1|^2
w\, dx \Big)^{1/2}\!\!\!\frac{ dt}{(1+t)^{3/2}}\Big)
\end{equation*}
where $C_N$ denotes a constant that depends only on $N$.
\end{lemma}
\begin{proof}
We start with the estimate
$$
\int_0^T\int |Z^I F^{\,0}|\,|\pa Z^I h^1| w\, dx dt\le
\int_0^T\Big (\int |Z^I F^{\,0}|^2 w\, dx \Big )^{\frac 12}\,
\Big (\int |\pa Z^I h^1|^2 w\, dx \Big )^{\frac 12} dt
$$
Now by Lemma \ref{subtractoffschwarz}
$$
|Z^I F^{\,0}|\leq \begin{cases} C_N\varepsilon^2(t+|q|+1)^{-4}
,\quad q>0,\\
C_N\varepsilon(t+|q|+1)^{-3} ,\quad q<0,
\end{cases}+\frac{C_N\varepsilon }{(t+|q|+1)^3}
\sum_{|J|\leq |I|}
  |Z^{J}h^1|.
$$
Therefore, once again using a Hardy type inequality of
Corollary \ref{cor:Poinc}, we obtain
\begin{align*}
\int |Z^I F^{\,0}|^2 w\, dx &\leq C_N\int_0^t \frac{\ve^2 r^2
dr}{(1+t)^6} +C_N \int_t^\infty \frac{\ve^4 r^2
(1+r)^{1+2\gamma}\,
dr}{(1+r)^8}+\frac{C_N\ve^2}{(1+t)^4}\int \frac{|Z^I h_1|^2}{(1+|q|)^2}\, w\, dx\\
&\leq \frac{C_N \ve^2}{(1+t)^3}+ \frac{C_N \ve^2}{(1+t)^4}\int|\pa
Z^I h^1|^2 w\, dx
\end{align*}
and the result follows.
\end{proof}
The estimate for the term containing the commutator term $D^I=\Boxr_g Z^I h^1-\hat{Z}^I\Boxr_g h^1$ follows from the following
 \begin{lemma} \label{lem:commutatorenergy}
 Under the assumptions of Theorem \ref{energyest}
\begin{multline}
\ve^{-1}\int_0^T\int \big|\Boxr_g
Z^I h^1-\hat{Z}^I\Boxr_g h^1\big|^2 (1+t)\,w \, dx \, dt\\
\les \ve\sum_{|J|\leq |I|}\int_0^T\int\Big( \frac{|\pa Z^J
h^1|^2}{1+t} w\,
 +\, |\pab Z^J h^1|^2 w^\prime \Big)\, dx dt
 +\ve\sum_{|J|\leq |I|-1}\int_0^T\int \frac{|\pa Z^J
 h^1|^2}{(1+t)^{1-2C\ve}}w
 \, dx dt+\ve^3
\end{multline}
\end{lemma}
We postpone the proof of Lemma \ref{lem:commutatorenergy} for a moment
and finish the proof of Theorem \ref{energyest}.

Using \eqref{eq:energyuseh1I} together with Lemmas \ref{lem:Fest}-\ref{lem:commutatorenergy} yields
\begin{multline}
 \int_{\Si_{t}} |\pa Z^I h^1|^{2}\,w + \int_{0}^{t} \int_{\Si_{\tau}}
|\pab Z^I h^1|^{2}\,w^{\,\prime}\leq
  8\int_{\Si_{0}} |\pa h^1|^{2}\,w
  +C_N\ve\sum_{|J|\leq |I|}
  \int_0^T\frac{1}{(1+t)^{3/2}}\Big(\int|\pa Z^J h^1|^2 w\, dx
\Big)^{1/2}\!\! dt\\
  +C\ve\sum_{|J|\leq |I|}\int_0^T\int\Big( \frac{|\pa Z^J
h^1|^2}{1+t} w\,
 +\, |\pab Z^J h^1|^2 w^\prime \Big)\, dx dt
 +C\ve\sum_{|J|\leq |I|-1}\int_0^T\int \frac{|\pa Z^J
 h^1|^2}{(1+t)^{1-2C\ve}}w
 \, dx dt+C\ve^3
\end{multline}
where $C_N$ depends only on $|I|\leq N$.
As before we denote
$$
\E_k(t)=\sup_{0\le \tau\le t} \sum_{Z\in {\cal Z},\,|I|\le k}
\int_{\Si_\tau} |\pa Z^I h^1|^2 w\, dx
 $$
 and let
 $$
 S_k(t):=\sum_{Z\in {\cal Z},\,|I|\le k}
\int_0^t\int_{\Si_\tau} |\pab Z^I h^1|^2 w'\, dx
 $$
 It therefore follows that
\begin{equation}\label{eq:finalenergyestimate}
E_k(t)+S_k(t)\leq 8 E_k(0)+C\ve S_k(t)+\int_0^t \frac{C\ve
E_k(\tau)}{1+\tau} d\tau +\int_0^t \frac{C_N\ve
E_k(\tau)^{1/2}}{(1+\tau)^{3/2}} d\tau + C\ve^3+\int_0^t
\frac{C\ve E_{k-1}(\tau)}{(1+\tau)^{1-C\ve}} d\tau
\end{equation}
For $C\ve$
sufficiently small we can absorb the space-time integral $S_k(t)$
into the one on the left hand-side at the expense of at most doubling
all the constants on the right hand-side.  Similarly,  since
$16 C_N \ve 1/4E_k^{1/2}\leq E_k+ 64^2 C_N^2 \ve^2$ and $E_k(t)$ is increasing,
we can absorb $1/4\int_0^t E_k(\tau)\,(1+\tau)^{-3/2}\leq 1/2E_k(t)$.
 If we also
use the assumption $E_N(0)\leq \ve^2$ we
obtain for $\ve>0$ sufficiently small
\begin{equation}\label{eq:finalenergyestimate}
E_k(t)+S_k(t)\leq C_N \ve^2 +\int_0^t \frac{C\ve
E_k(\tau)}{1+\tau} d\tau +\int_0^t \frac{C\ve
E_{k-1}(\tau)}{(1+\tau)^{1-C\ve}} d\tau
\end{equation}
 where the last term is absent if $k=0$ and $C_N$ is a constant
 that depending only $N$.

For $k=0$ this yields the estimate
 $$
 E_0(t)\leq C_N \ve^2+ \int_0^t \frac{c_0\ve
E_0(\tau)}{1+\tau} d\tau
 $$
 and the Gronwall inequality gives the bound
 $$
 E_0(t)\leq C_N(1+t)^{c_0\ve}
 $$
which prove \eqref{eq:energyestk} for $k=0$.

Assuming \eqref{eq:energyestk} for $k$ replaced by $k-1$ we get
from \eqref{eq:finalenergyestimate}
\begin{equation}\label{eq:finalenergyestimate'}
E_k(t)\leq C_N \ve^2 +\int_0^t \frac{c\ve E_k(\tau)}{1+\tau} d\tau
+\int_0^t \frac{c\ve^3 d\tau  }{(1+\tau)^{1-c\ve}}
\end{equation}
which leads to the bound
$$
E_{k}(t)\leq C_N\ve^2 (1+t)^{2c\ve}
$$
This concludes the induction and the proof of the theorem.
\end{proof}
\subsection {Proof of Lemma \ref{lem:commutatorenergy}}
\begin{proof}
Define the tensor
$$
H_1^{\mu\nu}:=H^{\mu\nu} - H_0^{\mu\nu}, \qquad H_0^{\mu\nu}=
-\chi(\frac rt)\chi(r) \frac Mr \delta^{\mu\nu}
$$
Observe that $H_0$ coincides with the tensor $-h^0$.
 Define also the wave operator $\Boxr_1=\Box+H_1^{\alpha\beta}\pa_\alpha\pa_\beta$.
 According to \eqref{eq:curvwaveeqcommutest2-5} of Proposition
\ref{prop:commut3}
 \begin{align}
\big|\Boxr_1 Z^I h^1&-\hat{Z}^I\Boxr_1 h^1\big|\les \sum_{|K|\le
|I|}\,\,\,\,\sum_{|J|+(|K|-1)_+\leq |I|} \Big (\frac{|Z^J
H_1|}{1+t+|q|}
+\frac{|Z^J H_1|_{\cal LL }}{1+|q|}\Big )|\pa Z^K h^1|\label{eq:commutatorH1}\\
&+ \sum_{|K|\le |I|}\Big (\sum_{|J|+(|K|-1)_+\leq |I|-1} \frac{|Z^J
H_1|_{\cal LT}}{1+|q|}+\sum_{|J|+(|K|-1)_+\leq |I|-2} \frac{|Z^{J}
H_1|}{1+|q|}\Big)|\pa Z^K h^1|.\nn
 \end{align}
 Our goal is to obtain the estimate for the quantity
$$
\sum_{|I|\le N} \int_0^T \int \big|\Boxr_1 Z^I h^1-\hat{Z}^I\Boxr_1
h^1\big|^2 (1+t) w\, dx dt
$$
Let us first deal with the terms in \eqref{eq:commutatorH1} with
$|K|\leq N/2+1$. In this case we use the decay estimates
\eqref{eq:generalcomp}
$$
| \pa Z^K h|\le
C\varepsilon
(1+t+|q|)^{-1+C\varepsilon}(1+|q|)^{-1/2+\mu^\prime} \le
C\varepsilon
(1+t+|q|)^{-1+C\varepsilon}(1+|q|)^{-C\ve -\mu}
$$
guaranteed by the assumptions of Theorem \ref{energyest}, provided
that $\mu<1/2-\mu'$. It is clear that in this case it suffices to consider
the expression
\begin{multline}
\int_0^T \int \Big(\frac{|Z^J
H_1|^2}{(1+t+|q|)^2}+ \frac{|Z^J H_1|_{\cal
LL}^2+ |Z^{J'} H_1|_{\cal LT}^2 + |Z^K H_1|^2}{(1+|q|)^2}\Big)\frac{\ve^2
(1+|q|)^{-2C\ve}}{(1+t+|q|)^{1-2C\ve}}
\frac{w\, dx dt}{(1+q_-)^{2\mu}}\\
\les  \int_0^T \int \frac{|Z^J
H_1|^2}{(1+|q|)^2}\frac{\ve^2\, w\,dx dt}{1+t}+  \int_0^T \int\frac{|Z^K
H_1|^2}{(1+|q|)^2}\frac{\ve^2 w\,dx dt }{(1+t)^{1-2C\ve}}
\\+  \int_0^T \int\frac{|Z^J H_1|_{\cal LL}^2+ |Z^{J'} H_1|_{\cal LT}^2}{(1+|q|)^2}\frac{\ve^2(1+|q|)^{-2C\ve}}{(1+t+|q|)^{1-2C\ve}}
\frac{w\,dx dt}{(1+q_-)^{2\mu}}
\end{multline}
with $|J|\le k, |J'|\le k-1, |K|\le k-2$, where $k=|I|$, After
applying the Hardy type inequalities of Corollary \ref{cor:Poinc}
the above expression is bounded by
$$
\ve^2\Big (\int_0^T \int \frac{|\pa Z^J
H_1|^2}{1+t}\, w\,dx dt +  \int_0^T \int \frac{|\pa Z^K
H_1|^2}{(1+t)^{1-2C\ve}}\, w\,dx dt
+  \int_0^T \int (|\pa Z^J H_1|_{\cal LL}^2+ |\pa Z^{J'} H_1|_{\cal LT}^2)  \, {\tilde w\,dx dt}\Big ),
$$
where
$$
\tilde w=\min\, (w', \frac{w}{(1+t+|q|)^{1-2C\eps}}).
$$
Ignoring the difference, which we shall comment on at the end of the
proof, between the tensors
$H_1^{\mu\nu}=g^{\mu\nu}-m^{\mu\nu}-H_0^{\mu\nu}$ and
$h^1_{\mu\nu}=g_{\mu\nu} - m_{\nu\mu} - h^0_{\mu\nu}$, we see that
the first two terms are as claimed in the statement of the lemma. We
now recall that according to Lemma \ref{lem:Wave-H1} of Appendix D
\begin{multline}
\sum_{|J|\le k} \big| \pa Z^J H_1\big|_{\cal LL} + \sum_{|J|\le k-1}
\big| \pa Z^J H_1\big|_{\cal LT} \les\sum_{|J|\leq k} |\pab Z^{J}
H_1|+\sum_{|J'|\le k-1} |\pa Z^{J'} H_1|\nn\\+ \ve\frac
1{1+t+|q|}\sum_{|J|\leq k} \Big( |\pa Z^{J} H_1|+\frac{|Z^{J}
H_1|}{1+t+|q|}\Big)\nn \\
+ \sum_{|J_1|+|J_2|\le k} |Z^{J_1}H_1| |\pa
Z^{J_2} H_1| +\frac{C\ve\chi_0(1/2<r/t<3/4)}{(1+t+|q|)^2}+\frac{C\ve^2}{(1+|t|+|q|)^3}\nn,
\end{multline}
It is clear that in the sum $\sum_{|J_1|+|J_2|\le k}$ above at least one of the
indices is $\le N/2$ and therefore we can use the decay estimates
\eqref{eq:generalcomp}: for $|J'|\le N/2+2$
$$
\frac {| Z^{J'} h|}{1+|q|}+ | \pa Z^{J'} h|\le
C\varepsilon
(1+t+|q|)^{-1+C\varepsilon}(1+|q|)^{-1/2+\mu^\prime} \le
C\varepsilon
(1+t+|q|)^{-1+C\varepsilon}(1+|q|)^{-C\ve}.
$$
Thus, with $|J|\le k, |J'|\le k-1$,
 \begin{align}
\int_0^T \int\Big (|\pa Z^J H_1|_{\cal LL}^2 + |\pa Z^{J'}
H_1|_{\cal LT}^2\Big)\, &{\tilde w\,dx dt}\les \int_0^T
\int \sum_{|J|\leq k}|\pab Z^{J} H_1|^2 w^\prime\,
dx dt\label{eq:neck} \\+ \ve^2\int_0^T\int
\sum_{|J|\leq k} \Big( |\pa Z^{J} H_1|^2&+\frac{|Z^{J}
H_1|^2}{(1+|q|)^2}\Big)
\Big(\frac{1+|q|}{1+t+|q|}\Big)^{2-2C\ve}\, w^\prime\, dx dt \,+\nn\\
\sum_{|J'|\leq k-1}\int_0^T \int
\frac {|\pa Z^{J'} H_1|^2}{(1+t)^{1-2C\eps}} w\, dx dt+\int_0^T\int&\Big(\frac{C\ve^2\chi_0(1/2<r/t<3/4)}{(1+t+|q|)^4}
+\frac{C\ve^4}{(1+|t|+|q|)^6}\Big)\, w^\prime\, dx dt\nn
\end{align}
where $\chi_0(1/2<r/t<3/4)$ is the characteristic function of the
set where $t/2<r<3t/4$. Using the properties of the function $w'$, in particular that
 $w^\prime\les w/(1+|q|)$ we obtain that the above has a bound
\begin{align*}
\eqref{eq:neck}\,\,\les\,\,\sum_{|J|\leq k} \int_0^T \int |\pab Z^{J} H_1|^2
w^\prime\, dx dt + \ve^2\int_0^T\int \sum_{|J|\leq k} \Big( |\pa
Z^{J} H_1|^2+\frac{|Z^{J} H_1|^2}{(1+|q|)^2}\Big)
\frac{ w\, dx dt}{1+t}&\\+\sum_{|J'|\leq k-1}\int_0^T \int
\frac {|\pa Z^{J'} H_1|^2}{(1+t)^{1-2C\eps}} w\, dx dt
+\int_0^T\int\frac{C\ve^2\chi^2_0(1/2<r/t<3/4)}{(1+t+|q|)^4(1+|q|)^{1+2\mu}}\,
dx dt &\\\hskip -2pc+\int_0^T\frac{C\ve^4}{(1+|t|+|q|)^6}\,
(1+|q|)^{1+2\gamma}\, dx dt&
\end{align*}
Once again we use the  Hardy type inequality of Corollary \ref{cor:Poinc}
with $a=0$ to conclude that
\begin{align*}
\int_0^T \int\Big (|\pa Z^I H_1|_{\cal LL}^2 + |\pa Z^J H_1|_{\cal
LT}^2\Big)\, {w'\,dx dt}&\les\sum_{|J|\leq k} \int_0^T \int
\Big(|\pab Z^{J} H_1|^2 \, w^\prime\,+ \frac{|\pa Z^{J}
H_1|^2}{1+t}\, w\Big)dx dt \\ &+\sum_{|J'|\leq k-1}\int_0^T \int
\frac {|\pa Z^{J'} H_1|^2}{(1+t)^{1-2C\eps}} w\, dx dt+C\ve^2
\end{align*}
as desired, modulo the difference between $H_1$ and $h^1$, in Lemma \ref{lem:commutatorenergy}.
To finish the proof of the Lemma for the case $|K|\le N/2+1$ it remains
address the difference between the tensors
$$
H_1^{\mu\nu} = g^{\mu\nu} - m^{\mu\nu} + {h^0}^{\mu\nu}
\qquad {\text and }\quad
h^1_{\mu\nu} = g_{\mu\nu} - m_{\mu\nu} - {h^0}_{\mu\nu}
$$
The above expressions imply that
$$
(m+h^0+h^1)_{\mu\nu}^{-1}=(m-h^0+ H_1)^{\mu\nu}
$$
Therefore,
$
H_1^{\mu\nu} = - {h^1}^{\mu\nu} + O^{\mu\nu}\Big ((h^0+h^1)^2)$ and
it follows that
\begin{align*}
&|\pa Z^J H_1|\les |\pa Z^J h^1| +\sum_{|J_1|+|J_2|\le |J|}
\Big (|Z^{J_1} h^0| |\pa Z^{J_2} h^1| + |Z^{J_1} h^1| |\pa Z^{J_2}h^1| +
|\pa Z^{J_1} h^0| |Z^{J_2}h^1| + |Z^{J_1} h^0| |\pa Z^{J_2} h^0| \Big ),\\
&|\pab Z^J H_1|\les |\pab Z^J h^1| +\sum_{|J_1|+|J_2|\le |J|}
\Big (|Z^{J_1} h^0| |\pa Z^{J_2} h^1| + |Z^{J_1} h^1| |\pa Z^{J_2} h^1| +
|\pa Z^{J_1} h^0| |Z^{J_2} h^1| + |Z^{J_1} h^0| |\pa Z^{J_2} h^0| \Big )
\end{align*}
Taking into account that
for  $i=0,1$ and $|J'|\le N/2+2$
$$
\frac {| Z^{J'} h^i|}{1+|q|}+ | \pa Z^{J'} h|\le
C\varepsilon
(1+t+|q|)^{-1+C\varepsilon}(1+|q|)^{-1/2+\mu^\prime} \le
C\varepsilon
(1+t+|q|)^{-1+C\varepsilon}(1+|q|)^{-C\ve}.
$$
we obtain that
\begin{align*}
&|\pa Z^J H_1|\les |\pa Z^J h^1| +\Big (\frac {1+|q|}{1+t+|q|}\Big)^{1-C\ve}
\sum_{|K|\le |J|}
 \Big (|\pa Z^K h^1|+ \frac {|Z^K h^1|}{1+|q|}\Big )  + \frac {\ve^2}{(1+t+|q|)^3},\\
&|\pab Z^J H_1|\les |\pab Z^J h^1| +\Big (\frac {1+|q|}{1+t+|q|}\Big)^{1-C\ve}
\sum_{|K|\le |J|}
 \Big (|\pa Z^K h^1|+ \frac {|Z^K h^1|}{1+|q|}\Big )  + \frac {\ve^2}{(1+t+|q|)^3}
 \end{align*}
Thus, with the help of the inequality $w'\le w/(1+|q|)$,
\begin{align*}
&\sum_{|J|\le k, |K|\le k-1} \int_0^T \int \Big (\frac{|\pa Z^J
H_1|^2}{1+t}\, w\,dx dt +\frac{|\pa Z^K H_1|^2}{(1+t)^{1-2C\ve}}\Big
)\, w\,dx dt+ \int_0^T \int |\pab Z^J H_1|^2\, w'\,dx dt \\ &\les
\sum_{|J|\le k, |K|\le k-1} \int_0^T \int \Big (\frac{|\pa Z^J
h^1|^2}{1+t}\, w\,dx dt +\frac{|\pa Z^K h^1|^2}{(1+t)^{1-2C\ve}}\Big
)\, w\,dx dt+ \int_0^T \int |\pab Z^J h^1|^2\, w'\,dx dt \\
&\qquad\qquad +\int_0^T \int \frac 1{(1+t)}\frac{|Z^J
h^1|^2}{(1+|q|)^2}\, w\,dx dt + \int_0^T \int \frac
1{(1+t)^{1-2C\ve}}\frac{|Z^K h^1|^2}{(1+|q|)^2}\, w\,dx dt +\ve^2\\
&\les \sum_{|J|\le k, |K|\le k-1} \int_0^T \int \Big (\frac{|\pa Z^J
h^1|^2}{1+t}\, w\,dx dt +\frac{|\pa Z^K h^1|^2}{(1+t)^{1-2C\ve}}\Big
)\, w\,dx dt+ \int_0^T \int |\pab Z^J h^1|^2\, w'\,dx dt +\ve^2,
\end{align*}
where to pass to the last inequality we once again used the Hardy type inequality
of Corollary \ref{cor:Poinc}.

Returning to \eqref{eq:commutatorH1} we now deal with the case $|K|\ge N/2$,
which implies that $|J|\le N/2+1$ and allows us to use the decay estimates
\eqref{eq:specialcomp}-\eqref{eq:generalcomp} for $H_1=-h^1 + O(h^2)$.
Therefore, the contribution of the terms with $|K|\ge N/2$ to
$\big|\Boxr_1 Z^I h^1-\hat{Z}^I\Boxr_1 h^1\big|$ can be bounded by
\begin{align*}
\sum_{|K|=|I|}\sum_{|J|=1} &\Big (\frac{|Z^J H_1|}{1+t+|q|} +
\frac{|Z^J H_1|_{\cal LL} + |H_1|_{\cal LT}}{1+|q|}\Big ) |\pa Z^K h^1|+
\sum_{|K|<|I|} \frac{|Z^J H_1|}{1+|q|}|\pa Z^K h^1|\\ &\les
\sum_{|K|=|I|} \frac {|\pa Z^K h^1|}{1+t} +
\sum_{|K|<|I|} \frac{|\pa Z^K h^1|}{(1+t)^{1-C\ve}}
\end{align*}
and the desired result follows.

To estimate the commutator $\widetilde\Box_g Z^I h^1 -
\hat Z^I\widetilde \Box_g h^1$ it remains to address the term
$$
|H_0^{\a\b}  \pa_\a\pa_\b Z^I h^1 -
\hat Z^I \big (H_0^{\a\b}  \pa_\a\pa_\b Z^I h^1 \big )|\le
\ve \frac 1{1+t+|q|} \sum_{|J|\le |I|} C_J \frac {|\pa Z^J h^1|}{1+|q|}
$$
Therefore
\begin{align*}
\int_0^T \int |H_0^{\a\b}  \pa_\a\pa_\b Z^I h^1 - \hat Z^I \big
(H_0^{\a\b}  \pa_\a\pa_\b Z^I h^1\big )|^2 (1+t) \,w\,dxdt \les
\ve^2\sum_{|J|\le |I|} \int_0^T \int  \frac{|\pa Z^J h^1|^2}{1+t}
\,w\,dxdt
\end{align*}

\end{proof}

\section{Appendix A. Commutators}
Recall the family of vector fields ${\cal Z}=\{\pa_\a,
\Omega_{\a\b}=-x_\a\pa_b+x_\b\pa_\a, S=t\pa_t + r\pa_r\}$. In this
section we address commutation properties of the family ${\cal Z}$
with various differential structures. Recall that for any $Z\in
{\cal Z}$ we have $[Z,\Box]=-c_Z \Box$, where $c_Z$ is different
from zero only for the scaling vector field $S$, $c_S=2$.
\begin{lemma} \label{commut1}
Let $Z\in {\cal Z}$ and let the constants
$c_\alpha^\mu$ be defined by
$$
[\pa_\alpha,Z]=c_\alpha^{\,\,\, \mu}\pa_\mu ,\qquad
c_\alpha^{\,\,\, \mu}=\pa_\alpha Z^\mu
$$
Then  $c_{LL}=c^{\Lb\Lb}=0$.
In addition, if $Q$ is a null form, then
\beq\label{eq:nullformcommut} Z Q (\pa \phi,\pa \psi)=Q(\pa
\phi,\pa Z \psi)+Q(\pa Z \phi,\pa \psi) + \tilde Q (\pa \phi,\pa
\psi) \eq for some null form $\tilde Q$ on the right hand-side.
\end{lemma}
\begin{proof}
Since $Z=Z^\alpha\pa_\alpha$ is a Killing or conformally Killing
vector field we have \beq\label{eq:confkill} \pa_{\a} Z_{\beta} +
\pa_{\beta} Z_{\alpha} = f m_{\a\b} \eq where
$Z_\alpha=m_{\alpha\beta} Z^\beta$. In fact, for the vector fields
above, $f=0$ unless $Z=S$ in which case $f=2$. In particular,
$$
L^{\alpha} L^{\beta}\pa_\alpha Z_\beta =0.
$$
If $c_{\alpha}^{\,\,\,\mu}$ is  as defined above and
$c_{\alpha\beta} = c_{\alpha}^{\,\,\,\mu}m_{\mu\beta}=\pa_\alpha
Z_\beta$ the above simply means that $c_{LL}=c^{\Lb\Lb}=0$. which
proves the first part of the lemma. To verify
\eqref{eq:nullformcommut} we first consider the null form
$Q=Q_{\a\b}$ We have
\begin{align*}
Z Q_{\a\b}(\pa\phi,\pa\psi) &= Q_{\a\b}(\pa Z \phi,\pa\psi) +
 Q_{\a\b}(\pa\phi,\pa Z \psi) \\ &+ [Z,\pa_\a] \phi \pa_\b \psi -
\pa_\b \phi [Z,\pa_\a]\psi + [Z,\pa_\b] \phi \pa_\a \psi - \pa_\a
\phi [Z,\pa_\b]\psi\\ & =
 Q_{\a\b}(\pa Z \phi,\pa\psi) +
 Q_{\a\b}(\pa\phi,\pa Z \psi) - c_\a^\mu ( \pa_\mu \phi \pa_\b \psi
- \pa_\b\phi \pa_\mu \psi ) - c_\b^\mu ( \pa_\mu \phi \pa_\a \psi
- \pa_\a\phi \pa_\mu \psi)\\ & =
 Q_{\a\b}(\pa Z \phi,\pa\psi) +
 Q_{\a\b}(\pa\phi,\pa Z \psi) - c_\a^\mu Q_{\mu\b} (\pa\phi,\pa\psi)-
c_\b^\mu Q_{\mu\a} (\pa\phi,\pa\psi)
\end{align*}
The calculation for the null form $Q_0(\pa\phi,\pa\psi) = m^{\a\b}
\pa_\a\phi \pa_\b \psi$ is similar and we leave it to the reader.
\end{proof}
For any symmetric 2-tensor $\pi$ and a vector field $Z\in{\cal Z}$
define \beq
\pi_Z^{\alpha\beta}\pa_\alpha\pa_\beta=\pi^{\alpha\beta}[\pa_\alpha\pa_\beta,
Z], \qquad\text{i.e.,}\qquad \pi_Z^{\alpha\beta}:=\pi^{\alpha\gamma}
c_{\gamma}^{\,\,\, \beta}
+\pi^{\gamma\beta}c_{\gamma}^{\,\,\,\alpha}. \eq
\begin{lemma} \label{commut2}
The tensor coefficients $\pi_Z^{\a\b}$ verify the following estimate:
\beq \label{eq:kLT}
|\pi_Z|_{LL}\leq 2|\pi|_{L\cal T}. \eq
In general,
\begin{align}
&[\pi^{\alpha\beta}\pa_\alpha\pa_\beta,
Z^I]=\sum_{I_1+I_2=I, \,|I_2|<|I|}
\pi^{I_1\alpha\beta}\pa_\alpha\pa_\beta Z^{I_2}, \label{eq:kformula}\\
&\pi^{J\alpha\beta}:= \!\sum_{ |K|\leq |J|}c^{J\alpha\beta}_{K\mu\nu}
Z^K (\pi^{\mu\nu})= -Z^J (\pi^{\alpha\beta}) -\!\sum_{K+Z=J} Z^K
\pi_Z^{\alpha\beta} +\!\!\sum_{ |K|\leq |J|-2}
d^{J\alpha\beta}_{K\mu\nu} Z^K (\pi^{\mu\nu}) \nn
\end{align}
for some constants
$c^{J\alpha\beta}_{M\mu\nu}$ and $d^{J\alpha\beta}_{M\mu\nu}$.
Here the sum \eqref{eq:kformula} means the sum over all possible
order preserving partitions of the multi-index $I$ into
multi-indices $I_1$, $I_2$.
\end{lemma}
\begin{proof}
First observe that since the vector fields $Z$ are linear in $t$
and $x$ we have
$$
[\pa^{2}_{\a\b}, Z] = [\pa_{\b}, Z] \pa_{\a} + [\pa_{\a}, Z]
\pa_{\b}= c_{\beta}^{\,\,\gamma}\pa_\gamma\pa_\alpha
+c_{\alpha}^{\,\,\,\gamma} \pa_\gamma\pa_\beta,
$$
which proves the first statement, while the second follows since
$c_{L}^{\,\,\Lb} =0$.

To prove \eqref{eq:kformula} we first write
$$
Z^I \big(\pi^{\alpha\beta}\pa_\alpha\pa_\beta \phi\big)
=\sum_{K+J=I} (Z^K \pi^{\alpha\beta}) Z^J \big( \pa_\alpha\pa_\beta
\phi\big)
$$
Then we observe that \beq\label{eq:Zcom2} Z^{J}
\pa_\alpha\pa_\beta\phi=\sum_{J_1+J_2=J,\,
J_1=(\iota_{\ell_1},...,\iota_{\ell_n})}
\left[Z^{\iota_{\ell_1}},\left[Z^{\iota_{\ell_2}},\left[...,
\left[Z^{\iota_{\ell_{n-1}}},[Z^{\iota_{\ell_n}},
\pa^2_{\a\b}]\,\right]...\right] \,\right]\,\right] Z^{J_2}\phi,
\eq where the sum is over all order preserving partitions of $(1,...,k)$ into two ordered
sequences $(\ell_1,...,\ell_n)$ and
$(\ell_{n+1},...,\ell_k)$ such that $J_2=(\iota_{\ell_{n+1}},...,\iota_{\ell_k})$.
It therefore follows that
$$
\pi^{J\a\b} = - \sum_{K+L=J,\, L=(\iota_1,...,\iota_l)} (Z^K
\pi^{\a\b}) \left[Z^{\iota_1},\left[Z^{\iota_2},\left[...,
\left[Z^{\iota_{l-1}},[Z^{\iota_l},
\pa^2_{\a\b}]\,\right]...\right] \,\right]\,\right]
$$
The desired representation follows after taking into account that
$$
(Z^K \pi^{\a\b}) [Z,\pa^2_{\a\b}] = - (Z^K
\pi^{\a\b}_Z)\pa_\alpha\pa_\beta
$$
\end{proof}
For a symmetric 2-tensor $H$ and a vector field $Z\in {\cal Z}$ we
set $\hat Z= Z+ c_Z$, where $c_Z$ is the constant in the commutator
$[Z,\Box]=-c_Z \Box$, and \beq \hat H^{J\alpha\beta}= \!\sum_{
|M|\leq |J|} {c}^{J\alpha\beta}_{M\mu\nu} \hat{Z}^M H^{\mu\nu}=
-\hat Z^J H^{\alpha\beta} - \!\sum_{M+Z=J } \hat Z^M
H_{Z}^{\alpha\beta} +\!\!\sum_{ |M|\leq |J|-2}
{d}^{J\alpha\beta}_{M\mu\nu} \hat Z^M H^{\mu\nu} \eq
\begin{cor}\label{commut3}
 Let $\Boxr_g=\Box+H^{\alpha\beta}\pa_\alpha\pa_\beta$.
 Then
\begin{align}
&\Boxr_g Z \phi- \hat{Z}\Boxr_g
\phi= -({\hat Z} H^{\alpha\beta}+
H_{{Z}}^{\alpha\beta})\pa_\alpha\pa_\beta \phi, \label{eq:curvwaveeqcommut1} \\
&\big|\Boxr_g Z \phi- \hat{Z}\Boxr_g \phi|\les \bigg (\frac {|Z
H|+|H|}{1+t+|q|} +\frac {|Z H|_{LL}+|H|_{L\cal T}}{1+|q|}\bigg )
\sum_{|I|\le 1}|\pa Z^I \phi| \label{eq:curvwaveeqcommutest1}
\end{align}
In general,
\begin{align}
&\Boxr_g Z^I \phi- \hat{Z}^I\Boxr_g
\phi= -\sum_{I_1+I_2=I, \, |I_2|<|I|}\hat H^{{I}_1\alpha\beta}
\pa_\alpha\pa_\beta Z^{I_2}\phi, \label{eq:curvwaveeqcommut2} \\
&|\Boxr_g Z^I \phi-\hat{Z}^I \Boxr_g \phi| \les \frac 1{1+t+|q|}
\,\,\,\sum_{|K|\leq |I|,}\,\, \sum_{|J|+(|K|-1)_+\le |I|} \,\,\,
|Z^{J} H|\,\, {|\pa Z^{K} \phi|}
\label{eq:curvwaveeqcommutest2} \\
+ \frac 1{1+|q|}
& \sum_{|K|\leq |I|}\Big(\sum_{|J|+(|K|-1)_+\leq |I|} \!\!\!\!\!|Z^{J} H|_{LL}
+\!\!\!\!\!\sum_{|J^{\prime}|+(|K|-1)_+\leq
|I|-1}\!\!\!\!\!|Z^{J^{\prime}} H|_{L\cal T}
+\!\!\!\!\!\sum_{|J^{\prime\prime}|+(|K|-1)_+\leq |I|-2}\!\!\!\!\!
|Z^{J^{\prime\prime}} H|\Big) {|\pa Z^{K} \phi|}\nn
\end{align}
where $(|K|-1)_+=|K|-1$ if $|K|\geq 1$ and $(|K|-1)_+=0$ if
$|K|=0$.
\end{cor}
\begin{proof}
First observe that
\begin{align*}
\hat Z \Boxr_g \phi &= (Z+ c_Z) \Box \phi +
(Z+c_Z)H^{\a\b}\pa^2_{\a\b} \phi
\\ &= \Box Z \phi + H^{\a\b} \pa^2_{\a\b} Z\phi +
(Z H^{\a\b} ) \pa^2_{\a\b} \phi + (H^{\a\b}_Z + c_Z
H^{\a\b})\pa^2_{\a\b}\phi \\ &= \Boxr_g Z\phi + (Z H^{\a\b} )
\pa^2_{\a\b} \phi + (H^{\a\b}_Z + c_Z H^{\a\b})\pa^2_{\a\b}\phi
\end{align*}
Recall now that the constant $c_Z$ is different from $0$ only in
the case of the scaling vector field $S$. Moreover, in that case
$$
H^{\a\b}_S + c_S H^{\a\b} =0
$$
The inequality \eqref{eq:curvwaveeqcommutest1} now follows from
\eqref{eq:curvwaveeqcommut1}, \eqref{eq:kLT} and the estimate
\eqref{eq:derframeZ}. The general commutation formula
\eqref{eq:curvwaveeqcommut2} follows from the following
calculation, similar to the one in Lemma \ref{commut2}. We have
$$
\hat Z^{I}\Boxr_{g}\phi = \hat Z^{I}\Box \phi + \hat Z^{I}
H^{\a\b}\pa^{2}_{\a\b}\phi = \Box Z^{I}\phi + \sum_{J+K=I} \hat
Z^{J} H^{\a\b} Z^{K}\pa^{2}_{\a\b}\phi
$$
If we now use \eqref{eq:Zcom2} we get \eqref{eq:curvwaveeqcommut2}
as in the proof of Lemma \ref{commut2}. The inequality
\eqref{eq:curvwaveeqcommutest2} now follows from
\eqref{eq:curvwaveeqcommut2}, \eqref{eq:kLT} and the estimate
\eqref{eq:derframeZ}.
\end{proof}

\section{Appendix B. Hardy type inequality}\label{section:poincare}
In this section we prove a version of the classical three dimensional
Hardy inequality
$$
\int_{\R^3} \frac {|f(x)|^2}{|x|^2}\, dx\le 4 \int_{\R^3} |\nab f(x)|^2 dx.
$$
The Hardy inequality converts the weighted $L^2$ norm of the function
$f$ into the $L^2$ norm of its gradient. These type of estimates prove to
be useful in the context of energy estimates for solutions of a quasilinear
wave equation $\Box_{g(\phi)}\phi=F$, where the energy norms contain
only the derivatives of $\phi$ while the error terms, generated by the
metric $g(\phi)$ also depend on the solution $\phi$ itself. The disadvantage
of the classical Hardy inequality in this context is that it requires a costly weight
$r^{-2}$. The "cost" here refers to the rate of decay of the weight in the wave zone
$r\approx t$. Therefore we seek a version of the Hardy inequality with the
weight dependent on the distance to the cone $r=t$ rather than the origin $r=0$.
\begin{lemma} Let $0\leq \alpha\leq 2$, $1+\mu>0$ and $\gamma>0$. Then
for any function $u\in C^1_0([0,\infty))$ and an arbitrary $t\ge 0$
there is a constant $C$,
depending on a lower bound for $\gamma>0$ and $1+\mu>0$,
such that
\begin{align}
\int_{0}^t \frac{u^2}{(1+|r-t|)^{2+\mu}}&
\frac{r^2\,dr }{(1+t+r)^{\alpha}}
+\int_{t}^\infty \frac{u^2}{(1+|r-t|)^{1-\gamma}}
\frac{r^2\,d r }{(1+t+r)^{\alpha}}\label{eq:Hardy}\\
&\leq C\int_{0}^t \frac {|\pa_r u|^2}{(1+|r-t|)^{\mu}}\frac{r^2}{(1+t+r)^{\alpha}}\,d r
+C\int_{t}^\infty {|\pa_r u|^2}\frac{(1+|r-t|)^{1+\gamma}}{(1+t+r)^{\alpha}}r^2\, d r\, \nn
\end{align}
\end{lemma}
\begin{remark}
The inequality \eqref{eq:Hardy} with $\alpha=0$ appears to be the
precise analogue of the classical (spherically-symmetric) Hardy
inequality. The presence of the additional growing weight
$(1+|t-r|)^{1+\gamma}$ seems to be necessary and in fact fits
perfectly in the context of our weighted energy estimates. In the
previous works \cite{L-R2} the estimates used to convert a weighted
norm of a solution into a norm of its derivative, e.g.,
$$
\int_{0}^\infty \frac{u^2}{(1+|r-t|)^{2+\mu}} r^2 dr \leq C \Big ( (1+t)^2| u(1+t)|^2
+ \int_{0}^\infty \frac{|\pa_ru|^2}{(1+|r-t|)^{2+\mu}} r^2 dr\Big )
$$
 were more reminiscent of the classical
Poincar\'e inequality.
\end{remark}
\begin{proof} Set
$$
m(q)=\begin{cases} (1+q)^{\gamma},\quad q>0\\
             (1-q)^{-1-\mu},\quad q\le 0\end{cases}
$$
Then
$$
m^\prime(q)=\begin{cases} \gamma(1+q)^{-1+\gamma},\quad q\ge 0\\
           (1+\mu)  (1-q)^{-2-\mu},\quad q\le 0\end{cases}
$$
Since $\alpha\leq 2$
$$
\pa_r \Big( r^2(1+t+r)^{-\alpha} m(r-t)\Big)=
\Big(\frac{2}{r}-\frac{\alpha}{1+t+r}+\frac{m^\prime(r-t)}{m(r-t)}\Big)
 \frac{r^2 m(r-t)}{(1+t+r)^{\alpha}}\geq
\frac{r^2 \,m^\prime(r-t)}{(1+t+r)^{\alpha}}
$$
Hence
$$
\pa_r\Big( r^2(1+t+r)^{-\alpha} m(r-t)\phi^2\Big)\geq
m^\prime(r-t)r^2(1+t+r)^{-\alpha}\phi^2+ 2r^2(1+t+r)^{-\alpha}
m(r-t)\phi \pa_r\phi
$$
If we integrate the above from $0$ to $\infty$ and use that
$\phi$ has compact support we see that
$$
\int_0^\infty m^\prime(r-t)\phi^2 \, \frac{r^2 dr}{(1+t+r)^{\alpha}
} \leq 2\int_0^\infty \, m(r-t) \phi \pa_r\phi\frac{r^2
dr}{(1+t+r)^{\alpha} }
$$
Since $m^\prime>0$ and $m\geq 0$ it follows from Cauchy-Schwarz inequality that
$$
\int_0^\infty m^\prime(r-t)\phi^2 \, \frac{r^2 dr}{(1+t+r)^{\alpha}
} \leq \sqrt{2}\int_0^\infty \, \frac{m(r-t)^2}{m^\prime(r-t)}
 |\pa_r\phi|^2\frac{r^2 dr}{(1+t+r)^{\alpha} }
$$
from which the lemma follows.
\end{proof}

\begin{cor} \label{cor:Poinc}
Let $\gamma>0$ and $\mu>0$ and set, for $q=r-t$,
$$
w(q)=\begin{cases} 1+(1+|q|)^{1+2\gamma},\quad q>0,\\
                   1+(1+|q|)^{-2\mu},\quad q<0
                   \end{cases}
$$
Then for any $-1\leq a\leq 1$ and any $\phi\in C^\infty_0(\R^3)$
$$\label{eq:poincareone}
\int \frac{|\phi|^2}{(1+|q|)^2}\,\frac{ w\, dx}{(1+t+|q|)^{1-a}}
\les \int |\pa \phi|^2\, \frac{ w\, dx}{(1+t+|q|)^{1-a}}
$$
If in addition  $a<2\min{(\gamma, \mu)}$,
$$\label{eq:poincaretwo}
\int \frac{|\phi|^2}{(1+|q|)^2}\,\frac{(1+|q|)^{-a} \,
}{(1+t+|q|)^{1-a}}\, \frac {w\, dx}{(1+q_-)^{2\mu}}\les \int {|\pa
\phi|^2}\min\, (w', \frac{w}{(1+t+|q|)^{1-a}}) \, dx
$$
where $q_-=|q|$, when $q<0$ and $q_-=0$, when $q>0$.
\end{cor}

\section{Appendix C. Weighted Klainerman-Sobolev inequalities}\label{eq:globalsobolev}
In this section we provide a straightforward generalization of the
Klainerman-Sobolev inequalities, expressing pointwise decay in terms
of the bounds on $L^2$ norms involving vector field $Z\in {\cal Z}$.
We consider energy norms with the following weight function
$$
w=w(q)=\begin{cases} 1+(1+|q|)^{1+2\gamma},\quad\text{when }\quad q>0\\
         1+(1+|q|)^{-2\mu}\,\quad\text{when }\quad q<0\end{cases}
$$
for some $0<\gamma<1$. Therefore,
\begin{align}
&w^{\,\prime}=:w^{\,\prime}(q)
=\begin{cases} (1+2\gamma)(1+|q|)^{2\gamma},\quad\text{when }\quad q>0\\
         2\mu(1+|q|)^{-1-2\mu},\quad\text{when }\quad q<0\end{cases}\nn\\
&w^\prime\leq 4w(1+|q|)^{-1}\leq 16\gamma^{-1}
w^\prime(1+q_-)^{2\mu}.\nn
\end{align}
We have the following global Sobolev inequality
\begin{prop} \label{prop:K-S}
For any function $\phi\in C^\infty_0(\R^3)$ and an arbitrary
$(t,x)$,
$$
|\phi(t,x)|(1+t+|q|)\big[ (1+|q|)w(q)\big]^{1/2} \leq
C\sum_{|I|\leq 3} \|w^{1/2} Z^I \phi(t,\cdot)\|_{L^2},\qquad q=t-r
$$
\end{prop}
\begin{proof} Fist note that it is sufficient to consider two cases when
the support of $\phi$ is in the sets $r\leq t/2$ and
$t/4\leq r$ respectively.  We argue as follows.
 Let $\chi(\tau)$ be a smooth cut-off function such
$\chi(\tau)=1$ when $\tau\leq -3/5$, $and \chi(\tau)=0$ for $\tau\geq -1/3$.
Define
$$
\phi_1(t,x) =\chi(\frac {r-t}{r+t}) \phi(t,x),\qquad
\phi_2(t,x) =(1-\chi(\frac {r-t}{r+t}) )\phi(t,x)
$$
The supports of functions $\phi_1, \phi_2$ then belong to the
desired regions
and $\sum_{|I|\leq 3}|Z^I\phi_i|\leq C\sum_{|I|\leq 3}|Z^I\phi|$,
since $|Z^I \chi_i|\leq C$.

First for $r\leq t/2$ consider the rescaled function
c$\phi_t(x)=\phi(t,tx)$, where $|x|\leq 1/2$ in the support.
By the standard Klainerman-Sobolev estimate
$$
\|\phi_t(x)\|_{L_x^\infty}
\leq C\sum_{|\alpha|\leq 2} \|\pa^\alpha \phi_t(x)\|_{L_x^2}
\leq C\sum_{|\alpha|\leq 2} \|(t^{|\alpha|} \pa^\alpha\phi)(t,t x)\|_{L_x^2}
$$
Since
$$
\sum_{|\alpha|\leq k} |t-r|^{|\alpha|} |\pa^\alpha\phi(t,x)|
\leq C\sum_{|I|\leq k} |Z^I\phi(t,x)|,
$$
we obtain  that
$$
\|\phi(t,x)\|_{L_x^\infty}\leq C t^{-3/2}
\sum_{|I|\leq k} \|Z^I\phi(t,x)\|_{L_x^2}.
$$
Multiplying both sides of the above inequality by $[w(-t)]^{1/2}$, which is
$\approx [w(q)]^{1/2}$ when
$r\leq t/2$, proves the proposition in the case when $\phi$ is supported in the region
$r\le t/2$.

When $r\geq t/4$ we can estimate the $L^\infty$
norm by the $L^1$ norms of derivatives:
$$
|\phi(t,x)|^2(1+t+|q|)^2 w(q)(1+|q|)
\leq C\sum_{|\alpha|\leq 2}\int_{S^2} \int\Big|\pa_\omega^\alpha
\pa_q \Big( w(q) (1+|q|)
(t+q)^2 \phi\big(t,(t+q)\omega\big)^2\Big)\Big|\, dq \,dS(\omega)
$$
Since $|w^\prime(q)|(1+|q|)\leq C w(q)$ and
$|q|\leq C|t+q|=Cr$ on the support of $\phi$  it follows that
we have the bound
\begin{align*}
\sum_{|\alpha|\leq 2, k=0,1}\int_{S^2} \int
w(q)
(t+q)^2  \Big|(q\pa_q)^k \pa_\omega^\alpha \phi\big(t,(t+q)\omega\big)\Big|^2\, dq \,dS(\omega)\leq C\sum_{|I|\leq 3} \int_{\bold{R}^3} w(q) |Z^I \phi(t,x)|^2\, dx
\end{align*}
\end{proof}

\section{Appendix D. Wave coordinate condition}

In this section we provide the details on the estimates following from the
wave coordinate condition
\beq\label{eq:WTR}
\pa_\mu \Big (g^{\mu\nu} \sqrt {|\det g|}\Big )=0
\eq
for a Lorentzian metric $g$ in a coordinate system $\{x^\mu\}_{\mu=0,...,3}$,
leading up to the proof of Proposition \ref{decaywavecZ}.
We recall the definition of the tensor $H^{\mu\nu} = g^{\mu\nu} -m^{\mu\nu}$
and below state the consequences of \eqref{eq:WTR} in terms of estimates for $H$.

We first observe that for any vector field $X$ and a collection of
our special vector fields $Z\in {\cal Z}$ we have that
\begin{equation}\label{eq:divergencecommutators}
Z^I \pa_\alpha X^\alpha
=\pa_\alpha\Big(Z^I X^\alpha +\sum_{|J|<|I|}
c_{J\,\,\gamma}^{\,I\,\,\,\alpha} Z^{J} X^{\gamma}\Big)
=\pa_\alpha\Big( \sum_{|J\le |I|} c_{J\,\gamma}^{\,I\,\,\,\alpha}
Z^{J} X^{\gamma}\Big),
\end{equation}
where $c_{J\,\gamma}^{\,\,\,\alpha}$ are constants such that
$$
 c_{J\, \gamma}^{\,I\,\,\,\alpha}=\delta^{\alpha}_{\,\, \gamma},
\quad\text{for}\quad |J|=|I|\qquad \text{and}\qquad c_{J\,
L}^{\,I\,\,\,\Lb}=0,\quad\text{for}\quad |J|=|I|-1
$$
The last identity is a consequence of the relation between $c_{J\,
\a}^{\,I\,\,\,\ga}$ and the commutator constants
$c_{\a\b}=[\pa_{\a}, Z]_{\b}$ for which we have established in Lemma
\ref{commut1} that
$c_{LL}=0$. It therefore follows from \eqref{eq:WTR} and
\eqref{eq:divergencecommutators} that
\beq \label{eq:defHJ}
H_{\mu\nu}^{[I]}:=Z^I \widetilde{H}_{\mu\nu}+\!\!\!\!\sum_{|J|<|I|}
c_{J\, \mu}^{I\,\,\,\,\ga } Z^{J} \widetilde{H}_{\gamma \nu},
\qquad\text{with} \quad
\widetilde{H}_{\mu\nu}:=H_{\mu\nu}-\frac{m_{\mu\nu}}{2}\tr H,
 \eq
 satisfies
 \beq
\pa_\mu {H^{[I]}}^{\mu\nu}+Z^I\pa_\mu O^{\mu\nu}(H^2)=0 , \qquad
\text{where}\quad O^{\mu\nu}(H^2)=O(|H|^2).
 \eq
 \begin{lemma}\label{HJest}
  Let $H$ be a 2-tensor and let $H^{[I]}$ be defined by
 \eqref{eq:defHJ}. Then, for $j=0,1$;
 \beq\label{eq:wavecoordianteestimmate}
 \sum_{|I|\leq k}|\pa^j Z^I H|_{LL}
 +\sum_{|J|\leq k-1}|\pa^j Z^J H|_{L\cal T} \les
 \sum_{|K|\leq k-2}|\pa^j Z^K H|
 +\sum_{|I|\leq  k} |\pa^j H^{[I]}|_{L\cal T}
 \eq
 \end{lemma}
 \begin{proof}
It is easy to see that \eqref{eq:defHJ} implies that for any $\alpha$
 \begin{align*}
 &\sum_{|I|\leq k}|\pa^\alpha Z^I H|_{LL}\les
 \sum_{|J|\leq k-1}|\pa^\alpha Z^J H|_{L\cal T}
  +\sum_{|K|\leq k-2}|\pa^\alpha Z^K H|
 +\sum_{|I|\leq  k} |\pa^\alpha H^{[I]}|_{LL},\\
& \sum_{|J|\leq k-1}|\pa^\alpha Z^J H|_{L\cal T}\leq
 \sum_{|K|\leq k-2}|\pa^\alpha Z^K H|
 +\sum_{|J|\leq  k-1} |\pa^\alpha H^{[J]}|_{L\cal T}
 \end{align*}
 and \eqref{eq:wavecoordianteestimmate} follows.
 \end{proof}

\begin{lemma}\label{decaywavectwo}
Let $Z$ represent Minkowski Killing or conformally Killing vector
fields from our family {\cal Z} and assume the tensor $H$ satisfies
the wave coordinate condition \eqref{eq:WTR}. Then for any
multi-index $I$ tensor $H^{[I]}$, defined in \eqref{eq:defHJ}
satisfies the following estimate
 \beq \label{eq:ZIHLL}
 \big| \pa
H^{[I]}\big|_{L\cal T} \les\sum_{|J|\leq |I|} |\pab Z^J H|+
\sum_{I_{1}+...+I_k=I,\, k\geq 2} |Z^{I_k} H|\cdot\cdot\cdot|Z^{I_2}
H| \, |\pa Z^{I_1} H| \nn .
 \eq
\end{lemma}
\begin{proof}
The wave coordinate condition \eqref{eq:wavecdef} can be written
in the form
$$
\pa_\mu \big( \widetilde{G}^{\mu\nu} \big)=0,\qquad
\text{where}\quad
\widetilde{G}^{\mu\nu}=(m^{\mu\nu}+H^{\mu\nu})\sqrt{|\det g|}.
$$
It follows from \eqref{eq:divergencecommutators}  that
$$
\pa_\mu \big(\sum_{|J|\le |I|} c_{J}^{\,I\,\,\,\mu\gamma}
Z^{J}\widetilde{G}_{\gamma\nu}\big)=0.
$$
Decomposing relative to the null frame $(L,\Lb, A, B)$ we obtain
$$
\pa_q \big( \sum_{|J|\le |I|} c_{J}^{\,I\,\,\,\Lb\,\gamma }
Z^{J}\widetilde{G}_{\gamma \nu} \big) =\pa_s  \Big(\sum_{|J|\le
|I|} c_{J}^{\,I\,\,\, L \gamma}
 Z^{J}\widetilde{G}_{\gamma\nu} \Big)
- A_{\mu} \pab_A \big( \sum_{|J|\le |I|} c_{J}^{\,I\,\,\,\mu\gamma
} Z^{J}\widetilde{G}_{\,\ga \nu} \big).
$$
We now contract the above identity with one of the tangential
vector fields $T^{\nu}$, $T\in \{ L,A,B\}$ to obtain
$$
\Big|L^\gamma \, T^\nu \pa_q  Z^I \widetilde{G}_{\gamma\nu}
 +\!\!\!\!\sum_{|J|<|I|}
c_{J}^{\,I\,\,\,\Lb\ga}  T^\nu\, \pa_q Z^{J}\widetilde{G}_{\gamma
\nu} \Big| \les \sum_{|J|\leq |I|} \big| \pab Z^{I}\tilde{G}\big|
$$
We now examine the expression
$$
L^\gamma \, T^\nu Z^J  \pa_q \widetilde{G}_{\gamma\nu}= L^\gamma
\, T^\nu \pa_q Z^J \bigg (
(m_{\gamma\nu}+H_{\gamma\nu})\sqrt{|\det g|} \bigg
)=\sum_{J_1+J_2=J} L^\gamma \, T^\nu \pa_q \bigg ( (Z^{J_1}
H_{\gamma\nu}) Z^{J_2} \sqrt{|\det g|}\bigg )
$$
since $m_{LT}= L^\gamma \, T^\nu m_{\gamma\nu}=0$. The desired
estimate now follows from the identity $\sqrt{|\det g|}=1+f(H)$,
which holds with a smooth function $f(H)$  such that $f(0)=0$ and
$f(H)=-\tr H/2+ O(H^2)$
\end{proof}
We now summarize the above results in the following
\begin{lemma}\label{lem:decaywavecZ}
Let $g$ be a Lorentzian metric satisfying the wave coordinate
condition \eqref{eq:WTR} relative to a coordinate system
$\{x^\mu\}_{\mu=0,...,3}$. Then the following estimates for the
tensor  $H^{\mu\nu}= g^{\mu\nu} -m^{\mu\nu}$ and a multi-index $I$
hold true under the assumption that $|Z^J H|\le C,\,\text{for all }
|J|\le |I|/2$:
\begin{align}
&|\pa Z^I H|_{ L\cal T}\les\Big (\sum_{|J|\leq |I|}|\overline{\pa} Z^J
H|+\!\!\!\!\sum_{|J|\leq |I|-1}\!\!\!|\pa Z^J H|\,\, +
 \sum_{\,\,\,\,|I_1|+||I_2| \leq |I|}|Z^{I_{2}} H| |\pa Z^{I_1} H|\Big )\label{eq:decaywavec6}\\
& |\pa Z^I H|_{ LL}\les  \Big (\sum_{|J|\leq |I|}|\overline{\pa} Z^J
H|+ \sum_{|J|\leq |I|-2} |\pa Z^{J} H| + \sum_{|I_1|+|I_2| \leq |I|}
|Z^{I_{2}} H| |\pa Z^{I_1} H|\Big )\label{eq:decaywavec5}.
\end{align}
\end{lemma}

We now establish that the tensor \beq\label{eq:def-H1}
H_1^{\mu\nu}=H^{\mu\nu} -H_0^{\mu\nu}, \qquad
H_0^{\mu\nu}:=-\chi(\frac rt) \chi (r) \frac Mr \delta^{\mu\nu} \eq
 obtained by subtracting the "Schwarzschild part" $H_0$ from $H$
 obeys similar structure estimates to those of $H$.
 \begin{lemma}\label{lem:Wave-H1}
 Let $g$ be a Lorentzian metric satisfying the wave coordinate condition
\eqref{eq:WTR} relative to a coordinate system $\{x^\mu\}_{\mu=0,...,3}$.
For a given integer $k\ge 0$
we have the following estimates for the tensor  $H_1^{\mu\nu}= g^{\mu\nu} -m^{\mu\nu}- H_0^{\mu\nu}$, under the assumption that
$|Z^J H|\le C,\,\forall |J|\le k/2$,
 \begin{multline}\label{eq:HJ1est2}
\sum_{|I|\le k}
\big| \pa Z^I H_1\big|_{\cal LL} + \sum_{|I|\le k-1}
\big| \pa Z^I H_1\big|_{\cal LT} \les\sum_{|I|\leq k} |\pab Z^{I}
H_1|+\sum_{|I|\le k-2} |\pa Z^I H_1|\\+ \ve\frac 1{1+t+|q|}\sum_{|I|\leq k} \Big( |\pa Z^{I}
H_1|+\frac{|Z^{I} H_1|}{1+t+|q|}\Big) \\+
\sum_{|I|+|J|\le k} |Z^I H_1| |\pa Z^J H_1|+
\frac{C\ve\chi_0(1/2<r/t<3/4)}{(1+t+|q|)^2}+\frac{C\ve^2}{(1+|t|+|q|)^3},
\end{multline}
where $\chi_0(1/2<r/t<3/4)$ is the characteristic function of the
set where $t/2<r<3t/4$.
\end{lemma}
\begin{proof}
We define the tensors
\begin{align*}
&H_{\mu\nu}^{[I]}=Z^I
\widetilde{H}_{\mu\nu}+\!\!\!\!\sum_{|J|<|I|} c_{J\,
\mu}^{I\,\,\,\,\ga } Z^{J} \widetilde{H}_{\gamma \nu},
\qquad
\widetilde{H}_{\mu\nu}=H_{\mu\nu}-\frac{m_{\mu\nu}}{2}\tr H,\\
&{H_0}_{\mu\nu}^{[I]}=Z^I
{\widetilde{H_0}}_{\mu\nu}+\!\!\!\!\sum_{|J|<|I|} c_{J\,
\mu}^{I\,\,\,\,\ga } Z^{J} \widetilde{H_0}_{\gamma \nu},
\qquad
{\widetilde{H_0}}_{\mu\nu}={H_0}_{\mu\nu}-\frac{m_{\mu\nu}}{2}\tr H_0,\\
&{H_1}_{\mu\nu}^{[I]}=Z^I
{\widetilde{H_1}}_{\mu\nu}+\!\!\!\!\sum_{|J|<|I|} c_{J\,
\mu}^{I\,\,\,\,\ga } Z^{J} {\widetilde{H_1}}_{\gamma \nu},
\qquad
{\widetilde{H_1}}_{\mu\nu}={H_1}_{\mu\nu}-\frac{m_{\mu\nu}}{2}\tr H_1.
\end{align*}
with the constants $c_{J\,
\mu}^{I\,\,\,\,\ga } $ such that
$c_{J\,L}^{I\,\,\,\,\underline{L} }=0$, if $|J|=|I|-1$.
Observe that the tensor ${H_0}^{[I]}_{\mu\nu}$ is defined in such a way that
$$
\pa^\mu {H_0}^{[I]}_{\mu\nu}= Z^I \Big (\pa^\mu {\widetilde {H_0}}_{\mu\nu}\Big )
$$
Since the tensor $g^{\mu\nu}$ satisfies the wave coordinate condition it follows that
 \beq
\pa^\mu{H^{[I]}_{\mu\nu}}+Z^I\pa_\mu O^{\mu\nu}(H^2)=0 , \qquad
\text{where}\quad O^{\mu\nu}(H^2)=O(|H|^2).
 \eq
On the other hand, calculating using the definition of $H_0$,
$$
\pa^\mu
{\widetilde{H_0}}_{\mu\nu}=2\chi^{\prime}(r/t)\chi(r){M}/{t^2}\delta_{\nu\,0}
$$
Therefore, since $H_{\mu\nu} = {H_0}_{\mu\nu} + {H_1}_{\mu\nu}$ and
$H_{\mu\nu}^{[I]} = {H_0}^{[I]}_{\mu\nu} + {H_1}^{[I]}_{\mu\nu}$
we obtain
$$
\pa^\mu {H_1}^{[I]}_{\mu\nu}=-Z^I \pa_\mu O^{\mu\nu}(H^2)+2Z^I
\big(\chi^{\prime}(r/t)\chi(r){M}/{t^2}\big)\delta_{\nu \,0}
$$
Arguing as in the proof of Lemma \ref{decaywavectwo} we derive
 \begin{align}
\sum_{|I|\le k} \big| \pa Z^I H_1\big|_{\cal LL} &+ \sum_{|I|\le k-1}
\big| \pa Z^I H_1\big|_{\cal LT} \les \sum_{|I|\le k}
\big| \pa H_1^{[I]}\big|_{L\cal T} + \sum_{|I|\le k-2} |\pa Z^I H_1|\\ &\les\sum_{|I|\leq k} |\pab Z^{I}
H_1|+ \sum_{|I|\le k}\sum_{I_{1}+...+I_n=I,\, n\geq 2} |Z^{I_n}
H|\cdot\cdot\cdot|Z^{I_2}  H| \, |\pa Z^{I_1} H|\\
&+\sum_{|I|\le k-2} |\pa Z^I H_1|+\frac{C\ve\chi_0(1/2<r/t<3/4)}{(1+t+|q|)^2}, \label{eq:ZHI}
\end{align}
where $\chi_0(1/2<r/t<3/4)$ is the characteristic function of  the
set where $t/2<r<3t/4$.  Using the condition $|Z^I H|\les 1,\,
\forall |I|\le k/2$ and the estimates
$$
|Z^I H_0|+ (1+t+|q|) |\pa Z^I H_0|\le \frac {\ve}{1+t+|q|}
$$
we obtain
\begin{multline}
\sum_{I_{1}+...+I_n=I,\, n\geq 2} |Z^{I_n}
H|\cdot\cdot\cdot|Z^{I_2}  H| \, |\pa Z^{I_1} H|\les
\sum_{|J|+|K|\leq |I| } |Z^J H|\,|\pa Z^K H|\\
\les \sum_{|J|+|K|\leq |I| } |Z^J H_0|\,|\pa Z^K H_1| +
 |Z^J H_1|\,|\pa Z^K H_0|+ |Z^J H_1|\,|\pa Z^K H_1|+
 |Z^J H_0|\,|\pa Z^K H_0|\\
 \les \sum_{|J|\leq |I| } \ve\frac{1}{1+t+|q|}
 \Big( |\pa Z^J H_1| + \frac{ |Z^J H_1|}{1+t+|q|}\Big)
 + \frac{\ve^2}{(1+t+|q|)^{3}}
\end{multline}
and  the lemma follows.
\end{proof}

\section{Appendix E: $L^1-L^\infty$ estimates and additional decay}
The bounds
\begin{equation}
\label{eq:enerNdef} \E_N(t)=\sum_{|I|\leq N} \|w^{1/2}\pa Z^I
h^1(t,\cdot)\|_{L^2} + \|w^{1/2}\pa Z^I \psi(t,\cdot)\|_{L^2} \leq
C_N (1+t)^{C_N\varepsilon},
\end{equation}
with \beq\label{eq:weight}
 w^{1/2}=\begin{cases} (1+|r-t|)^{1/2+\gamma},\,\,\, &r>t\\
1,\quad &r\leq t\end{cases} \eq imply, with the help of the
weighted Klainerman-Sobolev inequality, that
\begin{equation}\label{eq:globaldecay2}
| \pa Z^I h^1| + |\pa Z^I \psi|\leq \begin{cases}
C_N^\prime\varepsilon
(1+t+r)^{-1+C_N\varepsilon}(1+|t-r|)^{-1-\gamma},\quad r>t\\
C_N^\prime \varepsilon
(1+t+r)^{-1+C_N\varepsilon}(1+|t-r|)^{-1/2},\quad r\le t
\end{cases},\qquad  |I|\leq N-2
\end{equation}
 The exponent in the interior decay estimate \eqref{eq:globaldecay2}
 in $|t-r|$ can be improved from
$-1/2$ to $-1$ using the estimate
\begin{equation}\label{eq:globaldecay3}
| Z^I h^1| + |Z^I \psi|\leq C_N^{\prime\prime} \varepsilon
(1+t)^{-1+2C_N\varepsilon},\qquad |I|\leq N-3
\end{equation}
that we will now prove. We will divide the solution of Einstein's
equations into a linear and a nonlinear part with vanishing initial
data $h^1_{\mu\nu}=v_{\mu\nu}+w_{\mu\nu}$. The estimate for the
linear part follows from:
\begin{lemma} If $v$ is the solution of
 $$
\sq v=0,\quad \quad v\big|_{t=0}=v_0,\quad \pa_t v\big|_{t=0}=v_1
 $$
 then for any $\gamma>0$;
 \beq
 (1+t)|v(t,x)| \leq C\sup_x \big( (1+|x|)^{2+\gamma}
 ( |v_1(x)|+|\partial v_0(x)|) +(1+ |x|)^{1+\gamma}|v_0(x)|\big)
 \eq
\end{lemma}
\begin{proof} The proof is an immediate consequence of
the Kirchhoff's formula
 $$
v(t,x)=t\int_{|{\omega}|=1}{\big(v_1(x+t{\omega})+\langle
v^{\prime}_0(x+t{\omega}),{\omega}\rangle \big)\,dS({\omega})} +
 \int_{|{\omega}|=1}{v_0(x+t{\omega})\,dS({\omega})},
 $$
where $dS({\omega})$ is the normalized surface measure on $S^2$.
 Suppose that $x=r\bold{e}_1$, where $\bold{e}_1=(1,0,0)$.
 Then for $k=1,2$ we must estimate
 \beq
 \int\!\! \frac{dS(\omega)}{1+|r\bold{e}_1+t\omega|^{k+\gamma}}
 =\int_{-1}^1\frac{C d\omega_1}{1+\big(
 (r-t\omega_1)^2+t^2(1-\omega_1^2)\big)^{(k+\gamma)/2}}
 \leq \int_0^2\!\! \frac{C d s}{1+\big(
 (r-t+t s)^2+t^2 s\big)^{(k+\gamma)/2}}
 \eq
 If $k=2$ we make the change of variables $t^2 s=\tau$ to get
 an integral bounded by $Ct^{-2}$ and if $k=1$, we make the change
 of variables $t s=\tau$ to get an integral bounded by $t^{-1}$.
 \end{proof}
To estimate the nonlinear part we use H\"ormander's $L^1-L^\infty$
estimates for the fundamental solution of $\Box$, see \cite{H1, L1}:
 \begin{prop}\label{prop:hormander} If $v$
be the solution of
 $$
\sq w=g,\quad \quad w\big|_{t=0}=\pa_t w\big|_{t=0}=0
 $$
then
 \beq |w(t,x)|(1+t+|x|) \leq
C\sum_{|I|\leq
2}\int_0^t\int_{\bold{R}^3}\frac{|Z^Ig(s,y)|}{1+s+|y|}\, dy\,ds,
 \eq
\end{prop}

Let $h^1_{\mu\nu}=v_{\mu\ nu}+w_{\mu\nu}$ where, 
 \beq
 \Box w_{\mu\nu}= -H^{\alpha\beta}\pa_\alpha\pa_\beta h_{\mu\nu}
 +F_{\mu\nu}(h)(\pa h,\pa h),\qquad w_{\mu\nu}\big|_{t=0}=\pa_t
 w_{\mu\nu}\big|_{t=0}=0,
 \eq
 and 
 \beq
 \Box v_{\mu\nu}=0,\qquad v_{\mu\nu}\big|_{t=0}=h^1_{\mu\nu}\big|_{t=0},\qquad \pa_t
 v_{\mu\nu}\big|_{t=0}=\pa_t
 h^1_{\mu\nu}\big|_{t=0}.
 \eq
We have
 \beq
 |Z^I F_{\mu\nu}(h)(\pa h,\pa h)|\leq
 C\!\!\!\sum_{|J|+|K|\leq |I|\!\!\!\!\!\!}|\pa Z^J h|\, |\pa Z^K h|
 +C\!\!\! \sum_{|J|+|K|\leq |I|} \frac{|Z^J h|}{1+|q|}\, |\pa Z^K h|, 
 \eq
 and since $H^{\alpha\beta}=-h_{\alpha\beta}+O(h^2)$, 
  \beq
\big|Z^I\big( H^{\alpha\beta}\pa_\alpha\pa_\beta h_{\mu\nu}\big)\big|\leq C\!\!\!
\sum_{|J|+|K|\leq |I|+1,\, |J|\leq |I|} \frac{|Z^J h|}{1+|q|}\, |\pa
Z^K h|. 
  \eq
  Now
  \beq
\int |\pa Z^J h|\, |\pa Z^K h|(s,y)\, dy \leq \sum_{|I|\leq N}\|\pa
Z^I h(s,\cdot)\|_{L^2}^2\leq C(1+s)^{2C_N\varepsilon}.
  \eq
 We write $h=h^0+h^1$ and estimate
 \beq
 \int \frac{|h^0(s,y)|^2}{(1+|q|)^2}\,dy \leq M^2 \int_0^\infty
 \frac{r^2\, dr}{(1+|t+r|)^2(1+|t-r|)^2}\leq CM^2
 \eq
and by Corollary \ref{cor:Poinc}
 \beq
 \int \frac{|h^1(s,y)|^2}{(1+|q|)^2}\,dy
 \leq C\int |\pa h^1(s,y)|^2 w(q)\, dy\leq C_N^2\varepsilon^2 (1+t)^{2C_N\varepsilon}
 \eq
 where $w$ is as in \eqref{eq:weight}. Hence
 \beq
 \int \frac{|Z^J h|}{1+|q|}\, |\pa
Z^K h|(s,y)\, dy \leq C\varepsilon^2 (1+t)^{2C_N\varepsilon}
 \eq
 It now follows from Proposition \ref{prop:hormander} that
 \beq
 |w_{\mu\nu}(t,x)|(1+t+|x|)\leq \int_0^t
 \frac{\varepsilon^2 \, ds}{(1+s)^{1-2C_N\varepsilon}}\leq
 C\varepsilon (1+s)^{2C_N\varepsilon}
 \eq
 which proves \eqref{eq:globaldecay3}.

\end{document}